\documentclass[11pt]{article}
\usepackage{mathrsfs}
\usepackage{amsthm}
\usepackage{amssymb}
\usepackage{amsmath}
\usepackage{graphicx}
\usepackage{color}
\usepackage{amsfonts}
\usepackage{float}
\usepackage{cite}
\usepackage [latin1]{inputenc}
\usepackage[text={140mm,210mm},left=35mm,vmarginratio=1:1]{geometry}
\newtheorem{theorem}{Theorem}[section]

\newtheorem{lemma}[theorem]{Lemma}
\newtheorem{proposition}[theorem]{Proposition}

\numberwithin{equation}{section}
\normalsize

\begin{document}
\title{\textbf{Large and moderate deviations for empirical density fields of stochastic SEIR epidemics with vertex-dependent transition rates}}

\author{Xiaofeng Xue \thanks{\textbf{E-mail}: xfxue@bjtu.edu.cn \textbf{Address}: School of Mathematics and Statistics, Beijing Jiaotong University, Beijing 100044, China.}\\ Beijing Jiaotong University\\Xueting Yin \thanks{\textbf{E-mail}: 22121655@bjtu.edu.cn \textbf{Address}: School of Mathematics and Statistics, Beijing Jiaotong University, Beijing 100044, China.}\\ Beijing Jiaotong University}

\date{}
\maketitle

\noindent {\bf Abstract:} In this paper, we are concerned with stochastic susceptible-exposed-infected-removed epidemics on  complete graphs with vertex-dependent transition rates. Large and moderate deviations of empirical density fields of our models are given. Proofs of our main results utilize exponential martingale strategies. Mathematical difficulties are mainly in checks of exponential tightness of fluctuation density fields of our processes. As an application of our main results, moderate deviations of a family of hitting times of our processes are also given.

\quad


\section{Introduction}\label{section one}
In this paper we are concerned with large and moderate deviations of stochastic susceptible-exposed-infected-removed (SEIR) epidemics with vertex-dependent transition rates. We first introduce the definition of our model. Our stochastic SEIR model $\{\xi_t\}_{t\geq 0}$ is a continuous-time Markov processes with state space
\[
\mathbb{X}=\left\{0,1,2,3\right\}^N,
\]
where $N\geq 1$ is an integer. That is to say, for any configuration $\xi\in \mathbb{X}$, $\xi=\left(\xi(0),\xi(1),\ldots, \xi(N-1)\right)$, where $\xi(i)\in \{0,1,2,3\}$ is the $i$th coordinate of $\xi$ for $0\leq i\leq N-1$. Throughout this paper we denote by $\mathbb{T}$ the $1$-dimensional torus $[0, 1)$. Our process $\{\xi_t\}_{t\geq 0}$ evolves as follows. For each $0\leq i\leq N-1$ and $t\geq 0$,
\begin{align*}
&\xi_t(i)\text{~jumps from $k$ to $l$ at rate~}\\
&\begin{cases}
\phi\left(\frac{i}{N}\right) & \text{~if~}k=2\text{~and~}l=3,\\
\psi\left(\frac{i}{N}\right) & \text{~if~}k=1\text{~and~}l=2,\\
\frac{1}{N}\sum_{j\neq i}\lambda\left(\frac{i}{N}, \frac{j}{N}\right)1_{\{\xi_t(j)=2\}} & \text{~if~}k=0\text{~and~}l=1,\\
0 &\text{~else,}
\end{cases}
\end{align*}
where $1_A$ is the indicator function of the event $A$ and $\psi,\phi\in C^\infty(\mathbb{T}), \lambda\in C^\infty(\mathbb{T}^2)$ are strictly positive.

The process $\{\xi_t\}_{t\geq 0}$ describes the spread of a SEIR epidemic on a complete graph with $N$ vertices symbolled as $0,1,\ldots,N-1$. Each vertex is in one of four states $0,1,2,3$. Vertices in state $0$ are susceptible. Vertices in state $1$ are exposed. Vertices in state $2$ are infected and those in state $3$ are removed. The transition-rate functions of our process are vertex-dependent. In detail, if vertex $i$ is infected, then it becomes removed at rate $\phi\left(\frac{i}{N}\right)$. If $i$ is exposed, then it becomes infected at rate $\psi\left(\frac{i}{N}\right)$. If $i$ is susceptible and $j$ is infected, then $i$ is infected by $j$ to become an exposed one at rate $\frac{1}{N}\lambda\left(\frac{i}{N},\frac{j}{N}\right)$. A removed vertex is frozen in this state and will never be infected again.

We can define $\{\xi_t\}_{t\geq 0}$ equivalently via its generator. According to transition-rate functions given above, the generator $\mathcal{L}$ of $\{\xi_t\}_{t\geq 0}$ is given by
\begin{align}\label{equ generator}
\mathcal{L}f(\xi)=&\sum_{i=0}^{N-1}\phi\left(\frac{i}{N}\right)1_{\{\xi(i)=2\}}\left(f(\xi^{i,3})-f(\xi)\right)\notag\\
&+\sum_{i=0}^{N-1}\psi\left(\frac{i}{N}\right)1_{\{\xi(i)=1\}}\left(f(\xi^{i,2})-f(\xi)\right)\notag\\
&+\frac{1}{N}\sum_{i=0}^{N-1}\sum_{j=0}^{N-1}\lambda\left(\frac{i}{N}, \frac{j}{N}\right)1_{\{\xi(i)=0, \xi(j)=2\}}\left(f(\xi^{i,1})-f(\xi)\right)
\end{align}
for any $\xi\in \mathbb{X}$ and $f$ from $\mathbb{X}$ to $\mathbb{R}$, where $\xi^{i,k}\in \mathbb{X}$ is defined as
\[
\xi^{i,k}(j)=
\begin{cases}
\xi(j) &\text{~if~}j\neq i,\\
k & \text{~if~}j=i
\end{cases}
\]
for any $0\leq i\leq N-1$ and $k\in \{0,1,2,3\}$. From now on, to distinguish different $N$, we write $\xi_t$ as $\xi_t^N$ and write $\mathcal{L}$ as $\mathcal{L}^N$.

If $\psi\equiv \psi_0, \phi\equiv \phi_0$ and $\lambda\equiv \lambda_0$ for some $\psi_0, \phi_0, \lambda_0\in (0, +\infty)$, then our model reduces to the spatially homogeneous SEIR model which is an example of the density-dependent Markov chain introduced in \cite{Kurtz1978}. In detail, let
\[
A_t^N=\left(\sum_{i=0}^{N-1}1_{\{\xi_t^N(i)=0\}}, \sum_{i=0}^{N-1}1_{\{\xi_t^N(i)=1\}}, \sum_{i=0}^{N-1}1_{\{\xi_t^N(i)=2\}}\right),
\]
then $A_t^N$ jumps to $A_t^N+\varphi_k$ at rate $N\beta_k\left(\frac{A_t^N}{N}\right)$ for $k=1,2,3$, where
\[
\varphi_1=(0,0,-1), \varphi_2=(0,-1,1), \varphi_3=(-1,1,0)
\]
and
\[
\beta_1(u)=\phi_0u_3, \beta_2(u)=\psi_0u_2, \beta_3(u)=\lambda_0u_1u_2
\]
for any $u=(u_1, u_2, u_3)\in \mathbb{R}^3$. Reference \cite{Pardoux2017} gives large deviation principles of a family of spatially homogeneous epidemic models and Reference \cite{Xue2021} gives moderate deviation principles of a family of density dependent Markov chians. It is easy to check that $\{A_t^N\}_{t\geq 0}$ belongs to stochastic processes investigated in both aforesaid references. Therefore, large and moderate deviation principles of $\{A_t^N\}_{t\geq 0}$ are given by main theorems in \cite{Pardoux2017} and \cite{Xue2021} respectively. In detail, we have following two propositions.

\begin{proposition}\label{Proposition 1.1 LDP for homogeneous} (Pardoux and Samegni-Kepgnou, 2017, \cite{Pardoux2017})
Let $\Omega_1=\mathcal{D}\left([0, T], \mathbb{R}^3\right)$ be the set of c\`{a}dl\`{a}g functions from $[0, T]$ to $\mathbb{R}^3$. If $A_0^N/N=a_0$ for all $N\geq 1$, then
\[
\limsup_{N\rightarrow+\infty}\frac{1}{N}\log P\left(\left\{\frac{A_t^N}{N}\right\}_{0\leq t\leq T}\in C\right)\leq-\inf_{f\in C}K_{T, a_0}(f)
\]
for any closed $C\subseteq \Omega_1$ and
\[
\liminf_{N\rightarrow+\infty}\frac{1}{N}\log P\left(\left\{\frac{A_t^N}{N}\right\}_{0\leq t\leq T}\in O\right)\geq-\inf_{f\in O}K_{T, a_0}(f)
\]
for any open set $O\subseteq \Omega_1$, where $K_{T, a_0}$ is a good rate function.
\end{proposition}

\begin{proposition}\label{Proposition 1.2 MDP for homogeneous} (X., 2021, \cite{Xue2021}) Let $\{\gamma_N\}_{N\geq 1}$ be a positive sequence such that
\[
\lim_{N\rightarrow+\infty}\frac{\gamma_N}{N}=\lim_{N\rightarrow+\infty}\frac{\sqrt{N}}{\gamma_N}=0.
\]
If $A_0^N/N=a_0$ for all $N\geq 1$, then
\[
\limsup_{N\rightarrow+\infty}\frac{N}{\gamma_N^2}\log P\left(\left\{\frac{A_t^N-Na_t}{\gamma_N}\right\}_{0\leq t\leq T}\in C\right)\leq-\inf_{f\in C}Q_{T, a_0}(f)
\]
for any closed $C\subseteq \Omega_1$ and
\[
\liminf_{N\rightarrow+\infty}\frac{N}{\gamma_N^2}\log P\left(\left\{\frac{A_t^N-Na_t}{\gamma_N}\right\}_{0\leq t\leq T}\in O\right)\geq-\inf_{f\in O}Q_{T, a_0}(f)
\]
for any open set $O\subseteq \Omega_1$, where $\{a_t\}_{t\geq 0}$ is the solution to
\[
a_t=a_0+\int_0^t \sum_{k=1}^3\varphi_k\beta_k(a_u)du
\]
for any $t\geq 0$ and $Q_{T, a_0}$ is a good rate function.
\end{proposition}
Precise expressions of $K_{T, a_0}$ and $Q_{T, a_0}$ can be calculated directly by Theorems 4,5 of \cite{Pardoux2017} and Theorem 2.2 of \cite{Xue2021} respectively, which we omit here.

Inspired by Propositions \ref{Proposition 1.1 LDP for homogeneous} and \ref{Proposition 1.2 MDP for homogeneous}, we want to give large and moderate deviations for $\{\xi_t\}_{t\geq 0}$, i.e, the case where transition rates are vertex-dependent. According to the absence of spatial homogeneity of $\{\xi_t\}_{t\geq 0}$, it is not proper to only discuss large and moderate deviation principles for numbers of vertices in respective states. Instead, we are concerned with large and moderate deviations of empirical density fields generated by $\{\xi_t\}_{t\geq 0}$. In detail, let $\delta_a(du)$ be the Dirac measure concentrated on $a$, then in this paper we will give large and moderate deviations of the measured-value process
\[
\left(\sum_{i=0}^{N-1}1_{\{\xi_t^N(i)=0\}}\delta_{\frac{i}{N}}(du), \sum_{i=0}^{N-1}1_{\{\xi_t^N(i)=1\}}\delta_{\frac{i}{N}}(du),
\sum_{i=0}^{N-1}1_{\{\xi_t^N(i)=2\}}\delta_{\frac{i}{N}}(du)\right).
\]
For mathematical details, see Section \ref{section two}.

Law of large numbers and central limit theorems of empirical density fields of interacting particle systems are also called as hydrodynamic limits and fluctuations respectively. For a comprehensive survey of these two topics, see \cite{kipnis+landim99}. Hydrodynamic limits and fluctuations of stochastic epidemic models are investigated in previous literatures. Reference \cite{Xue2022} gives hydrodynamic limits and fluctuations of SIR epidemics with vertex-dependent transition rates by showing that states of different vertices are approximately independent as the total population $N$ grows to infinity. For the SEIR model discussed in this paper, it is not difficult to give analogues of results given in \cite{Xue2022} by utilizing the strategy introduced in \cite{Xue2022} with some details modified. So in this paper we do not discuss fluctuations of our SEIR models. As an application of our large deviation principle, we give hydrodynamic limits of our SEIR models in a different approach than that introduced in \cite{Xue2022}. According our new approach, we can show that the empirical density field of our model converges to its hydrodynamic limit at an exponential rate. For mathematical details, see Section \ref{section two}.

We are also inspired a lot by previous literatures about large and moderate deviations from hydrodynamic limits such as \cite{Kipnis1989} and \cite{Gao2003}, where large and moderate deviations of exclusion processes on lattices are given respectively. Proofs of our large and moderate deviation principles utilize exponential-martingale strategies introduced in \cite{Kipnis1989} and \cite{Gao2003} respectively. For the proof of large deviations of our models, the main new difficulty in mathematics is to show that our rate functions are good. For moderate deviations, main new difficulties are checking that our processes are exponentially tight in the proof of upper bounds and checking that our processes are tight under a transformed measure in the proof of lower bounds. For mathematical details, see Sections \ref{section four} and \ref{section five}.

\section{Main results}\label{section two}
In this section we give our main results. We first introduce some notations and definitions for later use. We denote by $\mathbb{M}$ the set of measures $\nu$ on $\mathbb{T}$ such that $\nu(\mathbb{T})\leq 1$. We let $\mathbb{M}$ be endowed with the weak topology, i.e., $\nu_n\rightarrow \nu$ in $\mathbb{M}$ if and only if
\[
\lim_{n\rightarrow+\infty}\int_\mathbb{T}f(u)\nu_n(du)=\int_\mathbb{T}f(u)\nu(du)
\]
for any $f\in C(\mathbb{T})$. We denote by $\mathcal{S}$ the dual of $C^\infty(\mathbb{T})$ endowed with the weak topology, i.e., $\mathcal{A}_n\rightarrow \mathcal{A}$ in $\mathcal{S}$ if and only if
\[
\lim_{n\rightarrow+\infty}\mathcal{A}_n(f)=\mathcal{A}(f)
\]
for any $f\in C^\infty(\mathbb{T})$. We can consider $\mathbb{M}$ as a subset of $(C(\mathbb{T}))^\prime$ and then $\mathbb{M}\subseteq\mathcal{S}$ by define
\[
\nu(f)=\int_{\mathbb{T}}f(u)\nu(du)
\]
for any $\nu\in \mathbb{M}$ and $f\in C(\mathbb{T})$. For $\mathcal{A}\in (C(T))^\prime$ and $h\in C(\mathbb{T}^2)$, we define $\mathcal{R}_{\mathcal{A}}h$ as the element in $C(\mathbb{T})$ such that
\[
\mathcal{R}_{\mathcal{A}}h(u)=\mathcal{A}(h_u)
\]
for any $u\in \mathbb{T}$, where $h_u(v)=h(u,v)$ for any $v\in \mathbb{T}$. Then, for $\mathcal{A}_1, \mathcal{A}_2\in (C(T))^\prime$, we define $\mathcal{A}_1\bigotimes \mathcal{A}_2$ as the element in $(C(\mathbb{T}^2))^\prime$ such that
\[
\mathcal{A}_1\bigotimes\mathcal{A}_2(h)=\mathcal{A}_1\left(\mathcal{R}_{\mathcal{A}_2}h\right)
\]
for any $h\in C(\mathbb{T}^2)$. Especially, for $\nu_1, \nu_2\in \mathbb{M}$ and $h\in C(\mathbb{T}^2)$,
\[
\nu_1\bigotimes\nu_2(h)=\int_{\mathbb{T}}\int_{\mathbb{T}}h(u,v)\nu_1(du)\nu_2(dv).
\]

For fixed $T>0$, we use $\Omega_2$ to denote $\mathcal{D}\left([0, T], \mathbb{M}^3\right)$, i.e., the set of c\`{a}dl\`{a}g functions from $[0,T]$ to $\mathbb{M}^3$. We use $\Omega_3$ to denote $\mathcal{D}\left([0, T], \mathcal{S}^3\right)$. We let $\Omega_2, \Omega_3$ be endowed with the Skorokhod topology. For any $\varpi\in \mathcal{S}^3$, the $k$-th coordinate of $\varpi$ is denoted by $\varpi_k$ for $k=1,2,3$. For any $W\in \Omega_3$ and $0\leq t\leq T$, the $k$-th coordinate of $W_t$ is denoted by $W_{t,k}$ for $k=1,2,3$.

For any $t\geq 0$ and $0\leq i\leq N-1$, we use $S_t^N(i)$ to denote $1_{\{\xi_t^N(i)=0\}}$. Similarly, we use $\mathcal{E}_t^N(i)$ to denote $1_{\{\xi_t^N(i)=1\}}$ and use $I_t^N(i)$ to denote $1_{\{\xi_t^N(i)=2\}}$. For any $N\geq 1$ and $0\leq t\leq T$, we define
\[
\mu_{t,1}^N(du)=\frac{1}{N}\sum_{i=0}^{N-1}S_t^N(i)\delta_{\frac{i}{N}}(du), \text{~} \mu_{t,2}^N(du)=\frac{1}{N}\sum_{i=0}^{N-1}\mathcal{E}_t^N(i)\delta_{\frac{i}{N}}(du)
\]
and $\mu_{t,3}^N(du)=\frac{1}{N}\sum_{i=0}^{N-1}I_t^N(i)\delta_{\frac{i}{N}}(du)$, where $\delta_a(du)$ is the Dirac measure concentrated on $a$. Then we use $\mu^N$ to denote $\left\{\left(\mu_{t,1}^N, \mu_{t,2}^N, \mu_{t,3}^N\right)\right\}_{0\leq t\leq T}$ and hence $\mu^N\in \Omega_2\subseteq \Omega_3$.

To give the large deviation rate function of the dynamic of our process, we need following definitions. For any $W\in \Omega_3$ and $F,G,H\in C^{1,0}\left([0, T]\times \mathbb{T}\right)$, we define
\begin{align*}
l_1(W,F,G,H)=&W_{T,1}(F_T)+W_{T,2}(G_T)+W_{T,3}(H_T)\\
&-W_{0,1}(F_0)-W_{0,2}(G_0)-W_{0,3}(H_0)\\
&-\int_0^TW_{s,1}(\partial_sF_s)+W_{s,2}(\partial_sG_s)+W_{s,3}(\partial_sH_s)ds.
\end{align*}
For any $W\in \Omega_3$ and $F,G, H\in C^{1,0}\left([0, T]\times \mathbb{T}\right)$, we define
\[
\mathcal{B}_1(W,F,G)=\int_0^T W_{s,1}\bigotimes W_{s,3}\left(h_{1,s}^{F,G}\right)ds,
\]
where $h^{F,G}_{1,s}\in C(\mathbb{T}^2)$ such that
\[
h^{F,G}_{1,s}(u,v)=\lambda(u,v)\left(\exp(-F_s(u)+G_s(u))-1\right)
\]
for any $u,v\in \mathbb{T}$. Further we define
\[
\mathcal{B}_2(W,G,H)=\int_0^TW_{s,2}\left(h_{2,s}^{G,H}\right)ds,
\]
where $h_{2,s}^{G,H}\in C(\mathbb{T})$ such that
\[
h_{2,s}^{G,H}(u)=\psi(u)\left(\exp(-G_s(u)+H_s(u))-1\right)
\]
for any $u\in \mathbb{T}$. Further we define
\[
\mathcal{B}_3(W,H)=\int_0^TW_{s,3}\left(h_{3,s}^{H}\right)ds,
\]
where $h_{3,s}^{H}\in C(\mathbb{T})$ such that
\[
h_{3,s}^{H}(u)=\phi(u)\left(\exp(-H_s(u))-1\right)
\]
for any $u\in \mathbb{T}$. Then we define
\[
\mathcal{I}_1(W, F, G,H)=l_1(W,F,G,H)-\mathcal{B}_1(W,F,G)-\mathcal{B}_2(W,G,H)-\mathcal{B}_3(W,H).
\]

Now we can give the large deviation rate function of the dynamic of our process. For any $W\in \Omega_2$, we define
\begin{equation}\label{equ 2.1 LDPdynamicRatefunction}
I_{dyn}(W)=\sup_{F,G,H\in C^{1,0}([0,T]\times \mathbb{T})}\left\{\mathcal{I}_1(W,F,G,H)\right\}.
\end{equation}
To give large deviation rate function of the initial state of our process, we first give the initial condition of our process. Throughout this paper, we adopt the following assumption.

\textbf{Assumption} (A): For all $N\geq 1$, $\{\xi_0^N(i)\}_{0\leq i\leq N-1}$ are independent. Furthermore, for all $0\leq i\leq N-1$,
\[
P\left(\xi_0(i)=0\right)=\rho_0\left(\frac{i}{N}\right), \text{~}P\left(\xi_0(i)=1\right)=\rho_1\left(\frac{i}{N}\right)
\]
and $P\left(\xi_0(i)=2\right)=\rho_2\left(\frac{i}{N}\right)$, where $\rho_0, \rho_1, \rho_2\in C^\infty(\mathbb{T})$ are strictly positive and
\[
\rho_0(u)+\rho_1(u)+\rho_2(u)<1
\]
for all $u\in \mathbb{T}$.

Let $\rho_0, \rho_1, \rho_2$ be defined as in Assumption (A). For any $f_1,f_2,f_3\in C(\mathbb{T})$ and $\varpi\in \mathbb{M}^3$, we define
\[
\mathcal{I}_2\left(\varpi, f_1, f_2, f_3\right)
=\sum_{k=1}^3\varpi_k(f_k)-\int_\mathbb{T}\log\left(1+\sum_{k=1}^3\rho_{k-1}(u)\left(\exp\left(f_k(u)\right)-1\right)\right)du.
\]
Now we give the large deviation rate function of the initial state of our process. For any $\varpi\in \mathbb{M}^3$, we define
\begin{equation}\label{equ 2.2 LDPinitialRatefunction}
I_{ini}(\varpi)=\sup_{f_1,f_2,f_3\in C(\mathbb{T})}\{\mathcal{I}_2(\varpi, f_1, f_2, f_3)\}.
\end{equation}
Before giving our large deviation principle, we first give two lemmas about properties of our rate functions.
\begin{lemma}\label{lemma 2.1 good LDP rate}
Let $I_{dyn}$ be defined as in Equation \eqref{equ 2.1 LDPdynamicRatefunction}, then $I_{dyn}$ is a good rate function.
\end{lemma}
Note that it is a trivial conclusion that $I_{ini}$ is good since $\varpi_k(\mathbb{T})\leq 1$ for any $\varpi\in \mathbb{M}^3$ and $k=1,2,3$.

To give the second property of our rate function, we use $D_0$ to denote the set of elements $W\in \Omega_2$ satisfies following three properties.

(1) For $k=1,2,3$ and $t\in [0, T]$,
\[
W_{t,k}(du)=w_{t,k}(u)du
\]
for some $w_k\in C^{1,0}\left([0, T]\times \mathbb{T}\right)$.

(2) For any $0\leq t\leq T$ and $u\in \mathbb{T}$,
\[
0<w_{t,1}(u)<w_{t,1}(u)+w_{t,2}(u)<w_{t,1}(u)+w_{t,2}(u)+w_{t,3}(u)<1.
\]

(3) For any $0\leq t\leq T$ and $u\in \mathbb{T}$,
\[
\frac{d}{dt}\sum_{r=1}^kw_{t,r}(u)<0
\]
for $k=1,2,3$.

Note that elements in $D_0$ are feasible weak limits of $\{\mu^N\}_{N\geq 1}$ since
\[
0\leq S_t^N(i)\leq S_t^N(i)+\mathcal{E}_t^N(i)\leq S_t^N(i)+\mathcal{E}_t^N(i)+I_t^N(i)\leq 1
\]
and
\[
S_t^N(i),\text{~}S_t^N(i)+\mathcal{E}_t^N(i) \text{~and~}S_t^N(i)+\mathcal{E}_t^N(i)+I_t^N(i)
\]
are all decreasing for all $N\geq 1$ and $0\leq i\leq N-1$.

Now we give the second property of our rate function.

\begin{lemma}\label{lemma 2.2 rate function on D0}
If $W\in D_0$, then $I_{ini}(W_0)+I_{dyn}(W)<+\infty$. Furthermore,
\begin{align}\label{equ 2.3 Ini RateInitialOnD0}
I_{ini}(W_0)=&\int_{\mathbb{T}}\Bigg(\sum_{k=1}^3w_{0,k}(u)\log\left(\frac{w_{0,k}(u)}{\rho_{k-1}(u)}\right)\\
&\text{\quad}-\left(1-\sum_{k=1}^3w_{0,k}(u)\right)\log\left(\frac{1-\sum_{k=1}^3\rho_{k-1}(u)}{1-\sum_{k=1}^3w_{0,k}(u)}\right)\Bigg)du \notag
\end{align}
and
\begin{align}\label{equ 2.4 Dyn rate function on D0}
I_{dyn}(W)=\mathcal{I}_1\left(W, F_W, G_W, H_W\right),
\end{align}
where
\[
H_W(t,u)=-\log\left(\frac{-\frac{d}{dt}\sum_{k=1}^3w_{t,k}(u)}{\phi(u)w_{t,3}(u)}\right),
\]
\[
G_W(t,u)=H_W(t,u)-\log\left(\frac{-\frac{d}{dt}\sum_{k=1}^2w_{t,k}(u)}{\psi(u)w_{t,2}(u)}\right)
\]
and
\[
F_W(t,u)=G_W(t,u)-\log\left(\frac{-\frac{d}{dt}w_{t,1}(u)}{w_{t,1}(u)\int_{\mathbb{T}}\lambda(u,v)w_{t,3}(v)dv}\right)
\]
for all $0\leq t\leq T$ and $u\in \mathbb{T}$.
\end{lemma}
Now we give our large deviation principle.

\begin{theorem}\label{Theorem 2.3 LDP}
Under Assumption (A),
\begin{equation}\label{equ 2.5 LDP upperbound}
\limsup_{N\rightarrow+\infty}\frac{1}{N}\log P(\mu^N\in C)\leq -\inf_{W\in C}\left(I_{ini}(W_0)+I_{dyn}(W)\right)
\end{equation}
for any closed set $C\subseteq \Omega_2$ and
\begin{equation}\label{equ 2.6 LDP lowerbound}
\liminf_{N\rightarrow+\infty}\frac{1}{N}\log P(\mu^N\in O)\geq -\inf_{W\in O\cap D_0}\left(I_{ini}(W_0)+I_{dyn}(W)\right)
\end{equation}
for any open set $O\subseteq \Omega_2$.
\end{theorem}
For any $W\in \Omega_2$, we guess that $I_{ini}(W_0)+I_{dyn}(W)<+\infty$ implies that there exist a series $\{W^n\}_{n\geq 1}$ in $D_0$ such that $W^n\rightarrow W$ and $I_{ini}(W_0^n)+I_{dyn}(W^n)\rightarrow I_{ini}(W_0)+I_{dyn}(W)$ as $n\rightarrow+\infty$. If this guess holds, then the righthand side of Equation \eqref{equ 2.6 LDP lowerbound} can be strengthen to $-\inf_{W\in O}\left(I_{ini}(W_0)+I_{dyn}(W)\right)$. We will work on this guess as a further investigation.

As an application of Theorem \ref{Theorem 2.3 LDP}, we can show that $\mu^N$ converges weakly to its hydrodynamic limit at an exponential rate. We first give the expression of this hydrodynamic limit. Let $\left\{(\theta_t^S(u), \theta_t^E(u), \theta_t^I(u))\right\}_{0\leq t\leq T, u\in \mathbb{T}}\in C^{1,0}\left([0, T]\times\mathbb{T}, [0, 1]^3\right)$ be the solution to the ordinary differential equation
\begin{equation}\label{equ 2.7 non linear hydrodynamic ODE}
\begin{cases}
&\frac{d}{dt}\theta_t^S(u)=-\theta_t^S(u)\int_{\mathbb{T}}\lambda(u,v)\theta_t^I(v)dv,\\
&\frac{d}{dt}\theta_t^E(u)=\theta_t^S(u)\int_{\mathbb{T}}\lambda(u,v)\theta_t^I(v)dv-\psi(u)\theta_t^E(u),\\
&\frac{d}{dt}\theta_t^I(u)=\psi(u)\theta_t^E(u)-\phi(u)\theta_t^I(u),\\
&\theta_0^S(u)=\rho_0(u),\\
&\theta_0^E(u)=\rho_1(u),\\
&\theta_0^I(u)=\rho_2(u),
\end{cases}
\end{equation}
where $\rho_0, \rho_1, \rho_2$ are defined as in Assumption (A). According to Grownwall's inequality, it is easy to check that the solution to Equation \eqref{equ 2.7 non linear hydrodynamic ODE} is unique. If $\lambda\equiv \lambda_0, \psi\equiv \psi_0$ and $\phi\equiv \phi_0$ for some $\lambda_0, \psi_0, \phi_0\in (0, +\infty)$, then
\[
(\hat{\theta}_t^S, \hat{\theta}_t^E, \hat{\theta}_t^I):=\int_\mathbb{T}\left(\theta_t^S(u), \theta_t^E(u), \theta_t^I(u)\right)du
\]
is driven by
\[
\begin{cases}
&\frac{d}{dt}\hat{\theta}_t^S=-\lambda_0\hat{\theta}_t^S\hat{\theta}_t^I,\\
&\frac{d}{dt}\hat{\theta}_t^E=\lambda_0\hat{\theta}_t^S\hat{\theta}_t^I-\psi_0\hat{\theta}_t^E,\\
&\frac{d}{dt}\hat{\theta}_t^I=\psi_0\hat{\theta}_t^E-\phi_0\hat{\theta}_t^I,
\end{cases}
\]
which is the classic deterministic version of spatially homogeneous SEIR model (see \cite{Anderson1991}).

Let $\mu\in \Omega_2$ such that
\[
\mu_{t,1}(du)=\theta_t^S(u)du, \text{~}\mu_{t,2}(du)=\theta_t^E(u)du\text{~and~}\mu_{t,3}(du)=\theta_t^I(u)du
\]
for any $0\leq t\leq T$, then we have the following conclusion.

\begin{theorem}\label{Theorem 2.4 exponentially rapid to hydrodynamic limit}
For any open set $O\subseteq \Omega_2$ such that $O\ni\mu$, there exists $\alpha_1=\alpha_1(O)>0$ independent of $N$ such that
\[
P\left(\mu^N\in O\right)\geq 1-\exp(\alpha_1N)
\]
for sufficiently large $N$.
\end{theorem}

By utilizing the strategy given in \cite{Xue2022} with some details modified, it is not difficult to show that $\mu^N$ converges to $\mu$ weakly. Theorem \ref{Theorem 2.4 exponentially rapid to hydrodynamic limit} gives a stronger conclusion that this convergence is exponentially rapid.

Now we introduce our result about the moderate deviation principle of $\mu^N$. We first introduce some notations and definitions for later use. Throughout this paper, let $\{\gamma_N\}_{N\geq 1}$ be a positive sequence such that
\[
\lim_{N\rightarrow+\infty}\frac{\gamma_N}{N}=\lim_{N\rightarrow+\infty}\frac{\sqrt{N}}{\gamma_N}=0.
\]
For any $0\leq t\leq T$, we define $\eta^N_{t,1}(du)=\frac{\sum_{i=1}^N\left(S_t^N(i)-\mathbb{E}S_t^N(i)\right)\delta_{\frac{i}{N}}(du)}{\gamma_N}$,
\[
\eta^N_{t,2}(du)=\frac{\sum_{i=1}^N\left(\mathcal{E}_t^N(i)-\mathbb{E}\mathcal{E}_t^N(i)\right)\delta_{\frac{i}{N}}(du)}{\gamma_N}
\]
and $\eta^N_{t,3}(du)=\frac{\sum_{i=1}^N\left(I_t^N(i)-\mathbb{E}I_t^N(i)\right)\delta_{\frac{i}{N}}(du)}{\gamma_N}$, where $\mathbb{E}$ is the expectation operator. Let $\eta^N=\left\{\left(\eta_{t,1}^N, \eta_{t,2}^N, \eta_{t,3}^N\right)\right\}_{0\leq t\leq T}$, then $\eta^N\in \Omega_3$.

To give the moderate deviation rate function of the dynamic of our process, we introduce following definitions. For any $W\in \Omega_3$ and $F,G,H\in C^{1,\infty}([0, T]\times\mathbb{T})$, we define
\[
\mathcal{B}_4(F,G)=\int_0^T\int_{\mathbb{T}^2}\lambda(u,v)\theta_s^S(u)\theta_s^I(v)\left(G_s(u)-F_s(u)\right)^2dsdudv,
\]
\[
\mathcal{B}_5(G,H)=\int_0^T\int_{\mathbb{T}}\psi(u)\theta_s^E(u)\left(-G_s(u)+H_s(u)\right)^2dsdu,
\]
and
\[
\mathcal{B}_6(H)=\int_0^T\int_\mathbb{T}\phi(u)\theta_s^I(u)H_s^2(u)dsdu,
\]
where $\theta^S, \theta^E, \theta^I$ are defined as in Equation \eqref{equ 2.7 non linear hydrodynamic ODE}. Furthermore, we define
\[
\mathcal{B}_7(W,F,G)=\int_0^TW_{s,1}\left(\mathcal{P}_{1,s}\left(-F_s+G_s\right)\right)+W_{s,3}\left(\mathcal{P}_{2,s}(-F_s+G_s)\right)ds,
\]
\[
\mathcal{B}_8(W,G,H)=\int_0^T W_{s,2}\left(\mathcal{P}_3(-G_s+H_s)\right)ds
\]
and
\[
\mathcal{B}_9(W,H)=\int_0^TW_{s,3}\left(\mathcal{P}_4(-H_s)\right)ds,
\]
where $\mathcal{P}_{1,s}, \mathcal{P}_{2,s}, \mathcal{P}_3, \mathcal{P}_4$ are linear operators on $C^\infty(\mathbb{T})$ such that
\[
\mathcal{P}_{1,s}f(u)=f(u)\int_{\mathbb{T}}\theta_s^I(v)\lambda(u,v)dv,
\]
\[
\mathcal{P}_{2,s}f(u)=\int_{\mathbb{T}}\theta_s^S(v)\lambda(v,u)f(v)dv,
\]
\[
\mathcal{P}_3f(u)=\psi(u)f(u),
\]
and
\[
\mathcal{P}_4f(u)=\phi(u)f(u)
\]
for any $f\in C^\infty(\mathbb{T})$ and $u\in \mathbb{T}$. Then we define
\[
l_2(W,F,G,H)=l_1(W,F,G,H)-\mathcal{B}_7(W,F,G)-\mathcal{B}_8(W,G,H)-\mathcal{B}_9(W,H)
\]
and
\[
\mathcal{J}_1(W,F,G,H)=l_2(W,F,G,H)-\frac{1}{2}\mathcal{B}_4(F,G)-\frac{1}{2}\mathcal{B}_5(G,H)-\frac{1}{2}\mathcal{B}_6(H).
\]
Now we give the moderate deviation rate function of the dynamic of our process. For any $W\in \Omega_3$, we define
\begin{equation}\label{equ 2.8 MDPdynamicRatefunction}
J_{dyn}(W)=\sup_{F,G,H\in C^{1,\infty}([0,T]\times \mathbb{T})}\left\{\mathcal{J}_1(W,F,G,H)\right\}.
\end{equation}

To give the moderate deviation rate function of the initial state of our process, we define
\begin{align*}
\mathcal{J}_2(\varpi, f_1, f_2, f_3)=&\sum_{k=1}^3\varpi_k(f_k)\\
&-\frac{1}{2}\int_{\mathbb{T}}\sum_{k=1}^3\rho_{k-1}(u)f_k^2(u)-\left(\sum_{k=1}^3f_k(u)\rho_{k-1}(u)\right)^2du
\end{align*}
for any $\varpi\in \mathcal{S}^3$ and $f_1, f_2, f_3\in C^\infty(\mathbb{T})$. Now we give the moderate deviation rate function of the initial state of our process. For any $\varpi\in \mathcal{S}^3$, we define
\begin{equation}\label{equ 2.9 MDPinitialRatefunction}
J_{ini}(\varpi)=\sup_{f_1,f_2,f_3\in C^\infty(\mathbb{T})}\left\{\mathcal{J}_2(\varpi, f_1, f_2, f_3)\right\}.
\end{equation}

According to some technical reasons, currently we can only prove our moderate deviation principles under the following assumption.

\textbf{Assumption} (B): there exist $\lambda_1, \lambda_2\in C^\infty(\mathbb{T})$ such that
\[
\lambda(u,v)=\lambda_1(u)\lambda_2(v)
\]
for any $u,v\in \mathbb{T}$.

Now we give our moderate deviation principles.

\begin{theorem}\label{Theorem 2.5 MDP}
Under Assumptions (A) and (B),
\begin{equation}\label{equ 2.10 MDP upperbound}
\limsup_{N\rightarrow+\infty}\frac{N}{\gamma_N^2}\log P\left(\eta^N\in C\right)\leq -\inf_{W\in C}\left(J_{ini}(W_0)+J_{dyn}(W)\right)
\end{equation}
for any closed set $C\subseteq \Omega_3$ and
\begin{equation}\label{equ 2.11 MDP lowerbound}
\liminf_{N\rightarrow+\infty}\frac{N}{\gamma_N^2}\log P\left(\eta^N\in O\right)\geq -\inf_{W\in O}\left(J_{ini}(W_0)+J_{dyn}(W)\right)
\end{equation}
for any open set $O\subseteq \Omega_3$.
\end{theorem}
It is natural to ask whether Theorem \ref{Theorem 2.5 MDP} holds without Assumption (B). We will work on this question as a further investigation.

As an application of Theorem \ref{Theorem 2.5 MDP}, moderate deviation principles of hitting times of $\mu^N$ can be given. In detail, let $f_1,f_2,f_3\in C^\infty(\mathbb{T})$ such that
\[
\frac{d}{dt}\sum_{k=1}^3\mu_{t,k}(f_k)>0
\]
for $0\leq t\leq T_1$ for some $T_1>0$. For $c\in \left(\sum_{k=1}^3\mu_{0,k}(f_k), \sum_{k=1}^3\mu_{T_1,k}(f_k)\right)$ and $N\geq 1$, let
\[
\tau_c=\inf\left\{t:~\sum_{k=1}^3\mu_{t,k}(f_k)=c\right\}
\]
and
\[
\tau_c^N=\inf\left\{t:~\sum_{k=1}^3\mu_{t,k}^N(f_k)\geq c\right\},
\]
then we have the following theorem.
\begin{theorem}\label{Theorem 2.6 MDPofHittingTimes}
For any $x\in \mathbb{R}$ and $t>0$, let
\[
J_{contra, t}(x)=\inf\left\{J_{ini}(W_0)+J_{dyn, t}(W):~\sum_{k=1}^3W_{t,k}(f_k)=x\right\},
\]
where $J_{dyn,t}$ is the $J_{dyn}$ defined as in Equation \eqref{equ 2.8 MDPdynamicRatefunction} with $T=t$, then
\[
J_{contra, t}(x)=x^2J_{contra, t}(1)
\]
for all $x\in \mathbb{R}$ and
$\frac{N}{\gamma_N}\left(\tau_c^N-\tau_c\right)$ satisfies a moderate deviation principle under Assumptions (A) and (B) at strength $\frac{\gamma_N^2}{N}$ with rate function $J_{hit}$ given by
\[
J_{hit}(x)=x^2J_{contra, \tau_c}(1)\left(\frac{d}{dt}\sum_{k=1}^3\mu_{t,k}(f_k)\Bigg|_{t=\tau_c}\right)^2
\]
for all $x\in \mathbb{R}$.
\end{theorem}

The rest of this paper is arranged as follows. In Section \ref{section three} we prove Lemmas \ref{lemma 2.1 good LDP rate} and \ref{lemma 2.2 rate function on D0}. The proof of Lemma \ref{lemma 2.1 good LDP rate} utilizes Arezl\`{a}-Ascoli theorem. The proof of Lemma \ref{lemma 2.2 rate function on D0} utilizes calculus of variation. In Section \ref{section four} we prove our large deviation principle Theorem \ref{Theorem 2.3 LDP} and its application Theorem \ref{Theorem 2.4 exponentially rapid to hydrodynamic limit}. The proof of Theorem \ref{Theorem 2.3 LDP} utilizes the exponential-martingale strategy introduced in \cite{Kipnis1989}. The proof of Theorem \ref{Theorem 2.4 exponentially rapid to hydrodynamic limit} relies heavily on Lemma \ref{lemma 2.1 good LDP rate}, which ensures that our rate functions are good. In Section \ref{section five} we prove our moderate deviation principle Theorem \ref{Theorem 2.5 MDP} and give an outline of the proof of Theorem \ref{Theorem 2.6 MDPofHittingTimes}. The proof of Theorem \ref{Theorem 2.5 MDP} utilizes the exponential-martingale strategy introduced in \cite{Gao2003}. Mathematical difficulties are mainly in checks of the exponential tightness of our processes and the tightness of our processes under a transformed probability measure. With Theorems \ref{Theorem 2.3 LDP} and \ref{Theorem 2.5 MDP}, the proof of Theorem \ref{Theorem 2.6 MDPofHittingTimes} follows from an analysis similar with that introduced in \cite{He2023}, which gives moderate deviation principles of hitting times of density-dependent Markov chains. 

\section{Proofs of Lemmas \ref{lemma 2.1 good LDP rate} and \ref{lemma 2.2 rate function on D0}}\label{section three}
In this section we prove Lemmas \ref{lemma 2.1 good LDP rate} and \ref{lemma 2.2 rate function on D0}. We first prove Lemma \ref{lemma 2.1 good LDP rate}.

\proof[Proof of Lemma \ref{lemma 2.1 good LDP rate}]
Since $\mathcal{I}_1(W,0,0,0)=0$, $I_{dyn}$ is non-negative. Since $\mathcal{I}_1$ is continuous, $I_{dyn}$ is lower-semicontinuous. As a result, we only need to show that $\mathcal{C}_a$ is compact for any $a>0$, where
\[
\mathcal{C}_a=\left\{W\in \Omega_2:\text{~}I_{dyn}(W)\leq a\right\}.
\]
According to Proposition 1.2 in Chapter 4 of \cite{kipnis+landim99}, which is an extension of Arzel\`{a}-Ascoli theorem, $\mathcal{C}_a$ is compact if we can show that
\[
\left\{\{W_{t,k}(f)\}_{0\leq t\leq T}:\text{~}W\in \mathcal{C}_a\right\}
\]
are equicontinuous for any given $f\in C(\mathbb{T})$, $k\in \{1,2,3\}$ and $a>0$. We only discuss the case where $k=1$ since other two cases can be discussed in the same way. For given $0<t_1<t_2<T$ and sufficiently large $n$, let $q_n\in L^1([0,T])$ such that $q_n(s)=n$ for $t_1\leq s\leq t_1+\frac{1}{n}$, $q_n(s)=-n$ for $t_2\leq s\leq t_2+\frac{1}{n}$ and $q_n(s)=0$ otherwise. Since $C([0, T])$ is dense in $L^1([0, T])$, there exist a series $\{q_{n,m}\}_{m\geq 1}$ in $C[0, T]$ such that $q_{n,m}$ converges to $q_n$ in $L^1$ as $m\rightarrow+\infty$. For any $0\leq s\leq T$, let $r_n(s)=\int_0^sq_n(s_1)ds_1$ and $r_{n,m}(s)=\int_0^sq_{n,m}(s_1)ds_1$, then $r_{n,m}$ converges to $r_n$ in $C([0, T])$ as $m\rightarrow+\infty$. For given $f\in C(\mathbb{T})$, let $F_{n,m}\in C^{1,0}\left([0, T]\times \mathbb{T}\right)$ such that
\[
F_{n,m}(s,u)=r_{n,m}(s)f(u)
\]
for all $(s,u)\in [0, T]\times \mathbb{T}$. Then, for any $W\in \mathcal{C}_a$,
\[
\mathcal{I}_1\left(W,F_{n,m},0,0\right)\leq a
\]
for any $m\geq 1$. Since $r_{n,m}\rightarrow r_n$ in $C([0, T])$ as $m\rightarrow+\infty$ and $W_{s,k}(\mathbb{T})\leq 1$ for $k=1,3$,
\begin{align*}
\alpha \geq &\lim_{m\rightarrow+\infty}\mathcal{I}_1\left(W,F_{n,m},0,0\right)\\
=&-n\int_{t_1}^{t_1+\frac{1}{n}}W_{s,1}(f)ds+n\int_{t_2}^{t_2+\frac{1}{n}}W_{s,1}(f)ds\\
&-\int_{t_1}^{t_1+\frac{1}{n}}W_{1,s}\bigotimes W_{3,s}\left(\lambda(\cdot,\ast)(e^{r_n(s)f(\cdot)}-1)\right)ds\\
&-\int_{t_2}^{t_2+\frac{1}{n}}W_{1,s}\bigotimes W_{3,s}\left(\lambda(\cdot,\ast)(e^{r_n(s)f(\cdot)}-1)\right)ds\\
&-\int_{t_1+\frac{1}{n}}^{t_2}W_{1,s}\bigotimes W_{3,s}\left(\lambda(\cdot,\ast)(e^{f(\cdot)}-1)\right)ds.
\end{align*}
Since $0\leq r_n(s)\leq 1$, let $n\rightarrow+\infty$, then
\[
\alpha\geq W_{t_2,1}(f)-W_{t_1,1}(f)-\int_{t_1}^{t_2}W_{1,s}\bigotimes W_{3,s}\left(\lambda(\cdot,\ast)(e^{f(\cdot)}-1)\right)ds.
\]
In the above equation, $f$ is arbitrary. Hence, let $f$ be replaced by $Kf$ for any $K>0$, then
\begin{equation}\label{equ 3.1}
W_{t_2,1}(f)-W_{t_1,1}(f)\leq \frac{\alpha}{K}+\int_{t_1}^{t_2}W_{1,s}\bigotimes W_{3,s}\left(\lambda(\cdot,\ast)(e^{Kf(\cdot)}-1)\right)ds.
\end{equation}
For any $\epsilon>0$, let $K=\frac{4\alpha}{\epsilon}$ and $\delta=\frac{\epsilon}{4\|\lambda\|_\infty\left(e^{K\|f\|_\infty}+1\right)}$, then by Equation \eqref{equ 3.1},
\[
W_{t_2,1}(f)-W_{t_1,1}(f)\leq \epsilon
\]
for any $t_1<t_2$ such that $t_2-t_1\leq \delta$. Let $f$ be replaced by $-Kf$ for any $K>0$, then a similar analysis implies that
\[
W_{t_2,1}(f)-W_{t_1,1}(f)\geq -\epsilon
\]
for any $t_1<t_2$ such that $t_2-t_1\leq \delta$. Since the choose of $\delta$ is independent of $W\in \mathcal{C}_a$,
\[
\left\{\{W_{t,k}(f)\}_{0\leq t\leq T}:\text{~}W\in \mathcal{C}_a\right\}
\]
are equicontinuous and the proof is complete.

\qed

Now we prove Lemma \ref{lemma 2.2 rate function on D0}.

\proof[Proof of Lemma \ref{lemma 2.2 rate function on D0}]

Our proof follows calculus of variation. We only prove Equation \eqref{equ 2.4 Dyn rate function on D0} since Equation \eqref{equ 2.3 Ini RateInitialOnD0} can be proved according to a similar analysis. For any $W\in D_0$, according to our definitions of $F_W, G_W, H_W$,
\begin{equation}\label{equ 3.2}
\begin{cases}
&\frac{-\frac{d}{dt}w_{t,1}(u)}{w_{t,1}(u)\int_{\mathbb{T}}\lambda(u,v)w_{t,3}(v)dv}=\exp\left(-F_W(t,u)+G_W(t,u)\right),\\
&\frac{-\frac{d}{dt}\sum_{k=1}^2w_{t,k}(u)}{\psi(u)w_{t,2}(u)}=\exp\left(-G_W(t,u)+H_W(t,u)\right),\\
&\frac{-\frac{d}{dt}\sum_{k=1}^3w_{t,k}(u)}{\phi(u)w_{t,3}(u)}=\exp\left(-H_W(t,u)\right)
\end{cases}
\end{equation}
for any $0\leq t\leq T$ and $u\in \mathbb{T}$. For $W\in D_0$, according to the formula of integration by parts,
\begin{align}\label{equ 3.3}
\mathcal{I}_1(W,F,G,H)=&\int_0^T\int_\mathbb{T}\frac{d}{ds}w_{s,1}(u)F_s(u)dsdu+\int_0^T\int_{\mathbb{T}}\frac{d}{ds}w_{s,2}(u)G_s(u)dsdu \notag\\
&+\int_0^T\int_{\mathbb{T}}\frac{d}{ds}w_{s,3}(u)H_s(u)dsdu \\
&-\int_0^T\int_{\mathbb{T}^2}w_{s,1}(u)w_{s,3}(v)h_{1,s}^{F,G}(u,v)dsdudv\notag\\
&-\int_0^T\int_{\mathbb{T}}w_{s,2}(u)h_{2,s}^{G,H}(u)du-\int_0^T\int_{\mathbb{T}}w_{s,3}(u)h_{3,s}^H(u)du \notag
\end{align}
for any $F,G,H\in C^{1,0}\left([0, T]\times\mathbb{T}\right)$. For any $\hat{F},\hat{G},\hat{H}\in C^{1,0}\left([0, T]\times\mathbb{T}\right)$, it is easy to check that
\[
\frac{d^2}{d\epsilon^2}\mathcal{I}_1(W,F_W+\epsilon \hat{F}, G_W+\epsilon \hat{G}, H_W+\epsilon \hat{H})\leq 0
\]
for all $\epsilon\in \mathbb{R}$ and
\[
\frac{d}{d\epsilon}\mathcal{I}_1(W,F_W+\epsilon \hat{F}, G_W+\epsilon \hat{G}, H_W+\epsilon \hat{H})\Big|_{\epsilon=0}=0
\]
according to Equations \eqref{equ 3.2} and \eqref{equ 3.3}. Hence,
\begin{equation}\label{equ 3.4}
\mathcal{I}_1(W,F_W, G_W, H_W)\geq\mathcal{I}_1(W,F_W+\epsilon \hat{F}, G_W+\epsilon \hat{G}, H_W+\epsilon \hat{H})
\end{equation}
for any $\epsilon\in \mathbb{R}$. For any $F, G, H\in C^{1,0}\left([0, T]\times\mathbb{T}\right)$, let $\hat{F}=F-F_W$, $\hat{G}=G-G_W$, $\hat{H}=H-H_W$ and $\epsilon=1$ in Equation \eqref{equ 3.4}, then
\[
\mathcal{I}_1(W,F_W, G_W, H_W)\geq \mathcal{I}_1(W,F, G, H).
\]
Hence
\[
I_{dyn}(W)=\sup_{F,G,H\in C^{1,0}([0,T]\times \mathbb{T})}\left\{\mathcal{I}_1(W,F,G,H)\right\}=\mathcal{I}_1(W,F_W, G_W, H_W)
\]
and the proof is complete.

\qed

\section{Proofs of Theorems \ref{Theorem 2.3 LDP} and \ref{Theorem 2.4 exponentially rapid to hydrodynamic limit}}\label{section four}
In this section we prove Theorems \ref{Theorem 2.3 LDP} and \ref{Theorem 2.4 exponentially rapid to hydrodynamic limit}. The proof of Theorem \ref{Theorem 2.3 LDP} follows the exponential-martingale strategy introduced in \cite{Kipnis1989}. The upper bound part of the proof is given in Subsection \ref{subsection 4.1} and the lower bound part is given in Subsection \ref{subsection 4.2}. As an application of Theorem \ref{Theorem 2.4 exponentially rapid to hydrodynamic limit}, the proof of Theorem \ref{Theorem 2.4 exponentially rapid to hydrodynamic limit} is given in Subsection \ref{subsection 4.3}, where Lemma \ref{lemma 2.1 good LDP rate} plays the key role.

\subsection{Proof of Equation \eqref{equ 2.5 LDP upperbound}}\label{subsection 4.1}
In this subsection we prove Equation \eqref{equ 2.5 LDP upperbound}. We first introduce some nations and definitions for later use. For any $N\geq 1$, $0\leq t\leq T$ and $F, G, H\in C^{1,0}\left([0, T]\times\mathbb{T}\right)$, we define
\[
V^N_{F,G,H}\left(t, \xi^N\right)=\exp\left(N\left(\mu_{t,1}^N(F_t)+\mu_{t,2}^N(G_t)+\mu_{t,3}^N(H_t)\right)\right).
\]
Further, we define
\[
\mathcal{U}^N_{F,G,H}\left(t, \xi^N\right)=\frac{V_{F,G,H}^N(t, \xi^N)}{V^N_{F,G,H}(0,\xi^N)}\exp\left(-\int_0^t
\frac{\left(\partial_s+\mathcal{L}^N\right)V_{F,G,H}^N(s, \xi^N)}{V_{F,G,H}^N(s, \xi^N)}ds\right),
\]
where $\mathcal{L}^N$ is the generator of our process given in \eqref{equ generator}. According to Feynman-Kac formula, $\{\mathcal{U}^N_{F,G,H}\left(t, \xi^N\right)\}_{0\leq t\leq T}$ is an exponential martingale. According to the definition of $\mathcal{L}^N$ and direct calculation,
\begin{equation}\label{equ 4.1.1}
\mathcal{U}^N_{F,G,H}\left(T, \xi^N\right)=\exp\left(N\left(\mathcal{I}_1\left(\mu^N, F, G, H\right)\right)\right).
\end{equation}

We first show that Equation \eqref{equ 2.5 LDP upperbound} holds when $C\subseteq \Omega_2$ is compact.
\begin{lemma}\label{lemma 4.1.1 LDPupperboundforCompact}
Under Assumption (A),
\begin{equation*}
\limsup_{N\rightarrow+\infty}\frac{1}{N}\log P(\mu^N\in C)\leq -\inf_{W\in C}\left(I_{ini}(W_0)+I_{dyn}(W)\right)
\end{equation*}
for any compact set $C\subseteq \Omega_2$.
\end{lemma}

\proof[Proof of Lemma \ref{lemma 4.1.1 LDPupperboundforCompact}]

For any $f_1, f_2, f_3\in C(\mathbb{T})$ and $F, G, H\in C^{1,0}\left([0, T]\times \mathbb{T}\right)$, since $\{\mathcal{U}^N_{F,G,H}\left(t, \xi^N\right)\}_{0\leq t\leq T}$ is an exponential martingale and $\mathcal{U}^N_{F,G,H}\left(0, \xi^N\right)=1$, we have
\begin{align}\label{equ 4.1.2}
\mathbb{E}\exp\left(N\sum_{k=1}^3\mu^N_{0,k}(f_k)\right)&=\mathbb{E}\left(\exp\left(N\sum_{k=1}^3\mu^N_{0,k}(f_k)\right)\mathcal{U}^N_{F,G,H}\left(T, \xi^N\right)\right)\\
&\geq \mathbb{E}\left(\exp\left(N\sum_{k=1}^3\mu^N_{0,k}(f_k)\right)\mathcal{U}^N_{F,G,H}\left(T, \xi^N\right)1_{\{\mu^N\in C\}}\right). \notag
\end{align}
By Assumption (A),
\begin{align}\label{equ 4.1.3}
&\mathbb{E}\exp\left(N\sum_{k=1}^3\mu^N_{0,k}(f_k)\right) \notag\\
&=\prod_{i=1}^N\mathbb{E}\exp\left(S_0^N(i)f_1\left(\frac{i}{N}\right)
+\mathcal{E}_0^N(i)f_2\left(\frac{i}{N}\right)+I_0^N(i)f_3\left(\frac{i}{N}\right)\right) \\
&=\prod_{i=1}^N\left(\sum_{k=1}^3e^{f_k\left(\frac{i}{N}\right)}\rho_{k-1}\left(\frac{i}{N}\right)+1-\sum_{k=1}^3\rho_{k-1}\left(\frac{i}{N}\right)\right) \notag\\
&=\exp\left(N\left(\int_\mathbb{T}\log\left(1+\sum_{k=1}^3\rho_{k-1}(u)\left(\exp\left(f_k(u)\right)-1\right)\right)du+o(1)\right)\right).\notag
\end{align}
By Equations \eqref{equ 4.1.1}, \eqref{equ 4.1.2} and \eqref{equ 4.1.3},
\begin{align*}
&\exp\left(N\left(\int_\mathbb{T}\log\left(1+\sum_{k=1}^3\rho_{k-1}(u)\left(\exp\left(f_k(u)\right)-1\right)\right)du+o(1)\right)\right)\\
&\geq P(\mu^N\in C)\times \\
&\text{\quad}\exp\left(N\left(\inf_{W\in C}\left(\sum_{k=1}^3W_{0,k}(f_k)+\mathcal{I}_1(W,F,G,H)\right)\right)\right).
\end{align*}
As a result,
\begin{align*}
&\limsup_{N\rightarrow+\infty}\frac{1}{N}\log P(\mu^N\in C)\\
&\leq -\inf_{W\in C}\left(\sum_{k=1}^3W_{0,k}(f_k)+\mathcal{I}_1(W,F,G,H)\right)\\
&\text{\quad}+\int_\mathbb{T}\log\left(1+\sum_{k=1}^3\rho_{k-1}(u)\left(\exp\left(f_k(u)\right)-1\right)\right)du\\
&=-\inf_{W\in C}\left(\mathcal{I}_2(W_0, f_1, f_2, f_3)+\mathcal{I}_1(W,F,G,H)\right).
\end{align*}
Since $f_1, f_2, f_3, F, G, H$ are arbitrary,
\begin{align}\label{equ 4.1.4}
&\limsup_{N\rightarrow+\infty}\frac{1}{N}\log P(\mu^N\in C)\\
&\leq -\sup_{f_1,f_2,f_3,F,G,H}\inf_{W\in C}\left(\mathcal{I}_2(W_0, f_1, f_2, f_3)+\mathcal{I}_1(W,F,G,H)\right). \notag
\end{align}
Since $\mathcal{I}_2(W_0, f_1, f_2, f_3)+\mathcal{I}_1(W,F,G,H)$ is convex of $W$ and concave of
\[
(f_1, f_2, f_3, F, G, H),
\]
according to the minimax theorem given in \cite{Sion1958},
\begin{align*}
&\sup_{f_1,f_2,f_3,F,G,H}\inf_{W\in C}\left(\mathcal{I}_2(W_0, f_1, f_2, f_3)+\mathcal{I}_1(W,F,G,H)\right)\\
&=\inf_{W\in C}\sup_{f_1,f_2,f_3,F,G,H}\left(\mathcal{I}_2(W_0, f_1, f_2, f_3)+\mathcal{I}_1(W,F,G,H)\right)
\end{align*}
when $C$ is compact. Since $(f_1, f_2, f_3)$ and $(F,G,H)$ do not have common variable,
\begin{align*}
&\sup_{f_1,f_2,f_3,F,G,H}\left(\mathcal{I}_2(W_0, f_1, f_2, f_3)+\mathcal{I}_1(W,F,G,H)\right)\\
&=\sup_{f_1, f_2, f_3}\left(\mathcal{I}_2(W_0, f_1, f_2, f_3)\right)+\sup_{F, G, H}\left(\mathcal{I}_1(W,F,G,H)\right)\\
&=I_{ini}(W_0)+I_{dyn}(W)
\end{align*}
and then Lemma \ref{lemma 4.1.1 LDPupperboundforCompact} follows from Equation \eqref{equ 4.1.4}.

\qed

At last, we prove Equation \eqref{equ 2.5 LDP upperbound}.

\proof[Proof of Equation \eqref{equ 2.5 LDP upperbound}]

By Lemma \ref{lemma 4.1.1 LDPupperboundforCompact}, we only need to check the exponential tightness of $\{\mu^N\}_{N\geq 1}$. According to the criterion given in \cite{Puhalskii1994}, we only need to check that
\begin{equation}\label{equ 4.1.5}
\limsup_{\delta\rightarrow 0}\limsup_{N\rightarrow+\infty}\frac{1}{N}
\log \sup_{\sigma\in \mathcal{T}}P\left(\sup_{0\leq t\leq \delta}\left|\mu_{t+\sigma,k}^N(f)-\mu_{\sigma,k}^N(f)\right|>\epsilon\right)=-\infty
\end{equation}
for any $\epsilon>0$, $f\in C(\mathbb{T})$ and $k=1,2,3$, where $\mathcal{T}$ is the set of stopping times of $\{\xi_t^N\}_{t\geq 0}$ bounded by $T$.

Since $0\leq S_t^N(i), \mathcal{E}_t^N(i), I_t^N(i)\leq 1$, $\left|\mu_{s+\sigma,k}^N(f)-\mu_{\sigma,k}^N(f)\right|$ is stochastically dominated from above by $\frac{\|f\|_\infty}{N}Y(NK_2s)$, where $\{Y(t)\}_{t\geq 0}$ is a Poisson process with rate $1$ and
\[
K_2=\|\lambda\|_\infty+\|\psi\|_\infty+\|\phi\|_\infty.
\]
Then by Markov's inequality, for any $K>0$,
\begin{align*}
&P\left(\sup_{0\leq t\leq \delta}\left|\mu_{t+\sigma,k}^N(f)-\mu_{\sigma,k}^N(f)\right|>\epsilon\right)\\
&\leq P\left(\|f\|_\infty Y(NK_2\delta)>N\epsilon\right)\leq e^{-KN\epsilon}\mathbb{E}e^{K\|f\|_\infty Y(NK_2\delta)}\\
&=e^{-KN\epsilon}e^{NK_2\delta\left(e^{K\|f\|_\infty}-1\right)}.
\end{align*}
As a result,
\[
\limsup_{\delta\rightarrow 0}\limsup_{N\rightarrow+\infty}\frac{1}{N}
\log \sup_{\sigma\in \mathcal{T}}P\left(\sup_{0\leq t\leq \delta}\left|\mu_{t+\sigma,k}^N(f)-\mu_{\sigma,k}^N(f)\right|>\epsilon\right)\leq -K\epsilon.
\]
Since $K$ is arbitrary, let $K\rightarrow+\infty$ in the above equation and then Equation \eqref{equ 4.1.5} holds, which completes the proof.

\qed

\subsection{Proof of Equation \eqref{equ 2.6 LDP lowerbound}}\label{subsection 4.2}
In this subsection we prove Equation \eqref{equ 2.6 LDP lowerbound}. We first introduce some notations and definitions for later use. For $f_1, f_2, f_3\in C(\mathbb{T})$, if
\[
0<f_1(u)<f_1(u)+f_2(u)<f_1(u)+f_2(u)+f_3(u)<1
\]
for all $u\in \mathbb{T}$, then we say that $(f_1, f_2, f_3)$ is a reasonable team. For a reasonable team $(f_1, f_2, f_3)$, we denote by $P^N_{f_1, f_2, f_3}$ the probability measure of our process $\{\xi_t^N\}_{t\geq 0}$ with the initial condition where
$\{\xi_0^N(i)\}_{0\leq i\leq N-1}$ are independent and $\xi_0^N(i)$ takes $0,1,2,3$ with respective probabilities
\[
f_1\left(\frac{i}{N}\right), \text{~}f_2\left(\frac{i}{N}\right),\text{~}f_3\left(\frac{i}{N}\right)\text{~and~}1-\sum_{k=1}^3f_k\left(\frac{i}{N}\right)
\]
for all $0\leq i\leq N-1$. For any $F, G, H\in C^{1,0}\left([0, T]\times\mathbb{T}\right)$, we denote by $\hat{P}^{N,F,G,H}_{f_1, f_2, f_3}$ the probability measure such that
\[
\frac{d\hat{P}^{N,F,G,H}_{f_1, f_2, f_3}}{dP^N_{f_1, f_2, f_3}}=\mathcal{U}^N_{F,G,H}\left(T, \xi^N\right),
\]
which is well-defined since $\{\mathcal{U}^N_{F,G,H}\left(t, \xi^N\right)\}_{0\leq t\leq T}$ is an exponential martingale with mean one. For any $0\leq t\leq T$, $f\in C(\mathbb{T})$ and $k=1,2,3$, we define
\[
\mathcal{M}^N_{f,k}(t)=\mu_{t,k}^N(f)-\mu_{0,k}^N(f)-\int_0^t\mathcal{L}^N\mu^N_{s,k}(f)ds.
\]
According to Dynkin's martingale formula, $\{\mathcal{M}^N_{f,k}(t)\}_{t\geq 0}$ is a martingale. Furthermore, the cross variation process
$\langle \mathcal{M}_{f, k_1}^N, \mathcal{M}_{g, k_2}^N \rangle(t)$ of $\mathcal{M}^N_{f,k_1}(t)$ and $\mathcal{M}^N_{g,k_2}(t)$ is given by
\begin{align*}
&\langle \mathcal{M}_{f, k_1}^N, \mathcal{M}_{g, k_2}^N \rangle(t)=\\
&\int_0^t\mathcal{L}^N\left(\mu_{s,k_1}^N(f)\mu_{s,k_2}^N(g)\right)-\mu_{s,k_1}^N(f)\mathcal{L}^N\left(\mu_{s,k_2}^N(g)\right)
-\mu_{s,k_2}^N(g)\mathcal{L}^N\left(\mu_{s,k_1}^N(f)\right) ds.
\end{align*}
For a reasonable team $(f_1, f_2, f_3)$ and $F, G, H\in C^{1,0}\left([0, T]\times \mathbb{T}\right)$, we denote by $\mu_{f_1, f_2, f_3}^{F, G, H}$ the element in $D_0$ such that $\mu_{f_1, f_2, f_3, t, k}^{F, G, H}(du)=\vartheta_{f_1, f_2, f_3, t, k}^{F, G, H}(u)du$, where $\{\vartheta_{f_1, f_2, f_3, t, k}^{F, G, H}(u)\}_{k=1,2,3}$ is the solution to the following ordinary differential equation.
\begin{equation}\label{equ 4.2.1 LLNODEunderTransformedMeasure}
\begin{cases}
&\frac{d}{dt}\vartheta_{t,1}(u)=-\vartheta_{t,1}(u)\int_\mathbb{T}\lambda(u,v)\vartheta_{t,3}(v)dv\exp\left(-F(t,u)+G(t,u)\right),\\
&\frac{d}{dt}\sum_{k=1}^2\vartheta_{t,k}(u)=-\psi(u)\vartheta_{t,2}(u)\exp\left(-G(t,u)+H(t,u)\right),\\
&\frac{d}{dt}\sum_{k=1}^3\vartheta_{t,k}(u)=-\phi(u)\exp\left(-H(t,u)\right),\\
&\vartheta_{0,k}(u)=f_k(u)\text{~for~}k=1,2,3,\\
&\left(\vartheta_{t,1}, \vartheta_{t,2}, \vartheta_{t,3}\right)\text{~is reasonable for all~}0\leq t\leq T.
\end{cases}
\end{equation}
By Grownwall's inequality, it is not difficult to check that the solution to Equation \eqref{equ 4.2.1 LLNODEunderTransformedMeasure} is unique. The motivation of the definition of $\mu_{f_1, f_2, f_3}^{F, G, H}$ is to find a $W\in D_0$ such that $F_W=F, G_W=G, H_W=H$ and $W_{0,k}(du)=f_k(u)du$ for any given $f_1, f_2, f_3, F, G, H$. Then by Equation \eqref{equ 3.2}, $W$ should be with densities as the solution to equation \eqref{equ 4.2.1 LLNODEunderTransformedMeasure}.

The proof of Equation \eqref{equ 2.6 LDP lowerbound} relies heavily on the following Lemma.

\begin{lemma}\label{lemma 4.2.1 LLNundertransformedMeasure}
Let $(f_1, f_2, f_3)$ be a reasonable team and $F, G, H\in C^{1,0}\left([0, T]\times \mathbb{T}\right)$.
As $N\rightarrow+\infty$, $\mu^N$ converges in $\hat{P}^{N,F,G,H}_{f_1, f_2, f_3}$-probability to $\mu_{f_1, f_2, f_3}^{F, G, H}$.
\end{lemma}

We prove Lemma \ref{lemma 4.2.1 LLNundertransformedMeasure} later. Now we utilize Lemma \ref{lemma 4.2.1 LLNundertransformedMeasure} to prove Equation \eqref{equ 2.6 LDP lowerbound}.

\proof[Proof of Equation \eqref{equ 2.6 LDP lowerbound}]

For open set $O\subseteq \Omega_2$, if $\inf_{W\in O\cap D_0}\left(I_{ini}(W_0)+I_{dyn}(W)\right)=+\infty$, then Equation \eqref{equ 2.6 LDP lowerbound} is trivial. If $\inf_{W\in O\cap D_0}\left(I_{ini}(W_0)+I_{dyn}(W)\right)<+\infty$, then for any $\epsilon>0$, there exists $\hat{W}\in O\cap D_0$ such that
\[
I_{ini}(\hat{W}_0)+I_{dyn}(\hat{W})\leq \inf_{W\in O\cap D_0}\left(I_{ini}(W_0)+I_{dyn}(W)\right)+\epsilon.
\]
Since $\hat{W}\in D_0$, $\hat{W}=\mu_{f_1, f_2, f_3}^{F, G, H}$ with $F=F_{\hat{W}}, G=G_{\hat{W}}$, $H=H_{\hat{W}}$ and $f_k=\frac{d\hat{W}_{0,k}}{du}$.
Then, by Lemma \ref{lemma 2.2 rate function on D0},
\begin{equation}\label{equ 4.2.2}
I_{dyn}(\hat{W})=\mathcal{I}_1(\hat{W}, F, G, H)
\end{equation}
and
\begin{align}\label{equ 4.2.3}
I_{ini}(\hat{W}_0)=&\int_{\mathbb{T}}\Bigg(\sum_{k=1}^3f_k(u)\log\left(\frac{f_k(u)}{\rho_{k-1}(u)}\right)\\
&\text{\quad}-\left(1-\sum_{k=1}^3f_k(u)\right)\log\left(\frac{1-\sum_{k=1}^3\rho_{k-1}(u)}{1-\sum_{k=1}^3f_k(u)}\right)\Bigg)du. \notag
\end{align}
According to our definition of $\hat{P}^{N,F,G,H}_{f_1, f_2, f_3}$,
\begin{equation}\label{equ 4.2.4}
P\left(\mu^N\in O\right)=\mathbb{E}_{\hat{P}^{N,F,G,H}_{f_1, f_2, f_3}}\left(\frac{dP}{dP^N_{f_1, f_2, f_3}}\frac{1}{\mathcal{U}^N_{F,G,H}\left(T, \xi^N\right)}1_{\{\mu^N\in O\}}\right).
\end{equation}
We denote by $D_{\epsilon,1}$ the event that
\[
\frac{dP}{dP^N_{f_1, f_2, f_3}}\geq \exp\left(N(-I_{ini}(\hat{W}_0)-\epsilon)\right)
\]
and denote by $D_{\epsilon, 2}$ the event that
\[
\left|\mathcal{I}_1(\mu^N, F, G, H)-\mathcal{I}_1(\hat{W}, F, G, H)\right|\leq \epsilon,
\]
then by Equations \eqref{equ 4.1.1}, \eqref{equ 4.2.2} and \eqref{equ 4.2.4},
\begin{align}\label{equ 4.2.5}
P\left(\mu^N\in O\right)\geq &\hat{P}^{N,F,G,H}_{f_1, f_2, f_3}\left(\{\mu^N\in O\}\cap D_{\epsilon,1}\cap D_{\epsilon,2}\right)\\
&\times \exp\left(-N\left(I_{ini}(\hat{W}_0)+I_{dyn}(\hat{W})+2\epsilon\right)\right).\notag
\end{align}
Since $\hat{W}=\mu_{f_1, f_2, f_3}^{F, G, H}$, $\mu^N$ converges in $\hat{P}^{N,F,G,H}_{f_1, f_2, f_3}$-probability to $\hat{W}$ as $N\rightarrow+\infty$ by Lemma \ref{lemma 4.2.1 LLNundertransformedMeasure}. Then, since $\hat{W}\in O$ and $\mathcal{I}_1(W, F, G, H)$ is continuous in $W$,
\[
\lim_{N\rightarrow+\infty}\hat{P}^{N,F,G,H}_{f_1, f_2, f_3}(\{\mu^N\in O\}\cap D_{\epsilon, 2})=1.
\]
According to the definition of $P^N_{f_1, f_2, f_3}$,
\begin{align*}
&\frac{1}{N}\log\frac{dP}{dP^N_{f_1, f_2, f_3}}=\frac{1}{N}\sum_{i=0}^{N-1}
\Bigg(\sum_{k=1}^31_{\{\xi_0^N(i)=k-1\}}\log\left(\frac{\rho_{k-1}\left(\frac{i}{N}\right)}{f_k\left(\frac{i}{N}\right)}\right)\\
&\text{\quad\quad\quad}+1_{\{\xi_0^N(i)=3\}}\log\left(\frac{1-\sum_{k=1}^3\rho_{k-1}\left(\frac{i}{N}\right)}{1-\sum_{k=1}^3f_k\left(\frac{i}{N}\right)}\right)\Bigg).
\end{align*}
Then, according to a mean-variance analysis, it is easy to check that
\begin{align*}
&\lim_{N\rightarrow+\infty}\frac{1}{N}\log\frac{dP}{dP^N_{f_1, f_2, f_3}}=\int_{\mathbb{T}}\Bigg(\sum_{k=1}^3f_k(u)\log\left(\frac{\rho_{k-1}(u)}{f_k(u)}\right)\\
&\text{\quad}+\left(1-\sum_{k=1}^3f_k(u)\right)\log\left(\frac{1-\sum_{k=1}^3\rho_{k-1}(u)}{1-\sum_{k=1}^3f_k(u)}\right)\Bigg)du=-I_{ini}(\hat{W}_0)
\end{align*}
in $\hat{P}^{N,F,G,H}_{f_1, f_2, f_3}$-probability by Equation \eqref{equ 4.2.3}. As a result,
\[
\lim_{N\rightarrow+\infty}\hat{P}^{N,F,G,H}_{f_1, f_2, f_3}(D_{\epsilon, 1})=1
\]
and hence
\begin{equation}\label{equ 4.2.6}
\lim_{N\rightarrow+\infty}\hat{P}^{N,F,G,H}_{f_1, f_2, f_3}\left(\{\mu^N\in O\}\cap D_{\epsilon,1}\cap D_{\epsilon,2}\right)=1.
\end{equation}
Then, by Equations \eqref{equ 4.2.5} and \eqref{equ 4.2.6},
\begin{align*}
\liminf_{N\rightarrow+\infty}\frac{1}{N}\log P(\mu^N\in O)&\geq-\left(I_{ini}(\hat{W}_0)+I_{dyn}(\hat{W})\right)-2\epsilon\\
&\geq -\inf_{W\in O\cap D_0}\left(I_{ini}(W_0)+I_{dyn}(W)\right)-3\epsilon.
\end{align*}
Since $\epsilon$ is arbitrary in the above equation, let $\epsilon\rightarrow 0$ and then the proof is complete.

\qed

Now we only need to prove Lemma \ref{lemma 4.2.1 LLNundertransformedMeasure}. The following three lemmas play key roles in the proof of Lemma \ref{lemma 4.2.1 LLNundertransformedMeasure}.

\begin{lemma}\label{lemma 4.2.2}
Let $N\rightarrow+\infty$, then
\[
\sum_{0\leq t\leq T}\left(\mu^N_{t,k}(f)-\mu^N_{t-,k}(f)\right)^2
\]
converges both in $P$-probability and $\hat{P}^{N,F,G,H}_{f_1, f_2, f_3}$-probability to $0$ for any $f\in C(\mathbb{T})$ and $k=1,2,3$, where $\mu^N_{t-,k}$ is the state of $\mu^N_{\cdot, k}$ at the moment just before $t$.

\end{lemma}

\begin{lemma}\label{lemma 4.2.3}
The sequence $\{\mu^N\}_{N\geq 1}$ are $\hat{P}^{N,F,G,H}_{f_1, f_2, f_3}$-tight.
\end{lemma}

\begin{lemma}\label{lemma 4.2.4}
Let $(f_1, f_2, f_3)$ be a reasonable team and $F, G, H\in C^{1,0}\left([0, T]\times \mathbb{T}\right)$, then the solution $W$ to the following $\Omega_2$-valued integral equation is unique. For any $f,g,h\in C(\mathbb{T})$ and $t\geq 0$,
\begin{equation}\label{equ integral equation}
\begin{cases}
&W_{t,1}(f)=W_{0,1}(f)-\int_0^tW_{s,1}\bigotimes W_{s,3}\left(\lambda(\cdot, \ast)f(\cdot)\left(e^{-F_s(\cdot)+G_s(\cdot)}\right)\right)ds,\\
&W_{t,2}(g)=W_{0,2}(g)+\int_0^tW_{s,1}\bigotimes W_{s,3}\left(\lambda(\cdot, \ast)g(\cdot)\left(e^{-F_s(\cdot)+G_s(\cdot)}\right)\right)ds\\
&\text{\quad\quad\quad\quad\quad}-\int_0^tW_{s,2}\left(\psi(\cdot)g(\cdot)\left(e^{-G_s(\cdot)+H_s(\cdot)}\right)\right)ds,\\
&W_{t,3}(h)=W_{0,3}(h)+\int_0^tW_{s,2}\left(\psi(\cdot)h(\cdot)\left(e^{-G_s(\cdot)+H_s(\cdot)}\right)\right)ds\\
&\text{\quad\quad\quad\quad\quad}-\int_0^tW_{s,3}\left(\phi(\cdot)h(\cdot)\left(e^{-H_s(\cdot)}\right)\right)ds,\\
&W_{0,k}(du)=f_{k}(u)du\text{~for~}k=1,2,3.
\end{cases}
\end{equation}
\end{lemma}

\proof[Proof of Lemma \ref{lemma 4.2.2}]

Under probability measure $P$, since $\mu^N_{t,k}(f)$ changes by at most $\frac{\|f\|_\infty}{N}$ at each jump moment $t$, $\sum_{0\leq t\leq T}\left(\mu^N_{t,k}(f)-\mu^N_{t-,k}(f)\right)^2$ is stochastically dominated from above by $\frac{\|f\|^2_\infty}{N^2}Y(NK_2T)$, where $Y(t)$ is the Poisson process at rate one and $K_2=\|\lambda\|_\infty+\|\psi\|_\infty+\|\phi\|_\infty$ defined as in the proof of Equation \eqref{equ 2.5 LDP upperbound}. Therefore,
by Markov's inequality,
\begin{align}\label{equ 4.2.7}
P\left(\sum_{0\leq t\leq T}\left(\mu^N_{t,k}(f)-\mu^N_{t-,k}(f)\right)^2\geq \epsilon\right)&\leq e^{-N^2\epsilon}\mathbb{E}\exp\left(\|f\|^2_\infty Y(NK_2 T)\right)\notag\\
&=e^{-N^2\epsilon}e^{NK_2T(e^{\|f\|^2_\infty}-1)}
\end{align}
for any $\epsilon>0$. As a result, $\sum_{0\leq t\leq T}\left(\mu^N_{t,k}(f)-\mu^N_{t-,k}(f)\right)^2$ converges in $P$-probability to $0$.

According to our definition of $\hat{P}^{N,F,G,H}_{f_1, f_2, f_3}$, it is easy to check that there exists $K_3>0$ independent of $N$ such that
\begin{equation}\label{equ 4.2.8}
\left(\frac{d\hat{P}^{N,F,G,H}_{f_1, f_2, f_3}}{dP}\right)^2\leq \exp(K_3N)
\end{equation}
for sufficiently large $N$. Hence, by Equation \eqref{equ 4.2.7} and Cauchy-Schwarz inequality, $\sum_{0\leq t\leq T}\left(\mu^N_{t,k}(f)-\mu^N_{t-,k}(f)\right)^2$ converges in $\hat{P}^{N,F,G,H}_{f_1, f_2, f_3}$-probability to $0$.

\qed

\proof[Proof of Lemma \ref{lemma 4.2.3}]

According to Aldous' criterion, we only need to check that
\begin{equation}\label{equ 4.2.9}
\lim_{\delta\rightarrow 0}\limsup_{N\rightarrow+\infty}\sup_{\sigma\in\mathcal{T}, t\leq \delta}\hat{P}^{N,F,G,H}_{f_1, f_2, f_3}\left(\left|\mu^N_{t+\sigma, k}(f)-\mu^N_{\sigma,k}(f)\right|>\epsilon\right)=0
\end{equation}
for any $f\in C(\mathbb{T})$, $k=1,2,3$ and $\epsilon>0$. In the proof of Equation \eqref{equ 2.5 LDP upperbound}, we have shown that
\[
P\left(\left|\mu^N_{t+\sigma, k}(f)-\mu^N_{\sigma,k}(f)\right|>\epsilon\right)
\leq \exp\left(-KN\epsilon+NK_2\delta\left(e^{K\|f\|_\infty}-1\right)\right)
\]
for any $t\leq \delta$ and $K>0$. Then, by Equation \eqref{equ 4.2.8},
\begin{align*}
&\hat{P}^{N,F,G,H}_{f_1, f_2, f_3}\left(\left|\mu^N_{t+\sigma, k}(f)-\mu^N_{\sigma,k}(f)\right|>\epsilon\right)\\
&\leq \exp\left(-KN\epsilon+NK_2\delta\left(e^{K\|f\|_\infty}-1\right)+\frac{1}{2}K_3N\right).
\end{align*}
Then for given $K>\frac{K_3}{\epsilon}$ and $\delta$ sufficiently small such that $K_2\delta\left(e^{K\|f\|_\infty}-1\right)<\frac{K_3}{4}$,
\[
\hat{P}^{N,F,G,H}_{f_1, f_2, f_3}\left(\left|\mu^N_{t+\sigma, k}(f)-\mu^N_{\sigma,k}(f)\right|>\epsilon\right)\leq \exp\left(N(-K\epsilon+K_3)\right)
\]
for $t\leq \delta$ and hence Equation \eqref{equ 4.2.9} holds.

\qed

\proof[Proof of Lemma \ref{lemma 4.2.4}]
Let $W$ and $\tilde{W}$ be two solutions to Equation \eqref{equ integral equation}, then we only need to show that $W=\tilde{W}$.

For any $\nu\in (C(\mathbb{T}))^\prime$, let
\[
\|\nu\|_\infty=\sup_{f\in C(\mathbb{T}):~\|f\|_\infty=1}|\nu(f)|.
\]
Then, for $\nu_1, \nu_2, \nu_3\in \mathbb{M}$,
\[
\left|\nu_1\bigotimes \nu_3\left(h(\cdot, \ast)\right)-\nu_2\bigotimes \nu_3\left(h(\cdot, \ast)\right)\right|
\leq \|\nu_1-\nu_2\|_\infty\|h\|_\infty
\]
and
\[
\left|\nu_1\bigotimes \nu_2\left(h(\cdot, \ast)\right)-\nu_1\bigotimes \nu_3\left(h(\cdot, \ast)\right)\right|
\leq \|\nu_2-\nu_3\|_\infty\|h\|_\infty
\]
for any $h\in C(\mathbb{T}^2)$. Then, since $W$ and $\tilde{W}$ are both solutions to Equation \eqref{equ integral equation},
\[
\|W_{t,1}-\tilde{W}_{t,1}\|_\infty\leq\|\lambda\|_\infty\exp\left(\|F\|_\infty+\|G\|_\infty\right)
\int_0^t\|W_{s,1}-\tilde{W}_{s,1}\|_\infty+\|W_{s,3}-\tilde{W}_{s,3}\|_\infty ds,
\]
\begin{align*}
&\|W_{t,2}-\tilde{W}_{t,2}\|_\infty\\
&\leq \|\lambda\|_\infty\exp\left(\|F\|_\infty+\|G\|_\infty\right)
\int_0^t\|W_{s,1}-\tilde{W}_{s,1}\|_\infty+\|W_{s,3}-\tilde{W}_{s,3}\|_\infty ds\\
&\text{\quad}+\|\psi\|_\infty\exp\left(\|H\|_\infty+\|G\|_\infty\right)\int_0^t\|W_{s,2}-\tilde{W}_{s,2}\|_\infty ds
\end{align*}
and
\begin{align*}
\|W_{t,3}-\tilde{W}_{t,3}\|_\infty\leq &\|\psi\|_\infty\exp\left(\|H\|_\infty+\|G\|_\infty\right)\int_0^t\|W_{s,2}-\tilde{W}_{s,2}\|_\infty ds\\
&+\|\phi\|_\infty\exp\left(\|H\|_\infty\right)\int_0^t\|W_{s,3}-\tilde{W}_{s,3}\|_\infty ds
\end{align*}
for any $t\geq 0$. As a result,
\[
\sum_{k=1}^3\|W_{t,k}-\tilde{W}_{t,k}\|_\infty
\leq K_1\int_0^t\sum_{k=1}^3\|W_{s,k}-\tilde{W}_{s,k}\|_\infty ds
\]
for any $t\geq 0$, where
\[
K_1=2\left(\|\lambda\|_\infty+\|\psi\|_\infty+\|\phi\|_\infty\right)\exp\left(\|F\|_\infty+\|G\|_\infty+\|H\|_\infty\right).
\]
Then, according to Grownwall's inequality,
\[
\sum_{k=1}^3\|W_{t,k}-\tilde{W}_{t,k}\|_\infty\leq 0\exp(K_1t)=0
\]
and hence $W=\tilde{W}$, which completes the proof.

\qed

At last, we prove Lemma \ref{lemma 4.2.1 LLNundertransformedMeasure}.

\proof[Proof of Lemma \ref{lemma 4.2.1 LLNundertransformedMeasure}]

Let $\hat{\mu}$ be a $\hat{P}^{N,F,G,H}_{f_1, f_2, f_3}$-weak limit of a subsequence of $\{\mu^N\}_{N\geq 1}$, then we only need to show that $\hat{\mu}=\mu_{f_1, f_2, f_3}^{F, G, H}$. Note that the existence of $\hat{\mu}$ follows from Lemma \ref{lemma 4.2.3}. For simplicity, we still denote the subsequence converging to $\hat{\mu}$ by $\{\mu^N\}_{N\geq 1}$.

According to It\^{o}'s formula,
\[
d\mathcal{U}^N_{F,G,H}(t, \xi^N)
=\frac{1}{V^N_{F,G,H}(0, \xi^N)}\exp\left(-\int_0^t
\frac{\left(\partial_s+\mathcal{L}^N\right)V_{F,G,H}^N(s, \xi^N)}{V_{F,G,H}^N(s, \xi^N)}ds\right)d\Lambda^N_{F, G, H}(t),
\]
where
\[
\Lambda^N_{F,G,H}(t)=V_{F,G,H}^N(t,\xi^N)-V_{F,G,H}^N(0,\xi^N)-\int_0^t(\partial_s+\mathcal{L}^N)V_{F,G,H}^N(s, \xi^N)ds
\]
is a martingale. Hence,
\begin{equation}\label{equ 4.2.10}
d\mathcal{U}^N_{F,G,H}(t, \xi^N)=\mathcal{U}^N_{F,G,H}(t, \xi^N)d\hat{\Lambda}^N_{F,G,H}(t),
\end{equation}
where
\[
d\hat{\Lambda}^N_{F,G,H}(t)=\frac{1}{V_{F,G,H}^N(t, \xi^N)}d\Lambda^N_{F, G, H}(t).
\]
By Equation \eqref{equ 4.2.10} and a generalized version of Girsanov's theorem introduced in \cite{Schuppen1974}, for any martingale $\{M_t\}_{0\leq t\leq T}$ under $P^N_{f_1, f_2, f_3}$, let
\[
\hat{M}_t=M_t-\langle M, \hat{\Lambda}^N_{F, G, H}\rangle_t,
\]
then $\{\hat{M}_t\}_{0\leq t\leq T}$ is a martingale under $\hat{P}^{N, F, G, H}_{f_1, f_2, f_3}$. Furthermore,
\[
\left[M, M\right]_t=\left[\hat{M}, \hat{M}\right]_t
\]
under both probability measures, where we use $[\cdot, \cdot]$ to denote the quadratic variation process. As a result, let
\[
\hat{\mathcal{M}}^N_{f,k}(t)=\mathcal{M}^N_{f,k}(t)-\langle \mathcal{M}^N_{f,k}, \hat{\Lambda}^N_{F, G, H}\rangle_t,
\]
then $\{\hat{\mathcal{M}}^N_{f,k}(t)\}_{0\leq t\leq T}$ is a martingale under $\hat{P}^{N, F, G, H}_{f_1, f_2, f_3}$ with
\[
\left[\hat{\mathcal{M}}^N_{f,k}, \hat{\mathcal{M}}^N_{f,k}\right]_t=\sum_{0\leq s\leq t}\left(\mu_{s,k}^N(f)-\mu_{s-,k}^N(f)\right)^2.
\]
Then by Lemma \ref{lemma 4.2.2} and Doob's inequality, as $N\rightarrow+\infty$,
\[
\sup_{0\leq t\leq T}\left|\hat{\mathcal{M}}^N_{f,k}(t)\right|
\]
converges in $\hat{P}^{N, F, G, H}_{f_1, f_2, f_3}$-probability to $0$. According to the definition of $\hat{\Lambda}^N_{F, G, H}(t)$,
\[
d\langle \mathcal{M}^N_{f,k}, \hat{\Lambda}^N_{F, G, H}\rangle_t
=\frac{1}{V_{F,G,H}^N(t, \xi^N)}d\langle \mathcal{M}^N_{f,k}, \Lambda^N_{F, G, H}\rangle_t.
\]
According to the definition of $\Lambda^N_{F, G, H}(t)$ and the Dynkin's martingale formula,
\begin{align*}
&d\langle \mathcal{M}^N_{f,k}, \Lambda^N_{F, G, H}\rangle_t=\Big(\mathcal{L}^N\left(\mu_{t,k}^N(f)V_{F,G,H}^N(t, \xi^N)\right)-V_{F,G,H}^N(t, \xi^N)\mathcal{L}^N\mu_{t,k}^N(f)\\
&\text{\quad\quad}-\mu_{t,k}^N(f)\mathcal{L}^NV_{F,G,H}^N(t, \xi^N)\Big)dt.
\end{align*}
As a result, by direct calculation,
\[
\langle\mathcal{M}^N_{f,1}, \hat{\Lambda}^N_{F, G, H}\rangle_t=-\int_0^t\mu_{s,1}^N\bigotimes\mu_{s,3}^N\left(\lambda(\cdot, \ast)f(\cdot)\left(e^{-F_s(\cdot)+G_s(\cdot)}-1\right)\right)ds,
\]
\begin{align*}
\langle\mathcal{M}^N_{g,2}, \hat{\Lambda}^N_{F, G, H}\rangle_t=&\int_0^t\mu_{s,1}^N\bigotimes\mu_{s,3}^N\left(\lambda(\cdot, \ast)g(\cdot)\left(e^{-F_s(\cdot)+G_s(\cdot)}-1\right)\right)ds\\
&-\int_0^t\mu_{s,2}^N\left(\psi(\cdot)g(\cdot)\left(e^{-G_s(\cdot)+H_s(\cdot)}-1\right)\right)ds
\end{align*}
and
\begin{align*}
\langle\mathcal{M}^N_{h,3}, \hat{\Lambda}^N_{F, G, H}\rangle_t=
&\int_0^t\mu_{s,2}^N\left(\psi(\cdot)h(\cdot)\left(e^{-G_s(\cdot)+H_s(\cdot)}-1\right)\right)ds\\
&-\int_0^t\mu_{s,3}^N\left(\phi(\cdot)h(\cdot)\left(e^{-H_s(\cdot)}-1\right)\right)ds
\end{align*}
for any $f, g, h\in C(\mathbb{T})$. Furthermore, by direct calculation,
\[
\mathcal{L}^N\mu_{s,1}^N(f)=-\mu_{s,1}^N\bigotimes\mu_{s,3}^N\left(\lambda(\cdot, \ast)f(\cdot)\right),
\]
\[
\mathcal{L}^N\mu_{s,2}^N(g)=\mu_{s,1}^N\bigotimes\mu_{s,3}^N\left(\lambda(\cdot, \ast)g(\cdot)\right)-\mu_{s,2}^N(\psi(\cdot)g(\cdot)),
\]
and
\[
\mathcal{L}^N\mu_{s,3}^N(h)=\mu_{s,2}^N(\psi(\cdot)h(\cdot))-\mu_{s,3}^N(\phi(\cdot)h(\cdot)).
\]
Consequently, since
\[
\mu_{t,k}^N(f)=\hat{\mathcal{M}}^N_{f,k}(t)+\mu_{0,k}^N(f)+\int_0^t\mathcal{L}^N\mu_{s,k}^N(f)ds+\langle\mathcal{M}^N_{f,k}, \hat{\Lambda}^N_{F, G, H} \rangle_t
\]
and
$\sup_{0\leq t\leq T}\left|\hat{\mathcal{M}}^N_{f,k}(t)\right|$ converges to $0$ in probability, let $N\rightarrow+\infty$ and then
\begin{equation}\label{equ 4.2.11}
\begin{cases}
&\hat{\mu}_{t,1}(f)=\hat{\mu}_{0,1}(f)-\int_0^t\hat{\mu}_{s,1}\bigotimes\hat{\mu}_{s,3}\left(\lambda(\cdot, \ast)f(\cdot)\left(e^{-F_s(\cdot)+G_s(\cdot)}\right)\right)ds,\\
&\hat{\mu}_{t,2}(g)=\hat{\mu}_{0,2}(g)+\int_0^t\hat{\mu}_{s,1}\bigotimes\hat{\mu}_{s,3}\left(\lambda(\cdot, \ast)g(\cdot)\left(e^{-F_s(\cdot)+G_s(\cdot)}\right)\right)ds,\\
&\text{\quad\quad\quad\quad\quad}-\int_0^t\hat{\mu}_{s,2}\left(\psi(\cdot)g(\cdot)\left(e^{-G_s(\cdot)+H_s(\cdot)}\right)\right)ds\\
&\hat{\mu}_{t,3}(h)=\hat{\mu}_{0,3}(h)+\int_0^t\hat{\mu}_{s,2}\left(\psi(\cdot)h(\cdot)\left(e^{-G_s(\cdot)+H_s(\cdot)}\right)\right)ds\\
&\text{\quad\quad\quad\quad\quad}-\int_0^t\hat{\mu}_{s,3}\left(\phi(\cdot)h(\cdot)\left(e^{-H_s(\cdot)}\right)\right)ds
\end{cases}
\end{equation}
for any $f,g,h\in C(\mathbb{T})$ and $t\geq 0$. According to the definition of $\hat{P}^{N,F,G,H}_{f_1, f_2, f_3}$ and Chebyshev's inequality, it is easy to check that $\hat{\mu}_{0,k}(du)=f_k(u)du$ for $k=1,2,3$. Hence, by Equation \eqref{equ 4.2.11}, $\hat{\mu}$ is a solution to Equation \eqref{equ integral equation}. It is easy to check that $\mu_{f_1, f_2, f_3}^{F, G, H}$ is also a solution to Equation \eqref{equ integral equation}. Therefore, $\hat{\mu}=\mu_{f_1, f_2, f_3}^{F, G, H}$ follows from Lemma \ref{lemma 4.2.4} and the proof is complete.

\qed

\subsection{Proof of Theorem \ref{Theorem 2.4 exponentially rapid to hydrodynamic limit}}\label{subsection 4.3}
In this subsection we prove Theorem \ref{Theorem 2.4 exponentially rapid to hydrodynamic limit}. We first give a lemma about when our rate function equals zero. \begin{lemma}\label{lemma 4.3.1}
If $I_{ini}(W_0)+I_{dyn}(W)=0$, then $W=\mu$.
\end{lemma}
We prove Lemma \ref{lemma 4.3.1} at the end of this subsection. Now we utilize Lemma \ref{lemma 4.3.1} to prove Theorem \ref{Theorem 2.4 exponentially rapid to hydrodynamic limit}.

\proof[Proof of Theorem \ref{Theorem 2.4 exponentially rapid to hydrodynamic limit}]

Since the complementary set of an open set $O$ is a closed set, to prove Theorem \ref{Theorem 2.4 exponentially rapid to hydrodynamic limit} we only need to show that $\inf_{W\in C}\left(I_{ini}(W_0)+I_{dyn}(W)\right)>0$ for any closed set $C$ such that $\mu\not\in C$ according to Equation \eqref{equ 2.5 LDP upperbound}. Equivalently, we only need to show that $\inf_{W\in C}\left(I_{ini}(W_0)+I_{dyn}(W)\right)=0$ for some closed set $C$ implies that $\mu\in C$. If a closed set $C$ makes $\inf_{W\in C}\left(I_{ini}(W_0)+I_{dyn}(W)\right)=0$, then there exist a series $\{W^n\}_{n\geq 1}$ in $C$ such that
\begin{equation}\label{equ 4.3.1}
\lim_{n\rightarrow+\infty}I_{ini}(W^n_0)+I_{dyn}(W^n)=0.
\end{equation}
Then, there exists $0<a<+\infty$ such that $\sup_{n\geq 1}\left\{I_{ini}(W^n_0)+I_{dyn}(W^n)\right\}<a$ and hence $W^n\in \mathcal{C}_a$ for all $n\geq 1$. By Lemma \ref{lemma 2.1 good LDP rate}, $\mathcal{C}_a$ is compact. Hence, a subsequence of $W^n$ converges to some $W\in \mathcal{C}_a\cap C$. Since our rate function is lower-semicontinuous,
\[
I_{ini}(W_0)+I_{dyn}(W)\leq 0
\]
by Equation \eqref{equ 4.3.1}. Since our rate function is non-negative,
\[
I_{ini}(W_0)+I_{dyn}(W)=0.
\]
Then, $W=\mu$ by Lemma \ref{lemma 4.3.1} and hence $\mu\in C$, which completes the proof.

\qed

At last we prove Lemma \ref{lemma 4.3.1}.

\proof[Proof of Lemma \ref{lemma 4.3.1}]

If $I_{ini}(W_0)+I_{dyn}(W)=0$, then $I_{ini}(W_0)=0$ and $I_{dyn}(W)=0$. For any given $f_1,f_2,f_3\in C(\mathbb{T})$ and $h_1, h_2, h_3\in C^1([0, T])$, let $F,G,H\in C^{1,0}\left([0,T]\times \mathbb{T}\right)$ such that
\[
F(s,u)=h_1(s)f_1(u), \text{~}G(s,u)=h_2(s)f_2(u)\text{~and~}H(s,u)=h_3(s)f_3(u)
\]
for all $(s,u)\in[0, T]\times \mathbb{T}$. Since $I_{dyn}(W)=0$, $\mathcal{I}_1(W, \epsilon F,0,0)$ gets maximum at $\epsilon=0$ and hence
\[
\frac{d}{d\epsilon}\mathcal{I}_1(W, \epsilon F,0,0)\Big|_{\epsilon=0}=0.
\]
Then, by direct calculation,
\begin{align*}
&h_1(T)W_{T,1}(f_1)-h_1(0)W_{0,1}(f_1)-\int_0^T\partial_sh_1(s)W_{s,1}(f_1)ds\\
&=-\int_0^Th_1(s)W_{s,1}\bigotimes W_{s,3}(\lambda(\cdot, \ast)f_1(\cdot))ds.
\end{align*}
Since $h_1$ in the above equation is arbitrary, $\{W_{s,1}(f_1)\}_{0\leq s\leq T}$ is absolutely continuous and
\[
\frac{d}{dt}W_{t,1}(f_1)=-W_{t,1}\bigotimes W_{t,3}(\lambda(\cdot, \ast)f_1(\cdot)).
\]
According to similar calculations on $\mathcal{I}_1(W, 0, \epsilon G,0)$ and $\mathcal{I}_1(W,0,0, \epsilon H)$,
\[
\frac{d}{dt}W_{t,2}(f_2)=W_{t,1}\bigotimes W_{t,3}(\lambda(\cdot, \ast)f_2(\cdot))-W_{t,2}(\psi f_2)
\]
and
\[
\frac{d}{dt}W_{t,3}(f_3)=W_{t,2}(\psi f_3)-W_{t,3}(\phi f_3).
\]
Similarly, since $I_{ini}(W_0)=0$,
\[
\frac{d}{d\epsilon}\mathcal{I}_2(W_0, \epsilon f_1, 0,0)\Big|_{\epsilon=0}=0,\text{~}\frac{d}{d\epsilon}\mathcal{I}_2(W_0, 0, \epsilon f_2,0)\Big|_{\epsilon=0}=0
\]
and $\frac{d}{d\epsilon}\mathcal{I}_2(W_0, 0, 0, \epsilon f_3)\Big|_{\epsilon=0}=0$. Then, by direct calculation,
\[
W_{0,k}(f_k)=\int_{\mathbb{T}}\rho_{k-1}(u)f_k(u)du
\]
for $k=1,2,3$ and hence $W_{0,k}(du)=\rho_{k-1}(u)du$ for $k=1,2,3$. In conclusion, $W$ is a solution to the measure-valued ordinary differential equation
\begin{equation}\label{equ 4.3.2 measure valued ODE}
\begin{cases}
&\frac{d}{dt}W_{t,1}(f_1)=-W_{t,1}\bigotimes W_{t,3}(\lambda(\cdot, \ast)f_1(\cdot)) \text{~for any~}f_1\in C(\mathbb{T}),\\
&\frac{d}{dt}W_{t,2}(f_2)=W_{t,1}\bigotimes W_{t,3}(\lambda(\cdot, \ast)f_2(\cdot))-W_{t,2}(\psi f_2) \text{~for any~}f_2\in C(\mathbb{T}),\\
&\frac{d}{dt}W_{t,3}(f_3)=W_{t,2}(\psi f_3)-W_{t,3}(\phi f_3) \text{~for any~}f_3\in C(\mathbb{T}),\\
&W_{0,k}(du)=\rho_{k-1}(u)du \text{~for~} k=1,2,3.
\end{cases}
\end{equation}
According to the definition of $\mu$, it is easy to check that $\mu$ is a solution to Equation \eqref{equ 4.3.2 measure valued ODE}. Now we only need to show that the solution to Equation \eqref{equ 4.3.2 measure valued ODE} is unique. Note that Equation \eqref{equ 4.3.2 measure valued ODE} is a special case of Equation \eqref{equ integral equation} where $(f_1, f_2, f_3)=(\rho_0, \rho_1, \rho_2)$ and $(F, G, H)=(0,0,0)$. Hence, the uniqueness of the solution to Equation \eqref{equ 4.3.2 measure valued ODE} follows from Lemma \ref{lemma 4.2.4} and then the proof is complete.

\qed

\section{Proofs of Theorems \ref{Theorem 2.5 MDP} and \ref{Theorem 2.6 MDPofHittingTimes}}\label{section five}
In this section we prove Theorems \ref{Theorem 2.5 MDP} and \ref{Theorem 2.6 MDPofHittingTimes}. The proof of Theorem \ref{Theorem 2.5 MDP} follows the exponential-martingale strategy introduced in \cite{Gao2003}. Main new mathematical difficulties are mainly in checks of exponential tightness of $\{\eta^N\}_{N\geq 1}$ under $P $ and tightness of $\{\eta^N\}_{N\geq 1}$ under a transformed measure.

\subsection{Preliminaries}\label{subsection 5.1}
In this subsection we give some lemmas as preliminaries of the proof of Theorem \ref{Theorem 2.5 MDP}. We first introduce some notations and definitions. For any $F, G, H, \hat{F}, \hat{G}, \hat{H}\in C^{1,\infty}([0, T]\times \mathbb{T})$, we define
\begin{align*}
&\mathcal{B}_{10}\left((F,G,H), (\hat{F}, \hat{G}, \hat{F})\right)\\
&=\int_0^T\int_{\mathbb{T}^2}\lambda(u,v)\theta^S_s(u)\theta^I_s(v)\left(G_s(u)-F_s(u)\right)\left(\hat{G}_s(u)-\hat{F}_s(u)\right)dsdudv\\
&\text{\quad}+\int_0^T\int_{\mathbb{T}}\psi(u)\theta_s^E(u)\left(-G_s(u)+H_s(u)\right)\left(-\hat{G}_s(u)+\hat{H}_s(u)\right)dsdu\\
&\text{\quad}+\int_0^T\int_{\mathbb{T}}\phi(u)\theta_s^I(u)H_s(u)\hat{H}_s(u)dsdu.
\end{align*}
Then $\mathcal{B}_{10}$ is a nonnegative definite quadratic form on $\left(C^{1,\infty}([0, T]\times \mathbb{T})\right)^3$ and
\[
\mathcal{B}_{10}\left((F, G, H), (F, G, H)\right)=\mathcal{B}_4(F, G)+\mathcal{B}_5(G,H)+\mathcal{B}_6(H).
\]
We write $(F, G, H)\simeq (\hat{F}, \hat{G}, \hat{H})$ when
\[
\mathcal{B}_{10}\left((F, G, H), (\hat{F}, \hat{G}, \hat{H})\right)=0.
\]
We denote by $\mathbb{H}$ the completion of $\left(C^{1,\infty}([0, T]\times \mathbb{T})\right)^3/\simeq$ under $\mathcal{B}_{10}$, then $\mathbb{H}$ is a Hilbert space with inner product $\mathcal{B}_{10}$. We have the following lemma about elements $W\in \Omega_3$ with finite moderate deviation rate function value.

\begin{lemma}\label{lemma 5.1.1 moderateRatefuntionEquivalent}
If $W\in \Omega_3$ makes $J_{ini}(W_0)+J_{dyn}(W)<+\infty$, then there exists $(\tilde{F}_W, \tilde{G}_W, \tilde{H}_{W})\in \mathbb{H}$ and $h^W_1, h^W_2, h^W_3\in L^2(\mathbb{T})$ such that
\[
l_2(W,F,G,H)=\mathcal{B}_{10}\left((F, G, H), (\tilde{F}_W, \tilde{G}_W, \tilde{H}_W)\right)
\]
for any $F, G, H\in C^{1,\infty}([0, T]\times \mathbb{T})$ and
\[
W_{0,k}(du)=h^W_k(u)du
\]
for $k=1,2,3$. Furthermore,
\begin{equation}\label{equ 5.1.1 equivexpressionofJdyn}
J_{dyn}(W)=\frac{1}{2}\mathcal{B}_{10}\left((\tilde{F}_W, \tilde{G}_W, \tilde{H}_W), (\tilde{F}_W, \tilde{G}_W, \tilde{H}_W)\right)
=\mathcal{J}_1(W, \tilde{F}_W, \tilde{G}_W, \tilde{H}_W)
\end{equation}
and
\begin{equation}\label{equ 5.1.2 euivexpressionofJini}
J_{ini}(W_0)=\frac{1}{2}\int_{\mathbb{T}}\sum_{k=1}^3\frac{(h^W_k(u))^2}{\rho_{k-1}(u)}+\frac{(\sum_{k=1}^3h^W_k(u))^2}{1-\sum_{k=1}^3\rho_{k-1}(u)}du.
\end{equation}
\end{lemma}

\proof[Proof of Lemma \ref{lemma 5.1.1 moderateRatefuntionEquivalent}]

We first prove the dynamic part. If $J_{dyn}(W)<+\infty$, then
\[
l_2(W, cF, cG, cH)-\frac{1}{2}\mathcal{B}_{10}\left((cF, cG, cH), (cF, cG, cH)\right)\leq J_{dyn}(W)<+\infty
\]
for any $c\in \mathbb{R}$ and $F, G, H\in C^{1,0}([0, T]\times\mathbb{T})$. Let $c=\frac{l_2(W, F, G, H)}{\mathcal{B}_{10}\left((F, G, H), (F, G, H)\right)}$, then
\[
l_2^2(W, F, G, H)\leq 2J_{dyn}(W)\mathcal{B}_{10}\left((F, G, H), (F, G, H)\right).
\]
Therefore, $l_2$ can be extended to a bounded linear operator $\tilde{l}_2$ from $\mathbb{H}$ to $\mathbb{R}$. Hence, existences of $\tilde{F}_W, \tilde{G}_W, \tilde{H}_W$ follow from Riesz representation theorem and Equation \eqref{equ 5.1.1 equivexpressionofJdyn} follows from Cauchy-Schwarz inequality.

Now we prove the initial state part. According to an analysis similar with that leading to Equation \eqref{equ 5.1.1 equivexpressionofJdyn}, there exists $g^W_1, g^W_2, g^W_3$ such that
\[
W_{0,k}(du)=\rho_{k-1}(u)\left(g^W_k(u)-\sum_{l=1}^3g^W_l(u)\rho_{l-1}(u)\right)du
\]
for $k=1,2,3$ and
\[
J_{ini}(W_0)=\frac{1}{2}\int_\mathbb{T}\sum_{k=1}^3\rho_{k-1}(u)(g^W_k(u))^2-\left(\sum_{k=1}^3\rho_{k-1}(u)g^W_k(u)\right)^2du.
\]
Denote $\rho_{k-1}(u)\left(g^W_k(u)-\sum_{l=1}^3g^W_l(u)\rho_{l-1}(u)\right)$ by $h_k^W(u)$ for $k=1,2,3$, then we can solve that
\[
g_k^W=\frac{h_1^W(u)+h_2^W(u)+h_3^W(u)}{1-\sum_{k=1}^3\rho_{k-1}(u)}+\frac{h_k(u)}{\rho_{k-1}(u)}
\]
for $k=1,2,3$ and
\[
\sum_{k=1}^3\rho_{k-1}(u)(g^W_k(u))^2-\left(\sum_{k=1}^3\rho_{k-1}(u)g^W_k(u)\right)^2
=\sum_{k=1}^3\frac{(h^W_k(u))^2}{\rho_{k-1}(u)}+\frac{(\sum_{k=1}^3h^W_k(u))^2}{1-\sum_{k=1}^3\rho_{k-1}(u)},
\]
which completes the proof.

\qed

To further investigate $W\in \Omega_3$ with finite moderate deviation rate function value, we give the following lemma.

\begin{lemma}\label{lemma 5.1.2}
If $W\in \Omega_3$ and $(\tilde{F}, \tilde{G}, \tilde{H})\in \mathbb{H}$ makes
\[
l_2(W, F, G, H)=\mathcal{B}_{10}\left((F, G, H), (\tilde{F}, \tilde{G}, \tilde{H})\right)
\]
for any $(F, G, H)\in \left(C^{1,\infty}([0, T]\times\mathbb{T})\right)^3$, then $W$ is the unique solution to the following $\Omega_3$-valued ordinary differential equation with given initial condition. For any $f,g,h\in C^{1,\infty}(\mathbb{T})$,
\begin{equation}\label{equ 5.1.3}
\begin{cases}
&\frac{d}{ds}W_{s,1}(f)=-W_{s,1}(\mathcal{P}_{1,s}f)-W_{s,3}(\mathcal{P}_{2,s}f)\\
&\text{\quad\quad\quad\quad\quad}-\int_{\mathbb{T}^2}\lambda(u,v)\theta_s^S(u)\theta_s^I(v)f(u)\left(\tilde{G}_s(u)-\tilde{F}_s(u)\right)dudv,\\
&\frac{d}{ds}W_{s,2}(g)=W_{s,1}(\mathcal{P}_{1,s}g)+W_{s,3}(\mathcal{P}_{2,s}g)-W_{s,2}(\mathcal{P}_3g)\\
&\text{\quad\quad\quad\quad\quad}+\int_{\mathbb{T}^2}\lambda(u,v)\theta_s^S(u)\theta_s^I(v)g(u)\left(\tilde{G}_s(u)-\tilde{F}_s(u)\right)dudv\\
&\text{\quad\quad\quad\quad\quad}-\int_{\mathbb{T}}\psi(u)\theta_s^E(u)g(u)(-\tilde{G}_s(u)+\tilde{H}_s(u))du,\\
&\frac{d}{ds}W_{s,3}(h)=W_{s,2}(\mathcal{P}_3h)-W_{s,3}(\mathcal{P}_4h)+\int_{\mathbb{T}}\psi(u)\theta_s^E(u)h(u)(-\tilde{G}_s(u)+\tilde{H}_s(u))du\\
&\text{\quad\quad\quad\quad\quad}+\int_{\mathbb{T}}\phi(u)\theta_s^I(u)h(u)\tilde{H}_s(u).
\end{cases}
\end{equation}
\end{lemma}

\proof[Proof of Lemma \ref{lemma 5.1.2}]

According to Grownwall's inequality, the uniqueness of the solution to Equation \eqref{equ 5.1.3} follows from an analysis similar with that given in the proof of Lemma \ref{lemma 4.2.4}. For any $y\in C^1([0, T])$ and $f\in C^\infty(\mathbb{T})$, let $F(s,u)=y_sf(u)$, then
\[
l_2(W, F, 0, 0)=\mathcal{B}_{10}\left((F,0,0), (\tilde{F}, \tilde{G}, \tilde{H})\right).
\]
By direct calculation,
\[
y_{T}W_{1,T}(f)-y_0W_{1,0}(f)-\int_0^T\partial_sy_sW_{1,s}(f)ds=\int_0^Ty_s\Xi(s)ds,
\]
where
\begin{align*}
\Xi(s)=&-W_{s,1}(\mathcal{P}_{1,s}f)-W_{s,3}(\mathcal{P}_{2,s}f)\\
&-\int_{\mathbb{T}^2}\lambda(u,v)\theta_s^S(u)\theta_s^I(v)f(u)\left(\tilde{G}_s(u)-\tilde{F}_s(u)\right)dudv.
\end{align*}
Since $y$ is arbitrary, $W_{1,s}(f)$ is absolutely continuous and $\frac{d}{ds}W_{1,s}(f)=\Xi(s)$. Hence, $W$ satisfies the first part of Equation \eqref{equ 5.1.3}. Let $G(s,u)=y_sg(u)$ and $H(s,u)=y_sh(u)$, then
\[
l_2(W, 0, G, 0)=\mathcal{B}_{10}\left((0,G,0), (\tilde{F}, \tilde{G}, \tilde{H})\right)
\]
and
\[
l_2(W, 0, 0, H)=\mathcal{B}_{10}\left((0,0,H), (\tilde{F}, \tilde{G}, \tilde{H})\right)
\]
implies that $W$ satisfies rest two parts of Equation \eqref{equ 5.1.3} and the proof is complete.

\qed

According to Lemmas \ref{lemma 5.1.1 moderateRatefuntionEquivalent} and \ref{lemma 5.1.2}, an element $W\in \Omega_3$ making $J_{dyn}(W)<+\infty$ is the solution to Equation \eqref{equ 5.1.3} with $(\tilde{F}, \tilde{G}, \tilde{H})$ taking $(\tilde{F}_W, \tilde{G}_W, \tilde{H}_W)$, which plays key role in the proof of Equation \eqref{equ 2.11 MDP lowerbound}. For mathematical details, see Subsection \ref{subsection 5.3}.

Our next lemma shows that states of different vertices are approximately independent when $N\rightarrow+\infty$.
\begin{lemma}\label{lemma 5.1.3 approximately independent}
There exists $K_4<+\infty$ independent of $N$ such that
\[
\left|P\left(\xi_t^N(i)=k_1, \xi_t^N(j)=k_2\right)-P\left(\xi_t^N(i)=k_1\right)P\left(\xi_t^N(j)=k_2\right)\right|\leq \frac{K_4}{N}
\]
for any $N\geq 1$, $0\leq i\neq j\leq N-1$, $0\leq t\leq T$ and $k_1, k_2\in \{0,1,2,3\}$.
\end{lemma}

The proof of the SIR version of Lemma \ref{lemma 5.1.3 approximately independent} is given in \cite{Xue2022} according to a graphic method, which is still valid for the SEIR model in this paper after some details modified. Hence we omit the proof of Lemma \ref{lemma 5.1.3 approximately independent} here.

\subsection{Proof of Equation \eqref{equ 2.10 MDP upperbound}}\label{subsection 5.2}
In this subsection we prove Equation \eqref{equ 2.10 MDP upperbound}. We first introduce some nations and definitions for later use. For any $F, G, H\in C^{1,\infty}([0,T]\times\mathbb{T})$, we define
\[
\Gamma_{F, G, H}^N(t, \xi^N)=\exp\left(\frac{\gamma_N^2}{N}\left(\eta_{t,1}^N(F_t)+\eta_{t,2}^N(G_t)+\eta_{t,3}^N(H_t)\right)\right).
\]
Furthermore, we define
\[
\mathcal{Y}_{F, G, H}^N(t, \xi^N)=\frac{\Gamma_{F, G, H}^N(t, \xi^N)}{\Gamma_{F, G, H}^N(0, \xi^N)}\exp\left(-\int_0^t\frac{\left(\partial_s+\mathcal{L}^N\right)\Gamma_{F, G, H}^N(s, \xi^N)}{\Gamma_{F, G, H}^N(s, \xi^N)}ds\right),
\]
then $\{\mathcal{Y}_{F, G, H}^N(t, \xi^N)\}_{0\leq t\leq T}$ is a martingale according to Dynkin's martingale formula.

We list some equations which will be repeatedly utilized in the proof of Theorem \ref{Theorem 2.5 MDP}. According to transition-rate functions of our SEIR model and Kolmogorov-Chapman equation, for any $0\leq i\leq N-1$,
\begin{equation}\label{equ 5.2.0}
\begin{cases}
&\frac{d}{dt}\mathbb{E}S_t^N(i)=-\frac{1}{N}\sum_{j=0}^{N-1}\lambda\left(\frac{i}{N},\frac{j}{N}\right)\mathbb{E}\left(S_t^N(i)I_t^N(j)\right),\\
&\frac{d}{dt}\mathbb{E}\mathcal{E}_t^N(i)=\frac{1}{N}\sum_{j=0}^{N-1}\lambda\left(\frac{i}{N},\frac{j}{N}\right)\mathbb{E}\left(S_t^N(i)I_t^N(j)\right)
-\psi\left(\frac{i}{N}\right)\mathbb{E}\mathcal{E}_t^N(i),\\
&\frac{d}{dt}\mathbb{E}I_t^N(i)=\psi\left(\frac{i}{N}\right)\mathbb{E}\mathcal{E}_t^N(i)-\phi\left(\frac{i}{N}\right)\mathbb{E}I_t^N(i).
\end{cases}
\end{equation}
According to the definition of $\mathcal{L}^N$,
\begin{align}\label{equ 5.2.0 two}
\frac{\mathcal{L}^N\Gamma_{F, G, H}^N(s, \xi^N)}{\Gamma_{F, G, H}^N(s, \xi^N)}
=&\frac{1}{N}\sum_{i=0}^{N-1}\sum_{j=0}^{N-1}\lambda\left(\frac{i}{N}, \frac{j}{N}\right)S_s^N(i)I_s^N(j)\left(e^{\frac{\gamma_N}{N}\left(-F_s\left(\frac{i}{N}\right)+G_s\left(\frac{i}{N}\right)\right)}-1\right)\notag\\
&+\sum_{i=0}^{N-1}\psi\left(\frac{i}{N}\right)\mathcal{E}_s^N(i)\left(e^{\frac{\gamma_N}{N}\left(-G_s\left(\frac{i}{N}\right)
+H_s\left(\frac{i}{N}\right)\right)}-1\right) \notag\\
&+\sum_{i=0}^{N-1}\phi\left(\frac{i}{N}\right)I_s^N(i)\left(e^{\frac{\gamma_N}{N}(-H_s\left(\frac{i}{N}\right))}-1\right).
\end{align}
By direct calculation, under Assumption (B),
\begin{align}\label{equ 5.2.0 three}
&\frac{1}{N\gamma_N}\sum_{i=0}^{N-1}\sum_{j=0}^{N-1}\lambda\left(\frac{i}{N}, \frac{j}{N}\right)\left(S_s^N(i)I_s^N(j)-\mathbb{E}\left(S_s^N(i)\right)\mathbb{E}\left(I_s^N(j)\right)\right)f\left(\frac{i}{N}\right)\notag\\
&=\eta_{s,3}^N(\lambda_2)\frac{1}{N}\sum_{i=0}^{N-1}\mathbb{E}S_s^N(i)\lambda_1\left(\frac{i}{N}\right)f\left(\frac{i}{N}\right)
+\eta_{s,1}^N(\lambda_1f)\frac{1}{N}\sum_{j=0}^{N-1}\mathbb{E}I_s^N(j)\lambda_2\left(\frac{j}{N}\right)\notag\\
&\text{\quad}+\eta_{s,1}^N\left(\lambda_1f\right)\frac{1}{N}\sum_{j=0}^{N-1}\lambda_2\left(\frac{j}{N}\right)\left(I_s^N(j)-\mathbb{E}I_s^N(j)\right)
\end{align}
for any $f\in C^\infty(\mathbb{T})$.

The following lemma is crucial for the proof of Equation \eqref{equ 2.10 MDP upperbound}

\begin{lemma}\label{lemma 5.2.1 exponential tightness}
Under Assumptions (A) and (B), $\{\eta^N\}_{N\geq 1}$ are exponentially tight.
\end{lemma}

Details of the proof of Lemma \ref{lemma 5.2.1 exponential tightness} are a little complex, which makes this proof a little long. So here we give a sketch of this proof. Our approach is inspired by the iteration strategy introduced in the proof of Grownwall's inequality given in Appendix 5 of \cite{Ethier1986}. According to this iteration strategy, we can relate the exponential tightness of our process to exponential tightness of the initial state and of an exponential martingale's logarithm. As a consequence, the exponential tightness of our process follows from Doob's inequality and Assumption (A).

\proof[Proof of Lemma \ref{lemma 5.2.1 exponential tightness}]

For any $f_1, f_2, f_3\in C^\infty(\mathbb{T})$, we use $\vec{f}$ to denote $(f_1, f_2, f_3)\in \left(C^\infty(\mathbb{T})\right)^3$. Furthermore, we use $\eta_t^N(\vec{f})$ to denote
\[
\eta_{t,1}^N(f_1)+\eta_{t,2}^N(f_2)+\eta_{t,3}^N(f_3).
\]
According to the criterion given in \cite{Puhalskii1994}, to check the exponential tightness of $\{\eta^N\}_{N\geq 1}$, we only need to show that
\begin{equation}\label{equ 5.2.1}
\limsup_{M\rightarrow+\infty}\limsup_{N\rightarrow+\infty}\frac{N}{\gamma^2_N}\log P\left(\sup_{0\leq t\leq T}|\eta_t^N(\vec{f})|>M\right)=-\infty
\end{equation}
and
\begin{equation}\label{equ 5.2.2}
\limsup_{\delta\rightarrow 0}\limsup_{N\rightarrow+\infty}\frac{N}{\gamma_N^2}
\log\sup_{\sigma\in \mathcal{T}}P\left(\sup_{0\leq t\leq \delta}\left(\eta^N_{t+\sigma}(\vec{f})-\eta^N_\sigma(\vec{f})\right)>\epsilon\right)=-\infty
\end{equation}
for any $\vec{f}\in \left(C^\infty(\mathbb{T})\right)^3$ and $\epsilon>0$, where $\mathcal{T}$ is the set of stopping times of $\{\xi_t^N\}_{t\geq 0}$ bounded by $T$.

We first check Equation \eqref{equ 5.2.1}. According to Theorem \ref{Theorem 2.4 exponentially rapid to hydrodynamic limit}, $\{\mu_{t,3}^N(\lambda_2)\}_{0\leq t\leq T}$ converges weakly to $\{\mu_{t,3}(\lambda_2)\}_{0\leq t\leq T}$ at an exponential rate and hence this convergence is also in $L^1$. As a result, let
\[
D_2^N=\left\{\max_{0\leq t\leq T}\left|\frac{1}{N}\sum_{j=0}^{N-1}\left(I_t^N(j)-\mathbb{E}I_t^N(j)\right)\lambda_2\left(\frac{j}{N}\right)\right|\geq 1\right\},
\]
then $\limsup_{N\rightarrow+\infty}\frac{1}{N}\log P(D_2^N)<0$ and hence
\[
\limsup_{N\rightarrow+\infty}\frac{N}{\gamma_N^2}\log P(D_2^N)=-\infty.
\]
Hence, to check Equation \eqref{equ 5.2.1}, we only need to check that
\begin{equation}\label{equ 5.2.6}
\limsup_{M\rightarrow+\infty}\limsup_{N\rightarrow+\infty}\frac{N}{\gamma^2_N}\log P\left(\sup_{0\leq t\leq T}|\eta_t^N(\vec{f})|>M, (D_2^N)^c\right)=-\infty,
\end{equation}
where $(D_2^N)^c$ is the complementary set of $D_2^N$.

For $\vec{f}\in \left(C(\mathbb{T})\right)^3$, we define
\[
\Gamma_{\vec{f}}^N(t, \xi^N)=\exp\left(\frac{\gamma_N^2}{N}\left(\eta_t^N(\vec{f})\right)\right)
\]
and
\[
\mathcal{Y}_{\vec{f}}^N(t, \xi^N)=\frac{\Gamma_{\vec{f}}^N(t, \xi^N)}{\Gamma_{\vec{f}}^N(0, \xi^N)}\exp\left(-\int_0^t\frac{\left(\partial_s+\mathcal{L}^N\right)\Gamma_{\vec{f}}^N(s, \xi^N)}{\Gamma_{\vec{f}}^N(s, \xi^N)}ds\right).
\]
Then by Feynman-Kac formula, $\{\mathcal{Y}_{\vec{f}}^N(t, \xi^N)\}_{0\leq t\leq T}$ is a martingale. For $0\leq s\leq T$, let $\mathcal{A}_s^N$ be the linear operator from $\left(C^\infty(\mathbb{T})\right)^3$ to $\left(C^\infty(\mathbb{T})\right)^3$ such that
\[
\left(\mathcal{A}_s^N\vec{f}\right)_1(u)=\frac{1}{N}\sum_{j=1}^N\mathbb{E}I_s^N(j)\lambda\left(u, \frac{j}{N}\right)\left(-f_1(u)
+f_2(u)\right),
\]
\[
\left(\mathcal{A}_s^N\vec{f}\right)_2(u)=\psi(u)\left(-f_2(u)+f_3(u)\right),
\]
and
\begin{align*}
\left(\mathcal{A}_s^N\vec{f}\right)_3(u)=\frac{1}{N}\sum_{j=1}^N\mathbb{E}S_s^N(j)\lambda\left(\frac{j}{N}, u\right)\left(-f_1\left(\frac{j}{N}\right)+f_2\left(\frac{j}{N}\right)\right)-\phi(u)f_3(u),
\end{align*}
for any $\vec{f}\in \left(C^\infty(\mathbb{T})\right)^3$ and $u\in \mathbb{T}$. For any $\vec{f}\in \left(C^\infty(\mathbb{T})\right)^3$, let
\[
\|\vec{f}\|_\infty=\max\left\{\|f_1\|_\infty, \|f_2\|_\infty, \|f_3\|_\infty\right\},
\]
then $\|\mathcal{A}_s^N\vec{f}\|_\infty\leq K_5\|\vec{f}\|_\infty$ for any $\vec{f}\in \left(C^\infty(\mathbb{T})\right)^3$, where
\[
K_5=2\|\lambda\|_\infty+2\|\psi\|_\infty+\|\phi\|_\infty.
\]
Let $\mathcal{B}_{11}$ be the linear operator from $\left(C^\infty(\mathbb{T})\right)^3$ to $\left(C^\infty(\mathbb{T})\right)^3$ such that
\[
\mathcal{B}_{11}\vec{f}=\left(\lambda_1(-f_1+f_2), 0, 0\right)
\]
for any $\vec{f}\in \left(C^\infty(\mathbb{T})\right)^3$, then
\[
\|\mathcal{B}_{11}\vec{f}\|_\infty\leq 2\|\lambda_1\|_\infty\|\vec{f}\|_\infty.
\]
By Equations \eqref{equ 5.2.0}-\eqref{equ 5.2.0 three} and Taylor's expansion of $\exp\left(\frac{\gamma_N}{N}x\right)$ up to the second order with Lagrange's remainder, we have
\begin{equation}\label{equ 5.2.3}
\mathcal{Y}_{\vec{f}}^N(t, \xi^N)=\exp\left(\frac{\gamma_N^2}{N}\left(\eta_t^N(\vec{f})-\eta_0^N(\vec{f})-\int_0^t\eta_s^N(\mathcal{A}_s^N\vec{f})ds
+\sum_{k=1}^4\varepsilon_{k,t}^N(\vec{f})\right)\right),
\end{equation}
where
\begin{align*}
\varepsilon_{1,t}^N(\vec{f})=\int_0^t&\frac{1}{N\gamma_N}\sum_{i=1}^N\sum_{j=1}^N\lambda_1\left(\frac{i}{N}\right)\lambda_2\left(\frac{j}{N}\right)\times\\
&\left(\mathbb{E}\left(S_s^N(i)I_s^N(j)\right)
-\mathbb{E}\left(S_s^N(i)\right)\mathbb{E}\left(I_s^N(j)\right)\right)\left(-f_1\left(\frac{i}{N}\right)+f_2\left(\frac{i}{N}\right)\right)ds
\end{align*}
and
\begin{align*}
|\varepsilon_{2,t}^N(\vec{f})|\leq T K_6\|\vec{f}\|_\infty^2
\end{align*}
for some $K_6<+\infty$ independent of $N, \vec{f}, t$ and
\[
|\varepsilon_{3,t}^N(\vec{f})|\leq \frac{\gamma_NT}{N}K_7\|\vec{f}\|_\infty^3\exp\left(2\|\vec{f}\|_\infty\frac{\gamma_N}{N}\right)
\]
for some $K_7<+\infty$ independent of $N, \vec{f}, t$ and
\[
\varepsilon_{4,t}^N(\vec{f})=-\int_0^t\eta_{s,1}^N\left(\lambda_1(\cdot)(-f_1(\cdot)+f_2(\cdot))\right)
\left(\frac{1}{N}\sum_{j=0}^{N-1}\left(I_s^N(j)-\mathbb{E}I_s^N(j)\right)\lambda_2\left(\frac{j}{N}\right)\right)ds.
\]

According to Lemma \ref{equ 5.1.3},
\[
|\varepsilon_{1,t}^N(f)|\leq \frac{K_8T}{\gamma_N}\|\vec{f}\|_\infty,
\]
where $K_8<+\infty$ is independent of $N, \vec{f}, t$. Let $\mathcal{Z}_t^N(\vec{f})=\frac{N}{\gamma_N^2}\log \mathcal{Y}_{\vec{f}}^N(t, \xi^N)$ and
\begin{align*}
\varepsilon_5^N(\vec{f})=&|\eta_0^N(\vec{f})|+K_6T\|\vec{f}\|_\infty^2+\frac{\gamma_NT}{N}K_7\|\vec{f}\|_\infty^3\exp\left(2\|\vec{f}\|_\infty\frac{\gamma_N}{N}\right)
\\
&+\frac{K_8T}{\gamma_N}\|\vec{f}\|_\infty+\max_{0\leq t\leq T}\mathcal{Z}_t^N(\vec{f}),
\end{align*}
then on $(D_2^N)^c$,
\begin{equation}\label{equ 5.2.4}
|\eta_t^N(\vec{f})|\leq \int_0^t|\eta_s^N(\mathcal{A}_s^N\vec{f})|+|\eta_s^N(\mathcal{B}_{11}\vec{f})|ds+\varepsilon_5^N(\vec{f})
\end{equation}
for any $t\geq 0$ and $\vec{f}\in \left(C^\infty(\mathbb{T})\right)^3$. By repeated utilizing Equation \eqref{equ 5.2.4} for $k$ times, we have
\begin{align*}
&|\eta_t^N(\vec{f})|\\
&\leq \int_0^t|\eta_{t_1}^N(\mathcal{A}_{t_1}^N\vec{f})|+|\eta_{t_1}^N(\mathcal{B}_{11}\vec{f})|dt_1+\varepsilon_5^N(\vec{f}) \\
&\leq \int_0^t\int_0^{t_1}|\eta_{t_2}^N\left(\mathcal{A}_{t_2}^N\mathcal{A}_{t_1}^N\vec{f}\right)|dt_1dt_2
+\int_0^t\int_0^{t_1}|\eta_{t_2}^N\left(\mathcal{B}_{11}\mathcal{B}_{11}\vec{f}\right)|dt_1dt_2\\
&\text{\quad}+\int_0^t\int_0^{t_1}|\eta_{t_2}^N\left(\mathcal{B}_{11}\mathcal{A}_{t_2}^N\vec{f}\right)|dt_1dt_2
+\int_0^t\int_0^{t_1}|\eta_{t_2}^N\left(\mathcal{A}_{t_1}^N\mathcal{B}_{11}\vec{f}\right)|dt_1dt_2\\
&\text{\quad}+\int_0^t\varepsilon_5^N(\mathcal{A}_{t_1}^N\vec{f})dt_1+\int_0^t\varepsilon_5^N(\mathcal{B}_{11}\vec{f})d_{t_1}+\varepsilon_5^N(\vec{f})\\
&\ldots\\
&\leq \int_0^t\int_0^{t_1}\ldots\int_0^{t_{k-1}}\sum_{\mathcal{C}_{t_j}^N\in \{A_{t_j}^N, \mathcal{B}_{11}\}\atop \text{~for~}
1\leq j\leq k}|\eta_{t_k}^N\left(\mathcal{C}_{t_k}^N\mathcal{C}_{t_{k-1}}^N\ldots\mathcal{C}_{t_1}^N\vec{f}\right)|dt_1\ldots dt_k\\
&\text{\quad\quad}+\varepsilon_5^N(\vec{f})+\sum_{l=1}^{k-1}\int_0^t\int_0^{t_1}\ldots\int_0^{t_{l-1}}
\sum_{\mathcal{C}_{t_j}^N\in \{A_{t_j}^N, \mathcal{B}_{11}\}\atop \text{~for~}
1\leq j\leq l}\varepsilon_5^N\left(\mathcal{C}_{t_{l}}^N\ldots\mathcal{C}_{t_1}^N\vec{f}\right)dt_1dt_2\ldots dt_{l}
\end{align*}
on $(D_2^N)^c$. There exists $K_9^N<\infty$ independent of $\vec{f}, t$ such that
\[
|\eta_t^N(\vec{f})|\leq K_9^N\|\vec{f}\|_\infty.
\]
Hence, on $(D_2^N)^c$,
\begin{align*}
&\int_0^t\int_0^{t_1}\ldots\int_0^{t_{k-1}}\sum_{\mathcal{C}_{t_j}^N\in \{A_{t_j}^N, \mathcal{B}_{11}\}\atop \text{~for~}
1\leq j\leq k}|\eta_{t_k}^N\left(\mathcal{C}_{t_k}^N\mathcal{C}_{t_{k-1}}^N\ldots\mathcal{C}_{t_1}^N\vec{f}\right)|dt_1\ldots dt_k\\
&\leq \frac{K_9^N(K_5+2\|\lambda\|_1)^kT^k\|\vec{f}\|_\infty}{k!}.
\end{align*}
Let $k\rightarrow+\infty$, then we have
\begin{equation}\label{equ 5.2.5}
|\eta_t^N(\vec{f})|\leq \varepsilon_5^N(\vec{f})+\sum_{l=1}^{+\infty}\int_0^t\ldots\int_0^{t_{l-1}}
\sum_{\mathcal{C}_{t_j}^N\in \{A_{t_j}^N, \mathcal{B}_{11}\}\atop \text{~for~}
1\leq j\leq l}\varepsilon_5^N\left(\mathcal{C}_{t_{l}}^N\ldots\mathcal{C}_{t_1}^N\vec{f}\right)dt_1dt_2\ldots dt_{l}
\end{equation}
on $(D_2^N)^c$. Without loss of generality, we assume that $K_5+2\|\lambda_1\|_\infty<1$. Otherwise, for sufficiently large $K$, $\{\xi^N_{\frac{t}{K}}\}_{t\geq 0}$ is a version of our SEIR model with parameters
\[
\left(\hat{\lambda}_1, \hat{\lambda}_2, \hat{\psi}, \hat{\phi}\right)=\left(\frac{1}{\sqrt{K}}\lambda_1, \frac{1}{\sqrt{K}}\lambda_2, \frac{1}{K}\psi, \frac{1}{K}\phi\right)
\]
and then $2\|\hat{\lambda}\|_\infty+2\|\hat{\psi}\|_\infty+\|\hat{\phi}\|_\infty+2\|\hat{\lambda}\|_1<1$. When $K_5+2\|\lambda_1\|_\infty<1$,
\[
\|\mathcal{C}_{t_{l}}^N\ldots\mathcal{C}_{t_1}^N\vec{f}\|_\infty\leq \|f\|_\infty
\]
for $l\geq 1$. Hence there exists $N_0$ independent of $l$ such that
\begin{align*}
&\varepsilon_5^N\left(\mathcal{C}_{t_{l}}^N\ldots\mathcal{C}_{t_1}^N\vec{f}\right)\\
&\leq |\eta_0^N(\mathcal{C}_{t_{l}}^N\ldots\mathcal{C}_{t_1}^N\vec{f})|
+\max_{0\leq t\leq T}\mathcal{Z}_t^N(\mathcal{C}_{t_{l}}^N\ldots\mathcal{C}_{t_1}^N\vec{f})+(K_6T+1)\|\mathcal{C}_{t_{l}}^N\ldots\mathcal{C}_{t_1}^N\vec{f}\|^2_\infty\\
&\leq |\eta_0^N(\mathcal{C}_{t_{l}}^N\ldots\mathcal{C}_{t_1}^N\vec{f})|
+\max_{0\leq t\leq T}\mathcal{Z}_t^N(\mathcal{C}_{t_{l}}^N\ldots\mathcal{C}_{t_1}^N\vec{f})+K_{10}\|\vec{f}\|^2_\infty
\end{align*}
for any $N\geq N_0$ and $l\geq 1$, where $K_{10}=K_6T+1$. As a result, on $(D_2^N)^c$,
\begin{equation}\label{equ 5.2.7}
\max_{0\leq t\leq T}|\eta_t^N(\vec{f})|\leq \varepsilon_6^N+\varepsilon_7^N+\varepsilon_8^N,
\end{equation}
where
\[
\varepsilon_6^N=\sum_{l=1}^{+\infty}\int_0^T\ldots\int_0^{t_{l-1}}
\sum_{\mathcal{C}_{t_j}^N\in \{A_{t_j}^N, \mathcal{B}_{11}\}\atop \text{~for~}
1\leq j\leq l}\max_{0\leq t\leq T}\mathcal{Z}_t^N\left(\mathcal{C}_{t_{l}}^N\ldots\mathcal{C}_{t_1}^N\vec{f}\right)dt_1dt_2\ldots dt_{l},
\]
\[
\varepsilon_7^N=\sum_{l=1}^{+\infty}\int_0^T\ldots\int_0^{t_{l-1}}
\sum_{\mathcal{C}_{t_j}^N\in \{A_{t_j}^N, \mathcal{B}_{11}\}\atop \text{~for~}
1\leq j\leq l}\left|\eta_0^N\left(\mathcal{C}_{t_{l}}^N\ldots\mathcal{C}_{t_1}^N\vec{f}\right)\right|dt_1dt_2\ldots dt_{l}
\]
and
\[
\varepsilon_8^N=\max_{0\leq t\leq T}\mathcal{Z}_t^N\left(\vec{f}\right)+|\eta_0^N(\vec{f})|+K_{10}\|\vec{f}\|_\infty^2e^T.
\]
For $l\geq 1$, let $D_{3,l}^N$ be the event that
\[
\int_0^T\ldots\int_0^{t_{l-1}}\sum_{\mathcal{C}_{t_j}^N\in \{A_{t_j}^N, \mathcal{B}_{11}\}\atop \text{~for~}
1\leq j\leq l}\max_{0\leq t\leq T}\mathcal{Z}_t^N\left(\mathcal{C}_{t_{l}}^N\ldots\mathcal{C}_{t_1}^N\vec{f}\right)dt_1dt_2\ldots dt_{l}
\geq \frac{(4T)^l}{l!}\frac{M}{e^{4T}}.
\]
Then $P(\varepsilon_6^N>M)\leq \sum_{l=1}^{\infty}P\left(D_{3,l}^N\right)$. According to Markov's inequality,
\begin{align}\label{equ 5.2.8}
&P\left(D_{3,l}^N\right)\notag\\
&=P\Bigg(\frac{l!}{2^lT^l}\int_0^T\ldots\int_0^{t_{l-1}}\sum_{\mathcal{C}_{t_j}^N\in \{A_{t_j}^N, \mathcal{B}_{11}\}\atop \text{~for~}
1\leq j\leq l}
2\max_{0\leq t\leq T}\frac{\gamma_N^2}{N}\mathcal{Z}_t^N\left(\mathcal{C}_{t_{l}}^N\ldots\mathcal{C}_{t_1}^N\vec{f}\right)dt_1\ldots dt_l
\notag\\
&\text{\quad\quad\quad\quad}\geq \frac{2^{l+1}M}{e^{4T}}\frac{\gamma_N^2}{N}\Bigg)\notag\\
&\leq \exp\left(-\frac{2^{l+1}M}{e^{4T}}\frac{\gamma_N^2}{N}\right)\times\\
&\text{\quad\quad}\mathbb{E}\exp\left(\frac{l!}{2^lT^l}\int_0^T\ldots\int_0^{t_{l-1}}\sum_{\mathcal{C}_{t_j}^N\in \{A_{t_j}^N, \mathcal{B}_{11}\}\atop \text{~for~}
1\leq j\leq l}
2\max_{0\leq t\leq T}\frac{\gamma_N^2}{N}\mathcal{Z}_t^N\left(\mathcal{C}_{t_{l}}^N\ldots\mathcal{C}_{t_1}^N\vec{f}\right)dt_1\ldots dt_l\right). \notag
\end{align}
Since $\int_0^T\ldots\int_0^{t_{l-1}}1dt_1\ldots dt_l=\frac{T^l}{l!}$ and the sum $\sum_{\mathcal{C}_{t_j}^N\in \{A_{t_j}^N, \mathcal{B}_{11}\}\atop \text{~for~}
1\leq j\leq l}$ has $2^l$ terms, by Jensen's inequality,
\begin{align}\label{equ 5.2.9}
&\mathbb{E}\exp\left(\frac{l!}{2^lT^l}\int_0^T\ldots\int_0^{t_{l-1}}\sum_{\mathcal{C}_{t_j}^N\in \{A_{t_j}^N, \mathcal{B}_{11}\}\atop \text{~for~}
1\leq j\leq l}
2\max_{0\leq t\leq T}\frac{\gamma_N^2}{N}\mathcal{Z}_t^N\left(\mathcal{C}_{t_{l}}^N\ldots\mathcal{C}_{t_1}^N\vec{f}\right)dt_1\ldots dt_l\right) \notag\\
&\leq \frac{l!}{2^lT^l}\int_0^T\ldots\int_0^{t_{l-1}}\sum_{\mathcal{C}_{t_j}^N\in \{A_{t_j}^N, \mathcal{B}_{11}\}\atop \text{~for~}
1\leq j\leq l}
\mathbb{E}e^{2\max_{0\leq t\leq T}\frac{\gamma_N^2}{N}\mathcal{Z}_t^N\left(\mathcal{C}_{t_{l}}^N\ldots\mathcal{C}_{t_1}^N\vec{f}\right)}dt_1\ldots dt_l.
\end{align}
According to the definition of $\mathcal{Z}_t^N$, $\{\exp\left(\frac{\gamma_N^2}{N}\mathcal{Z}_t^N(\vec{f})\right)\}_{0\leq t\leq T}$ is a martingale. Hence, by Doob's inequality,
\begin{equation*}
\mathbb{E}e^{2\max_{0\leq t\leq T}\frac{\gamma_N^2}{N}\mathcal{Z}_t^N\left(\mathcal{C}_{t_{l}}^N\ldots\mathcal{C}_{t_1}^N\vec{f}\right)}
\leq 4\mathbb{E}e^{2\frac{\gamma_N^2}{N}\mathcal{Z}_T^N\left(\mathcal{C}_{t_{l}}^N\ldots\mathcal{C}_{t_1}^N\vec{f}\right)}.
\end{equation*}
According to the definition of $\mathcal{Y}_{\vec{f}}^N(t, \xi_t^N)$ and Taylor's expansion of $\exp\left(\frac{\gamma_N^2}{N}x\right)$ up to the second order with Lagrange's remainder, we have
\begin{align*}
&e^{2\frac{\gamma_N^2}{N}\mathcal{Z}_T^N\left(\vec{f}\right)}\leq\\
&\mathcal{Y}_{2\vec{f}}^N(T, \xi_T^N)\exp\left(\frac{\gamma_N^2}{N}\left(2K_6T\|\vec{f}\|_\infty^2+\frac{6\gamma_NT}{N}K_7\|\vec{f}\|_\infty^3
\exp\left(4\|\vec{f}\|_\infty\frac{\gamma_N}{N}\right)\right)\right).
\end{align*}
Then, since $\{\mathcal{Y}_{2\vec{f}}^N(t, \xi_t^N)\}_{0\leq t\leq T}$ is a martingale and $\|\mathcal{C}_{t_{l}}^N\ldots\mathcal{C}_{t_1}^N\vec{f}\|_\infty\leq \|\vec{f}\|_\infty$, there exists $N_1\geq 1$ independent of $l$ and $\mathcal{C}_{t_{l}}^N,\ldots,\mathcal{C}_{t_1}^N$ such that
\begin{equation}\label{equ 5.2.9 two}
\mathbb{E}e^{2\max_{0\leq t\leq T}\frac{\gamma_N^2}{N}\mathcal{Z}_t^N\left(\mathcal{C}_{t_{l}}^N\ldots\mathcal{C}_{t_1}^N\vec{f}\right)}
\leq \exp\left(\frac{\gamma_N^2}{N}K_{11}\right)
\end{equation}
when $N\geq N_1$, where $K_{11}=(2K_6T+1)\|\vec{f}\|_\infty^2$. In conclusion, by Equations \eqref{equ 5.2.8}, \eqref{equ 5.2.9} and \eqref{equ 5.2.9 two},
\[
P(D_{3,l}^N)\leq 4\exp\left(\frac{\gamma_N^2}{N}\left(-\frac{2^{l+1}M}{e^{4T}}+K_{11}\right)\right)
\]
and hence
\[
P(\varepsilon_6^N>M)
\leq\sum_{l=1}^{+\infty}\exp\left(\frac{\gamma_N^2}{N}\left(-\frac{2^{l+1}M}{e^{4T}}+K_{11}\right)\right)
\leq 8\exp\left(\frac{\gamma_N^2}{N}\left(K_{11}-\frac{2M}{e^{4T}}\right)\right)
\]
for sufficiently large $N$. As a result,
\begin{equation}\label{equ 5.2.10}
\limsup_{M\rightarrow+\infty}\limsup_{N\rightarrow+\infty}\frac{N}{\gamma^2_N}\log P\left(\varepsilon_6^N>M\right)=-\infty.
\end{equation}
According to Assumption (A), Taylor's expansion up to the second order and the fact that $\exp(|x|)\leq \exp(x)+\exp(-x)$, there exists $K_{12}<+\infty$ independent of $N, \vec{f}$ such that
\begin{equation}\label{equ 5.2.11}
\mathbb{E}\exp\left(\frac{\gamma_N^2}{N}\left|\eta_0^N(\vec{f})\right|\right)\leq 2\exp\left(\frac{\gamma_N^2}{N}K_{12}\|\vec{f}\|^2_\infty\right)
\end{equation}
when $N$ is sufficiently large. By Equations \eqref{equ 5.2.9 two} and \eqref{equ 5.2.11}, an analysis similar with that leading to Equation \eqref{equ 5.2.10} implies that
\begin{equation}\label{equ 5.2.12}
\limsup_{M\rightarrow+\infty}\limsup_{N\rightarrow+\infty}\frac{N}{\gamma^2_N}\log P\left(\varepsilon_7^N+\varepsilon_8^N>M\right)=-\infty.
\end{equation}
Equation \eqref{equ 5.2.6} follows from Equations \eqref{equ 5.2.7}, \eqref{equ 5.2.10} and \eqref{equ 5.2.12}. Consequently, Equation \eqref{equ 5.2.1} follows from Equation \eqref{equ 5.2.6} as we have introduced.

Now we check Equation \eqref{equ 5.2.2}. We still utilize the iteration strategy. Hence, some details will be omitted when they are similar with those in the check of Equation \eqref{equ 5.2.1}.

According to an analysis similar with that leading to Equation \eqref{equ 5.2.7}, on $(D_2^N)^c$, we have
\begin{equation}\label{equ 5.2.13}
\sup_{0\leq t\leq \delta}\left|\eta_{t+\sigma}^N(\vec{f})-\eta_{\sigma}^N(\vec{f})\right|\leq \varepsilon_9^N+\varepsilon_{10}^N+\varepsilon_{11}^N,
\end{equation}
where
\[
\varepsilon_9^N=\sum_{l=1}^{+\infty}\int_\sigma^{\sigma+\delta}\int_0^{t_1}\ldots\int_0^{t_{l-1}}
\sum_{\mathcal{C}_{t_j}^N\in \{A_{t_j}^N, \mathcal{B}_{11}\}\atop \text{~for~}
1\leq j\leq l}\max_{0\leq t\leq T}\mathcal{Z}_t^N\left(\mathcal{C}_{t_{l}}^N\ldots\mathcal{C}_{t_1}^N\vec{f}\right)dt_1dt_2\ldots dt_{l},
\]
\[
\varepsilon_{10}^N=\sum_{l=1}^{+\infty}\int_{\sigma}^{\sigma+\delta}\int_0^{t_1}\ldots\int_0^{t_{l-1}}
\sum_{\mathcal{C}_{t_j}^N\in \{A_{t_j}^N, \mathcal{B}_{11}\}\atop \text{~for~}
1\leq j\leq l}\left|\eta_0^N\left(\mathcal{C}_{t_{l}}^N\ldots\mathcal{C}_{t_1}^N\vec{f}\right)\right|dt_1dt_2\ldots dt_{l},
\]
and
\[
\varepsilon_{11}^N=\max_{0\leq t\leq \delta}\left(\mathcal{Z}_{t+\sigma}^N(\vec{f})-\mathcal{Z}_\sigma^N(\vec{f})\right)+K_6\|\vec{f}\|_\infty^2\delta e^T.
\]
Let $D_{4,l}^N$ be the event that
\[
\int_\sigma^{\sigma+\delta}\int_0^{t_1}\ldots\int_0^{t_{l-1}}
\sum_{\mathcal{C}_{t_j}^N\in \{A_{t_j}^N, \mathcal{B}_{11}\}\atop \text{~for~}
1\leq j\leq l}\max_{0\leq t\leq T}\mathcal{Z}_t^N\left(\mathcal{C}_{t_{l}}^N\ldots\mathcal{C}_{t_1}^N\vec{f}\right)dt_1dt_2\ldots dt_{l}
\geq \frac{(20T)^l\epsilon}{l!e^{20T}},
\]
then $P(\varepsilon_9^N>\epsilon)\leq \sum_{l=1}^{+\infty}P(D_{4,l}^N)$ and $P(D_{4,l}^N)\leq P(\hat{D}_{4,l}^N)$, where $\hat{D}_{4,l}^N$ is the event
that
\begin{align*}
&\int_\sigma^{\sigma+\delta}dt_1\int_0^T dt_2\int_0^{t_2}\ldots\int_0^{t_{l-1}}\sum_{\mathcal{C}_{t_j}^N\in \{A_{t_j}^N, \mathcal{B}_{11}\}\atop \text{~for~}
1\leq j\leq l}\max_{0\leq t\leq T}\mathcal{Z}_t^N\left(\mathcal{C}_{t_{l}}^N\ldots\mathcal{C}_{t_1}^N\vec{f}\right)dt_3\ldots dt_l\\
&\geq  \frac{(20T)^l\epsilon}{l!e^{20T}}.
\end{align*}
Since $\int_\sigma^{\sigma+\delta}dt_1\int_0^T dt_2\int_0^{t_2}\ldots\int_0^{t_{l-1}}1dt_3\ldots dt_l=\frac{\delta T^{l-1}}{(l-1)!}$ and the sum $\sum_{\mathcal{C}_{t_j}^N\in \{A_{t_j}^N, \mathcal{B}_{11}\}\atop \text{~for~}
1\leq j\leq l}$ has $2^l$ terms, by Jensen's inequality, we have
\[
P(\hat{D}_{4,l}^N)\leq P(\tilde{D}_{4,l}^N),
\]
where $\tilde{D}_{4,l}^N$ is the event that
\begin{align*}
&\frac{(l-1)!}{2^l\delta T^{l-1}}\int_\sigma^{\sigma+\delta}dt_1\int_0^T dt_2\int_0^{t_2}\ldots\int_0^{t_{l-1}}\\
&\text{\quad\quad\quad}\sum_{\mathcal{C}_{t_j}^N\in \{A_{t_j}^N, \mathcal{B}_{11}\}\atop \text{~for~}
1\leq j\leq l}\exp\left(2\frac{\gamma_N^2}{N}\max_{0\leq t\leq T}\mathcal{Z}_t^N\left(\mathcal{C}_{t_{l}}^N\ldots\mathcal{C}_{t_1}^N\vec{f}\right)\right)dt_3\ldots dt_l\\
&\geq \exp\left(\frac{10^l\epsilon 2T}{le^{20T}\delta}\frac{\gamma_N^2}{N}\right).
\end{align*}
By Markov's inequality,
\begin{align*}
P(\tilde{D}_{4,l}^N)&\leq\exp\left(-\frac{10^l\epsilon 2T}{le^{20T}\delta}\frac{\gamma_N^2}{N}\right)\times\\
&\text{\quad\quad\quad}\mathbb{E}\Bigg(\frac{(l-1)!}{2^l\delta T^{l-1}}\int_\sigma^{\sigma+\delta}dt_1\int_0^T dt_2\int_0^{t_2}\ldots\int_0^{t_{l-1}}\\
&\text{\quad\quad\quad\quad}\sum_{\mathcal{C}_{t_j}^N\in \{A_{t_j}^N, \mathcal{B}_{11}\}\atop \text{~for~}
1\leq j\leq l}\exp\left(2\frac{\gamma_N^2}{N}\max_{0\leq t\leq T}\mathcal{Z}_t^N\left(\mathcal{C}_{t_{l}}^N\ldots\mathcal{C}_{t_1}^N\vec{f}\right)\right)dt_3\ldots dt_l\Bigg)\\
&\leq\exp\left(-\frac{10^l\epsilon 2T}{le^{20T}\delta}\frac{\gamma_N^2}{N}\right)\times\\
&\text{\quad\quad}\mathbb{E}\Bigg(\frac{(l-1)!}{2^l\delta T^{l-1}}\int_0^Tdt_1\int_0^T dt_2\int_0^{t_2}\ldots\int_0^{t_{l-1}}\\
&\text{\quad\quad\quad\quad}\sum_{\mathcal{C}_{t_j}^N\in \{A_{t_j}^N, \mathcal{B}_{11}\}\atop \text{~for~}
1\leq j\leq l}\exp\left(2\frac{\gamma_N^2}{N}\max_{0\leq t\leq T}\mathcal{Z}_t^N\left(\mathcal{C}_{t_{l}}^N\ldots\mathcal{C}_{t_1}^N\vec{f}\right)\right)dt_3\ldots dt_l\Bigg).
\end{align*}
Then by Equation \eqref{equ 5.2.9 two}, for sufficiently large $N$,
\[
P(\tilde{D}_{4,l}^N)\leq \frac{T}{\delta}\exp\left(\frac{\gamma_N^2}{N}\left(K_{11}-\frac{10^l}{l}\frac{2T\epsilon}{e^{20T}\delta}\right)\right)
\]
and hence
\begin{align*}
P(\varepsilon_9^N>\epsilon)&\leq \sum_{l=1}^{+\infty}\frac{T}{\delta}\exp\left(\frac{\gamma_N^2}{N}\left(K_{11}-\frac{10^l}{l}\frac{2T\epsilon}{e^{20T}\delta}\right)\right)\\
&\leq \frac{2T}{\delta}\exp\left(\frac{\gamma_N^2}{N}\left(K_{11}-\frac{20T\epsilon}{\delta e^{20T}}\right)\right).
\end{align*}
As a result,
\begin{equation}\label{equ 5.2.14}
\limsup_{\delta\rightarrow 0}\limsup_{N\rightarrow+\infty}\frac{N}{\gamma_N^2}
\log\sup_{\sigma\in \mathcal{T}}P\left(\varepsilon_9^N>\epsilon\right)=-\infty.
\end{equation}
By Equation \eqref{equ 5.2.11}, an analysis similar with that leading to Equation \eqref{equ 5.2.14} implies that
\begin{equation}\label{equ 5.2.15}
\limsup_{\delta\rightarrow 0}\limsup_{N\rightarrow+\infty}\frac{N}{\gamma_N^2}
\log\sup_{\sigma\in \mathcal{T}}P\left(\varepsilon_{10}^N>\epsilon\right)=-\infty.
\end{equation}
According to the definition of $\mathcal{Y}_{\vec{f}}^N(t, \xi_t^N)$ and Taylor's expansion of $\exp\left(\frac{\gamma_N^2}{N}x\right)$ up to the second order with Lagrange's remainder, we have
\begin{align*}
&e^{\frac{\gamma_N^2}{N}\left(\mathcal{Z}_{\sigma+t}^N(K\vec{f})-\mathcal{Z}_\sigma^N(K\vec{f})\right)}\leq \frac{\mathcal{Y}_{K\vec{f}}^N(\sigma+t, \xi_{\sigma+t}^N)}{\mathcal{Y}_{K\vec{f}}^N(\sigma, \xi_\sigma^N)}\times\\
&\text{\quad}\exp\left(\frac{\gamma_N^2}{N}\left(2K_6\delta(K^2-K)\|\vec{f}\|_\infty^2+\frac{6\gamma_N(K^2+K)\delta}{N}K_7\|\vec{f}\|_\infty^3
\exp\left(4\|\vec{f}\|_\infty\frac{\gamma_N}{N}\right)\right)\right)
\end{align*}
for any $K>0$, $0\leq t\leq \delta$ and $N$ sufficiently large. Hence, for sufficiently large $N$,
\begin{align*}
&P\left(\sup_{0\leq t\leq \delta}\left(\mathcal{Z}_{\sigma+t}^N(\vec{f})-\mathcal{Z}_\sigma^N(\vec{f})\right)\geq \epsilon\right)\\
&\leq P\left(\sup_{0\leq t\leq \delta}\frac{\mathcal{Y}_{K\vec{f}}^N(\sigma+t, \xi_{\sigma+t}^N)}{\mathcal{Y}_{K\vec{f}}^N(\sigma, \xi_\sigma^N)}
\geq \exp\left(\left(K\epsilon-3K_6\delta(K^2-K)\|\vec{f}\|_\infty^2\right)\frac{\gamma_N^2}{N}\right)\right).
\end{align*}

Since $\left\{\frac{\mathcal{Y}_{K\vec{f}}^N(\sigma+t, \xi_{\sigma+t}^N)}{\mathcal{Y}_{K\vec{f}}^N(\sigma, \xi_\sigma^N)}\right\}_{t\geq 0}$ is a martingale, by Doob's inequality,
\begin{align*}
&P\left(\sup_{0\leq t\leq \delta}\frac{\mathcal{Y}_{K\vec{f}}^N(\sigma+t, \xi_{\sigma+t}^N)}{\mathcal{Y}_{K\vec{f}}^N(\sigma, \xi_\sigma^N)}
\geq \exp\left(\left(K\epsilon-3K_6\delta(K^2-K)\|\vec{f}\|_\infty^2\right)\frac{\gamma_N^2}{N}\right)\right)\\
&\leq \exp\left(-\left(K\epsilon-3K_6\delta(K^2-K)\|\vec{f}\|_\infty^2\right)\frac{\gamma_N^2}{N}\right)
\end{align*}
and hence
\[
\limsup_{\delta\rightarrow 0}\limsup_{N\rightarrow+\infty}\frac{N}{\gamma_N^2}
\log\sup_{\sigma\in \mathcal{T}}P\left(\sup_{0\leq t\leq \delta}\left(\mathcal{Z}_{\sigma+t}^N(\vec{f})-\mathcal{Z}_\sigma^N(\vec{f})\right)\geq \epsilon\right)\leq -K\epsilon.
\]
Since $K$ is arbitrary, let $K\rightarrow+\infty$ and then
\[
\limsup_{\delta\rightarrow 0}\limsup_{N\rightarrow+\infty}\frac{N}{\gamma_N^2}
\log\sup_{\sigma\in \mathcal{T}}P\left(\sup_{0\leq t\leq \delta}\left(\mathcal{Z}_{\sigma+t}^N(\vec{f})-\mathcal{Z}_\sigma^N(\vec{f})\right)\geq \epsilon\right)=-\infty.
\]
As a result,
\begin{equation}\label{equ 5.2.16}
\limsup_{\delta\rightarrow 0}\limsup_{N\rightarrow+\infty}\frac{N}{\gamma_N^2}
\log\sup_{\sigma\in \mathcal{T}}P\left(\varepsilon_{11}^N>\epsilon\right)=-\infty.
\end{equation}
By Equations \eqref{equ 5.2.13}, \eqref{equ 5.2.14}, \eqref{equ 5.2.15} and \eqref{equ 5.2.16},
\begin{equation}\label{equ 5.2.17}
\limsup_{\delta\rightarrow 0}\limsup_{N\rightarrow+\infty}\frac{N}{\gamma_N^2}
\log\sup_{\sigma\in \mathcal{T}}P\left(\sup_{0\leq t\leq \delta}\left(\eta^N_{t+\sigma}(\vec{f})-\eta^N_\sigma(\vec{f})\right)>\epsilon, (D_2^N)^c\right)=-\infty.
\end{equation}
As we have introduced,
\[
\limsup_{N\rightarrow+\infty}\frac{N}{\gamma_N^2}\log P(D_2^N)=-\infty.
\]
Consequently, Equation \eqref{equ 5.2.2} follows from Equation \eqref{equ 5.2.17} and the proof is complete.

\qed

Now we prove Equation \eqref{equ 2.10 MDP upperbound}

\proof[Proof of Equation \eqref{equ 2.10 MDP upperbound}]

According to Equations \eqref{equ 5.2.0}-\eqref{equ 5.2.0 three} and Taylor's expansion of $\exp\left(\frac{\gamma_N}{N}x\right)$ up to the second order, we have
\begin{equation}\label{equ 5.2.18}
\mathcal{Y}_{F, G, H}^N(t, \xi^N)=\exp\left(\frac{\gamma_N^2}{N}\left(\mathcal{J}_1(\eta^N, F, G, H)+\varepsilon^N_{12}+o(1)\right)\right),
\end{equation}
where $\varepsilon_{12}^N=\sum_{l=13}^{19}\varepsilon_l^N$ with
\[
\varepsilon_{13}^N=\int_0^T\eta_{s,1}^N\left(\lambda_1(-F_s+G_s)\right)\left(\int_\mathbb{T}\theta_s^I(v)\lambda_2(v)dv-\frac{1}{N}
\sum_{j=0}^{N-1}\mathbb{E}I_s^N(j)\lambda_2\left(\frac{j}{N}\right)\right)ds,
\]
\begin{align*}
\varepsilon_{14}^N=&\int_0^T\eta_{s,3}^N(\lambda_2)\Bigg(\int_\mathbb{T}\theta_s^S(v)\lambda_1(v)\left(G_s(v)-F_s(v)\right)dv\\
&-\frac{1}{N}\sum_{j=0}^{N-1}\mathbb{E}S_s^N(j)\lambda_1\left(\frac{j}{N}\right)\left(G_s\left(\frac{j}{N}\right)-F_s\left(\frac{j}{N}\right)\right)\Bigg)ds,
\end{align*}
\[
\varepsilon_{15}^N=-\int_0^T\eta_{s,1}^N(\lambda_1(G_s-F_s))\frac{1}{N}\sum_{j=0}^{N-1}\lambda_2\left(\frac{j}{N}\right)
\left(I_s^N(j)-\mathbb{E}I_s^N(j)\right)ds,
\]
\begin{align*}
\varepsilon_{16}^N=\int_0^T&\frac{1}{N\gamma_N}\sum_{i=1}^N\sum_{j=1}^N\lambda_1\left(\frac{i}{N}\right)\lambda_2\left(\frac{j}{N}\right)\times\\
&\left(\mathbb{E}\left(S_s^N(i)I_s^N(j)\right)
-\mathbb{E}\left(S_s^N(i)\right)\mathbb{E}\left(I_s^N(j)\right)\right)\left(G_s\left(\frac{i}{N}\right)-F_s\left(\frac{i}{N}\right)\right)ds,
\end{align*}
\[
\varepsilon_{17}^N=\frac{1}{2}\int_0^T\left(\mu_{s,1}\bigotimes \mu_{s,3}-\mu_{s,1}^N\bigotimes\mu_{s,3}^N\right)\left(\lambda(\cdot, \ast)(G_s(\cdot)-F_s(\cdot))^2\right)ds,
\]
\[
\varepsilon_{18}^N=\frac{1}{2}\int_0^T\left(\mu_{s,2}-\mu^N_{s,2}\right)\left(\psi(H_s-G_s)^2\right)ds
\]
and
\[
\varepsilon_{19}^N=\frac{1}{2}\int_0^T\left(\mu_{s,3}-\mu_{s,3}^N\right)\left(\phi H_s^2\right)ds.
\]
By Theorem \ref{Theorem 2.4 exponentially rapid to hydrodynamic limit} and Lemma \ref{lemma 5.2.1 exponential tightness}, for any $\epsilon>0$,
\[
\limsup_{N\rightarrow+\infty}\frac{N}{\gamma_N^2}\log P(|\varepsilon^N_l|\geq \epsilon)=-\infty
\]
for $l=13, 14, 15$. By Theorem \ref{Theorem 2.4 exponentially rapid to hydrodynamic limit},
\[
\limsup_{N\rightarrow+\infty}\frac{1}{N}\log P(|\varepsilon^N_l|\geq \epsilon)<0
\]
for $l=17, 18, 19$. By Lemma \ref{lemma 5.1.3 approximately independent},
\[
\varepsilon_{16}^N=O(\gamma_N^{-1})=o(1).
\]
In conclusion,
\begin{equation}\label{equ 5.2.19}
\limsup_{N\rightarrow+\infty}\frac{N}{\gamma_N^2}\log P(|\varepsilon^N_{12}|\geq \epsilon)=-\infty
\end{equation}
According to Equations \eqref{equ 5.2.18} and \eqref{equ 5.2.19}, the fact that Equation \eqref{equ 2.10 MDP upperbound} holds for all compact $C\subseteq \Omega_3$ follows from an analysis similar with that given in the proof of Lemma \ref{lemma 4.1.1 LDPupperboundforCompact}. Then, according to Lemma \ref{lemma 5.2.1 exponential tightness}, Equation \eqref{equ 2.10 MDP upperbound} holds for all closed $C\subseteq \Omega_3$ and the proof is complete.

\qed

\subsection{Proof of Equation \eqref{equ 2.11 MDP lowerbound}}\label{subsection 5.3}
In this subsection we prove Equation \eqref{equ 2.11 MDP lowerbound}. We first introduce some notations and definitions for later use. For any $f_1, f_2, f_3\in C^\infty(\mathbb{T})$, we denote by $\tilde{P}^N_{f_1, f_2, f_3}$ the probability measure of our SEIR process with initial condition where $\{\xi_0^N(i)\}_{0\leq i\leq N-1}$ are independent and $\xi_0^N(i)$ takes $0, 1, 2, 3$ with respective probability
\[
\rho_0\left(\frac{i}{N}\right)+\frac{\gamma_N}{N}f_1\left(\frac{i}{N}\right), \rho_1\left(\frac{i}{N}\right)+\frac{\gamma_N}{N}f_2\left(\frac{i}{N}\right),
\rho_2\left(\frac{i}{N}\right)+\frac{\gamma_N}{N}f_3\left(\frac{i}{N}\right)
\]
and
\[
1-\sum_{k=1}^3\rho_{k-1}\left(\frac{i}{N}\right)-\frac{\gamma_N}{N}\sum_{k=1}^3f_k\left(\frac{i}{N}\right)
\]
for all $0\leq i\leq N-1$. Furthermore, for $\tilde{F}, \tilde{G}, \tilde{H}\in C^{1, \infty}\left([0, T]\times\mathbb{T}\right)$, we denote by $\widetilde{P}_{f_1, f_2, f_3}^{N, \tilde{F}, \tilde{G}, \tilde{H}}$ the probability measure such that
\[
\frac{d\widetilde{P}_{f_1, f_2, f_3}^{N, \tilde{F}, \tilde{G}, \tilde{H}}}{d\tilde{P}^N_{f_1, f_2, f_3}}=\mathcal{Y}_{\tilde{F}, \tilde{G}, \tilde{H}}^N(T, \xi^N).
\]
We give two lemmas which play key roles in the proof of Equation \eqref{equ 2.11 MDP lowerbound}.

\begin{lemma}\label{lemma 5.3.1 MDP tight under transformed measure}
Under Assumptions (A) and (B), the sequence $\{\eta^N\}_{N\geq 1}$ are $\widetilde{P}_{f_1, f_2, f_3}^{N, \tilde{F}, \tilde{G}, \tilde{H}}$-tight.
\end{lemma}

\begin{lemma}\label{lemma 5.3.2 MDP LLN under transformed measure}
As $N\rightarrow+\infty$, $\eta^N$ converges in $\widetilde{P}_{f_1, f_2, f_3}^{N, \tilde{F}, \tilde{G}, \tilde{H}}$-probability to the unique solution to Equation \eqref{equ 5.1.3} with initial condition where $W_{0,k}(du)=f_k(u)du$ for $u=1,2,3$.
\end{lemma}

\proof[Proof of Lemma \ref{lemma 5.3.1 MDP tight under transformed measure}]

According to Theorem 3.1 of \cite{Mitoma1983} and Aldous' criterion, we only need to check that
\begin{equation}\label{equ 5.3.1}
\lim_{M\rightarrow+\infty}\limsup_{N\rightarrow+\infty}\widetilde{P}_{f_1, f_2, f_3}^{N, \tilde{F}, \tilde{G}, \tilde{H}}
\left(|\eta_t^N(\vec{f})|\geq M\right)=0
\end{equation}
for any $t>0, \vec{f}\in \left(C^\infty(\mathbb{T})\right)^3$ and
\begin{equation}\label{equ 5.3.2}
\lim_{\delta\rightarrow 0}\limsup_{N\rightarrow+\infty}\sup_{\sigma\in\mathcal{T}, t\leq \delta}\widetilde{P}_{f_1, f_2, f_3}^{N, \tilde{F}, \tilde{G}, \tilde{H}}\left(\left|\eta^N_{t+\sigma}(\vec{f})-\eta^N_{\sigma}(\vec{f})\right|>\epsilon\right)=0
\end{equation}
for any $\vec{f}\in \left(C^\infty(\mathbb{T})\right)^3$.

 Now we check Equations \eqref{equ 5.3.1} and \eqref{equ 5.3.2}. According to the definition of $\tilde{P}^N_{f_1, f_2, f_3}$, the fact $\ln(1+u)\leq u$ and direct calculation show that there exists $K_{13}<+\infty$ independent of $N$ such that
\[
\mathbb{E}\left(\left(\frac{d\tilde{P}^N_{f_1, f_2, f_3}}{dP}\right)^2\right)\leq \exp\left(\frac{\gamma_N^2}{N}K_{13}\right).
\]
Then, by Cauchy-Schwarz inequality,
\begin{equation}\label{equ 5.3.3}
\tilde{P}^N_{f_1, f_2, f_3}\left(A\right)\leq \exp\left(\frac{K_{13}\gamma_N^2}{2N}\right)\sqrt{P\left(A\right)}.
\end{equation}
for any event $A$. According to the definition of $\mathcal{Y}_{\tilde{F}, \tilde{G}, \tilde{H}}^N(t, \xi^N)$ and Taylor's expansion of $\exp(\frac{\gamma_N}{N}x)$ up to the second order, there exists $K_{14}<+\infty$ independent of $N$ such that
\[
\left(\frac{d\widetilde{P}_{f_1, f_2, f_3}^{N, \tilde{F}, \tilde{G}, \tilde{H}}}{d\tilde{P}^N_{f_1, f_2, f_3}}\right)^2=(\mathcal{Y}_{\tilde{F}, \tilde{G}, \tilde{H}}^N(T, \xi^N))^2\leq \mathcal{Y}_{2\tilde{F}, 2\tilde{G}, 2\tilde{H}}^N(T, \xi^N)\exp\left(\frac{\gamma_N^2}{N}\left(K_{14}+o(1)\right)\right).
\]
Since $\{\mathcal{Y}_{2\tilde{F}, 2\tilde{G}, 2\tilde{H}}^N(t, \xi^N)\}_{t\geq 0}$ is a martingale with mean $1$, by Cauchy-Schwarz inequality,
\begin{equation*}
\widetilde{P}_{f_1, f_2, f_3}^{N, \tilde{F}, \tilde{G}, \tilde{H}}\left(A\right)
\leq \exp\left(\frac{\gamma_N^2}{2N}\left(K_{14}+o(1)\right)\right)\sqrt{\tilde{P}^N_{f_1, f_2, f_3}\left(A\right)}
\end{equation*}
and hence
\begin{align}\label{equ 5.3.4}
&\widetilde{P}_{f_1, f_2, f_3}^{N, \tilde{F}, \tilde{G}, \tilde{H}}\left(A\right)\\
&\leq \exp\left(\frac{\gamma_N^2}{N}\left(\frac{K_{13}}{4}+\frac{K_{14}}{2}+o(1)\right)\right)\left(P\left(A\right)\right)^{1/4}
\notag
\end{align}
for any event $A$ according to Equation \eqref{equ 5.3.3}. As a result, Equations \eqref{equ 5.3.1} and \eqref{equ 5.3.2} follows from Equations \eqref{equ 5.2.1}, \eqref{equ 5.2.2} and \eqref{equ 5.3.4}.

Since Equations \eqref{equ 5.3.1} and \eqref{equ 5.3.2} hold, the proof is complete.

\qed

\proof[Proof of Lemma \ref{lemma 5.3.2 MDP LLN under transformed measure}]

Lemma \ref{lemma 5.3.2 MDP LLN under transformed measure} follows from an analysis similar with that in the proof of Lemma \ref{lemma 4.2.1 LLNundertransformedMeasure}. So here we only give an outline of the proof to avoid repeating similar details with those in the proof of Lemma \ref{lemma 4.2.1 LLNundertransformedMeasure}. Let $\widetilde{\eta}$ be a $\widetilde{P}_{f_1, f_2, f_3}^{N, \tilde{F}, \tilde{G}, \tilde{H}}$- weak limit of a subsequence of $\{\eta^N\}_{N\geq 1}$, the existence of which follows from Lemma \ref{lemma 5.3.1 MDP tight under transformed measure}. According to the definition of $\widetilde{P}_{f_1, f_2, f_3}^{N, \tilde{F}, \tilde{G}, \tilde{H}}$, it is easy to check that  $\widetilde{\eta}_{0,k}(du)=f_k(u)du$ according to Chebyshev's inequality, then we only need to show that $\widetilde{\eta}$ is the solution to Equation \eqref{equ 5.1.3} to complete this proof.

For any $t\geq 0, k=1,2,3$ and $f\in C^\infty(\mathbb{T})$, let
\[
\mathfrak{M}_{f,k}^N(t)=\eta_{t,k}^N(f)-\int_0^t\left(\partial_s+\mathcal{L}^n\right)\eta_{s,k}^N(f)ds,
\]
then $\{\mathfrak{M}_{f,k}^N(t)\}_{t\geq 0}$ is a martingale according to Dynkin's martingale formula. Let
\[
\Theta^N_{\tilde{F}, \tilde{G}, \tilde{H}}(t)=\Gamma_{\tilde{F}, \tilde{G}, \tilde{H}}^N(t, \xi^N)-\Gamma_{\tilde{F}, \tilde{G}, \tilde{H}}^N(0, \xi^N)
-\int_0^t\left(\partial_s+\mathcal{L}_N\right)\Gamma_{\tilde{F}, \tilde{G}, \tilde{H}}^N(s, \xi^N)ds
\]
and $\{\widetilde{\Theta}^N_{\tilde{F}, \tilde{G}, \tilde{H}}(t)\}_{t\geq 0}$ be the martingale such that
\[
d\widetilde{\Theta}^N_{\tilde{F}, \tilde{G}, \tilde{H}}(t)=\frac{1}{\Gamma_{\tilde{F}, \tilde{G}, \tilde{H}}^N(t, \xi^N)}d\Theta^N_{\tilde{F}, \tilde{G}, \tilde{H}}(t).
\]
For any $t\geq 0$, we define
\[
\widetilde{\mathfrak{M}}_{f,k}^N(t)=\mathfrak{M}_{f,k}^N(t)-\langle \mathfrak{M}_{f,k}^N , \widetilde{\Theta}^N_{\tilde{F}, \tilde{G}, \tilde{H}}\rangle_t,
\]
then $\{\widetilde{\mathfrak{M}}_{f,k}^N(t)\}_{t\geq 0}$ is a martingale under $\widetilde{P}_{f_1, f_2, f_3}^{N, \tilde{F}, \tilde{G}, \tilde{H}}$ and
\[
\left[\widetilde{\mathfrak{M}}_{f,k}^N, \widetilde{\mathfrak{M}}_{f,k}^N\right]_t=\left[\mathfrak{M}_{f,k}^N, \mathfrak{M}_{f,k}^N\right]_t
=\sum_{0\leq t\leq T}\left(\eta_{t,k}^N(f)-\eta^N_{t-,k}(f)\right)^2
\]
under $\tilde{P}^N_{f_1, f_2, f_3}$ and $\widetilde{P}_{f_1, f_2, f_3}^{N, \tilde{F}, \tilde{G}, \tilde{H}}$ according to the generalized version of Girsanov's theorem introduced in \cite{Schuppen1974}. According to an analysis similar with that leading to Equation \eqref{equ 4.2.7}, we have
\begin{equation}\label{equ 5.3.6}
P\left(\sum_{0\leq t\leq T}\left(\eta_{t,k}^N(f)-\eta^N_{t-,k}(f)\right)^2\geq\epsilon\right)
\leq e^{-\gamma_N^2\epsilon}e^{NK_2T\left(e^{\|f\|_\infty^2}-1\right)}.
\end{equation}
According to Equations \eqref{equ 5.3.4} and \eqref{equ 5.3.6}, $\sum_{0\leq t\leq T}\left(\eta_{t,k}^N(f)-\eta^N_{t-,k}(f)\right)^2$ converges in $\widetilde{P}_{f_1, f_2, f_3}^{N, \tilde{F}, \tilde{G}, \tilde{H}}$-probability to $0$. Then, according to Doob's inequality, $\sup_{0\leq t\leq T}\left|\widetilde{\mathfrak{M}}_{f,k}^N(t)\right|$ converges in $\widetilde{P}_{f_1, f_2, f_3}^{N, \tilde{F}, \tilde{G}, \tilde{H}}$-probability to $0$. According to calculation similar with that in the proof of Lemma \ref{lemma 4.2.1 LLNundertransformedMeasure},
\begin{align*}
&\langle \mathfrak{M}_{f,1}^N , \widetilde{\Theta}^N_{\tilde{F}, \tilde{G}, \tilde{H}}\rangle_t=\\
&-\int_0^t\frac{1}{N\gamma_N}\sum_{i=1}^N\sum_{j=1}^NS_t^N(i)I_t^N(j)f(i/N)\lambda(i/N, j/N)
\left(e^{\frac{\gamma_N}{N}(-\tilde{F}(i/N)+\tilde{G}(i/N))}-1\right)ds.
\end{align*}
Then, according to the fact that $\exp(x\gamma_N/N)=1+x\gamma_N/N+o(x\gamma_N/N)$, we have
\begin{align}\label{equ 5.3.7}
\eta_{t,1}^N(f)=&\eta_{0,1}^N(f)+\widetilde{\mathfrak{M}}_{f,k}^N(t)-\int_0^tK_{15,s}^N(f)ds\\
&-\int_0^t\mu_{s,1}^N\bigotimes\mu_{s,3}^N\left(\lambda(\cdot, \ast)(\tilde{G}(\cdot)-\tilde{F}(\cdot))\right)ds+o(1),\notag
\end{align}
where
\[
K_{15,s}^N(f)=\frac{1}{N\gamma_N}\sum_{i=0}^{N-1}\sum_{j=0}^{N-1}\lambda\left(\frac{i}{N}, \frac{j}{N}\right)\left(S_s^N(i)I_s^N(j)-\mathbb{E}\left(S_s^N(i)I_s^N(j)\right)\right)f\left(\frac{i}{N}\right).
\]
According to Equation \eqref{equ 5.2.0 three} and Lemma \ref{lemma 5.1.3 approximately independent},
\begin{align*}
K_{15,s}^N(f)=&\eta_{s,3}^N(\lambda_2)\frac{1}{N}\sum_{i=0}^{N-1}\mathbb{E}S_s^N(i)\lambda_1\left(\frac{i}{N}\right)f\left(\frac{i}{N}\right)\\
&+\eta_{s,1}^N(\lambda_1f)\frac{1}{N}\sum_{j=0}^{N-1}\mathbb{E}I_s^N(j)\lambda_2\left(\frac{j}{N}\right)\notag\\
&+\eta_{s,1}^N\left(\lambda_1f\right)\frac{1}{N}\sum_{j=0}^{N-1}\lambda_2\left(\frac{j}{N}\right)\left(I_s^N(j)-\mathbb{E}I_s^N(j)\right)+o(1).
\end{align*}
According to Theorem \ref{Theorem 2.4 exponentially rapid to hydrodynamic limit}, Equations \eqref{equ 5.2.1} and \eqref{equ 5.3.4},
\[
\eta_{s,1}^N\left(\lambda_1f\right)\frac{1}{N}\sum_{j=0}^{N-1}\lambda_2\left(\frac{j}{N}\right)\left(I_s^N(j)-\mathbb{E}I_s^N(j)\right)
\]
converges in $\widetilde{P}_{f_1, f_2, f_3}^{N, \tilde{F}, \tilde{G}, \tilde{H}}$-probability to $0$. According to Theorem \ref{Theorem 2.4 exponentially rapid to hydrodynamic limit} and Equation \eqref{equ 5.3.4}, $\mu_{s,1}^N\bigotimes\mu_{s,3}^N$ converges in $\widetilde{P}_{f_1, f_2, f_3}^{N, \tilde{F}, \tilde{G}, \tilde{H}}$-probability to $\mu_{s,1}\bigotimes\mu_{s,3}$.

Consequently, let $N\rightarrow+\infty$ in Equation \eqref{equ 5.3.7}, we have
\begin{align}\label{equ 5.3.8}
\widetilde{\eta}_{t,1}(f)=&\widetilde{\eta}_{0,1}(f)-\int_0^t\int_\mathbb{T}\widetilde{\eta}_{s,3}(\lambda_2)\lambda_1(v)\theta_s^S(v)f(v)dvds \\
&-\int_0^t\int_{\mathbb{T}}\widetilde{\eta}_{s,1}(\lambda_1f)\lambda_2(v)\theta_s^I(v)dvds\notag\\
&-\int_0^t\mu_{s,1}\bigotimes\mu_{s,3}\left(\lambda(\cdot, \ast)(\tilde{G}(\cdot)-\tilde{F}(\cdot))\right)ds \notag
\end{align}
for any $t\geq 0$. According to Equation \eqref{equ 5.3.8} and definitions of $\mathcal{P}_{1,s}, \mathcal{P}_{2,s}$, $\widetilde{\eta}$ satisfies the first part of Equation \eqref{equ 5.1.3}. According to analyses similar with that leading to Equation \eqref{equ 5.3.8}, $\widetilde{\eta}$ satisfies other two parts of Equation \eqref{equ 5.1.3} and hence the proof is complete.

\qed

At last, we prove Equation \eqref{equ 2.11 MDP lowerbound}

\proof[Proof of Equation \eqref{equ 2.11 MDP lowerbound}]

If
\[
\inf_{W\in O}\left(J_{ini}(W_0)+J_{dyn}(W)\right)=\infty,
\]
then Equation \eqref{equ 2.11 MDP lowerbound} is trivial. If $\inf_{W\in O}\left(J_{ini}(W_0)+J_{dyn}(W)\right)<+\infty$, then, for any $\epsilon>0$, there exists $W^\epsilon\in O$ such that
\[
J_{ini}(W\epsilon_0)+J_{dyn}(W^\epsilon)\leq \inf_{W\in O}\left(J_{ini}(W_0)+J_{dyn}(W)\right)+\epsilon.
\]
Then, according to Lemmas \ref{lemma 5.1.1 moderateRatefuntionEquivalent} and \ref{equ 5.1.2 euivexpressionofJini}, there exist $(F^\epsilon, G^\epsilon, H^\epsilon) \in \mathbb{H}$ and $h_1^\epsilon, h_2^\epsilon, h_3^\epsilon\in L^2(\mathbb{T})$ such that
\[
J_{dyn}(W^\epsilon)=\frac{1}{2}\mathcal{B}_{10}\left((F^\epsilon, G^\epsilon, H^\epsilon), (F^\epsilon, G^\epsilon, H^\epsilon)\right),
\]
\[
J_{ini}(W_0^\epsilon)=\frac{1}{2}\int_{\mathbb{T}}\sum_{k=1}^3\frac{(h^\epsilon_k(u))^2}{\rho_{k-1}(u)}+\frac{(\sum_{k=1}^3h^\epsilon_k(u))^2}
{1-\sum_{k=1}^3\rho_{k-1}(u)}du
\]
and $W^\epsilon$ is the unique solution to Equation \eqref{equ 5.1.3} with $(\tilde{F}, \tilde{G}, \tilde{H})=(F^\epsilon, G^\epsilon, H^\epsilon)$ and initial condition where $W_{0,k}(du)=h_k^\epsilon(u)du$ for $k=1,2,3$.

Since $\left(C^{1,\infty}([0, T]\times \mathbb{T})\right)^3$ is dense in $\mathbb{H}$ and $C^\infty(\mathbb{T})$ is dense in $L^2(\mathbb{T})$, there exist a sequence $\{(F^n, G^n, H^n)\}_{n\geq 1}$ in $\left(C^{1,\infty}([0, T]\times \mathbb{T})\right)^3$ and a sequence $\{(h_1^n, h_2^n, h^n_3)\}_{n\geq 1}$ in $\left(C^\infty(\mathbb{T})\right)^3$ such that $(F^n, G^n, H^n)\rightarrow (F^\epsilon, G^\epsilon, H^\epsilon)$ in $\mathbb{H}$ and $h_k^n\rightarrow h_k^\epsilon$ in $L^2(\mathbb{T})$. Let $W^n$ is the solution to Equation \eqref{equ 5.1.3} with $(\tilde{F}, \tilde{G}, \tilde{H})=(F^n, G^n, H^n)$ and initial condition where $W_{0,k}(du)=h_k^n(u)du$, then by Grownwall's inequality, $W^n\rightarrow W^\epsilon$ in $\Omega_3$. Furthermore, by Lemmas \ref{lemma 5.1.1 moderateRatefuntionEquivalent} and \ref{lemma 5.1.2},
\begin{align*}
&J_{ini}(W^n_0)+J_{dyn}(W^n)\\
&=\frac{1}{2}\mathcal{B}_{10}\left((F^n, G^n, H^n), (F^n, G^n, H^n)\right)+\frac{1}{2}\int_{\mathbb{T}}\frac{(h^n_k(u))^2}{\rho_{k-1}(u)}+\frac{(\sum_{k=1}^3h^n_k(u))^2}{1-\sum_{k=1}^3\rho_{k-1}(u)}du\\
&=\mathcal{J}_1(W^n,F^n,G^n,H^n)+\frac{1}{2}\int_{\mathbb{T}}\frac{(h^n_k(u))^2}{\rho_{k-1}(u)}+\frac{(\sum_{k=1}^3h^n_k(u))^2}{1-\sum_{k=1}^3\rho_{k-1}(u)}du\\
\end{align*}
and hence
\[
\lim_{n\rightarrow+\infty}J_{ini}(W^n_0)+J_{dyn}(W^n)=J_{ini}(W_0^\epsilon)+J_{dyn}(W^\epsilon).
\]
Consequently, there exists $m\geq 1$ such that $W^m\in O$ and
\[
J_{ini}(W^m_0)+J_{dyn}(W^m)\leq \inf_{W\in O}\left(J_{ini}(W_0)+J_{dyn}(W)\right)+2\epsilon.
\]
According to the definition of $\widetilde{P}_{h_1^m, h_2^m, h_3^m}^{N, F^m, G^m, H^m}$,
\begin{align}\label{equ 5.3.9}
&P(\eta^N\in O)\\
&=\mathbb{E}_{\widetilde{P}_{h_1^m, h_2^m, h_3^m}^{N, F^m, G^m, H^m}}\left(\left(\mathcal{Y}_{F^m, G^m, H^m}^N(t, \xi^N)\right)^{-1}\frac{dP}{d\tilde{P}^N_{h_1^m, h_2^m, h_3^m}}1_{\{\eta^N\in O\}}\right). \notag
\end{align}
Let $D_{\epsilon, 3}$ be the event that
\[
\left|\mathcal{J}_1(\eta^N, F^m, G^m, H^m)-\mathcal{J}_1(W^m, F^m, G^m, H^m)\right|\leq \epsilon.
\]
By Lemma \ref{lemma 5.3.2 MDP LLN under transformed measure}, $\eta^N$ converges in $\widetilde{P}_{h_1^m, h_2^m, h_3^m}^{N, F^m, G^m, H^m}$-probability to $W^m$ as $N\rightarrow+\infty$. Since $W^m\in O$ and $\mathcal{J}_1(\cdot, F^m, G^m, H^m)$ is continuous on $\Omega_3$, we have
\begin{equation}\label{equ 5.3.10}
\lim_{N\rightarrow+\infty}\widetilde{P}_{h_1^m, h_2^m, h_3^m}^{N, F^m, G^m, H^m}\left(D_{\epsilon, 3}\cap\{\eta^N\in O\}\right)=1.
\end{equation}
By Equations \eqref{equ 5.2.19} and \eqref{equ 5.3.4},
\begin{equation}\label{equ 5.3.11}
\lim_{N\rightarrow+\infty}\widetilde{P}_{h_1^m, h_2^m, h_3^m}^{N, F^m, G^m, H^m}\left(|\varepsilon_{12}|\leq \epsilon\right)=1.
\end{equation}
Let $D_{4,\epsilon}$ be the event that
\[
\left(\mathcal{Y}_{F^m, G^m, H^m}^N(t, \xi^N)\right)^{-1}\geq \exp\left(-\frac{\gamma_N^2}{N}(\mathcal{J}_1(W^m, F^m, G^m, H^m)+3\epsilon)\right),
\]
then by Equations \eqref{equ 5.2.18} \eqref{equ 5.3.10} and \eqref{equ 5.3.11},
\begin{equation}\label{equ 5.3.12}
\lim_{N\rightarrow+\infty}\widetilde{P}_{h_1^m, h_2^m, h_3^m}^{N, F^m, G^m, H^m}\left(D_{4,\epsilon}\cap\{\eta^N\in O\}\right)=1.
\end{equation}
According to the definition of $\tilde{P}^N_{h_1^m, h_2^m, h_3^m}$, Chebyshev's inequality and the fact that $\ln(1+u)=u-u^2/2+o(u^2)$, $\frac{N}{\gamma_N^2}\log\frac{dP}{d\tilde{P}^N_{h_1^m, h_2^m, h_3^m}}$ converges in $\widetilde{P}_{h_1^m, h_2^m, h_3^m}^{N, F^m, G^m, H^m}$-probability to
\[
-\frac{1}{2}\int_{\mathbb{T}}\sum_{k=1}^3\frac{(h_k^m(u))^2}{\rho_{k-1}(u)}+\frac{(\sum_{k=1}^3h_k^m(u))^2}{1-\sum_{k=1}^3\rho_{k-1}(u)}du,
\]
which equals $-J_{ini}(W^m_0)$ according to Lemma \ref{lemma 5.1.1 moderateRatefuntionEquivalent}. Therefore, let $D_{5, \epsilon}$ be the event that
$|\frac{N}{\gamma_N^2}\log\frac{dP}{d\tilde{P}^N_{h_1^m, h_2^m, h_3^m}}+J_{ini}(W^m_0)|\leq \epsilon$, then
\[
\lim_{N\rightarrow+\infty}\widetilde{P}_{h_1^m, h_2^m, h_3^m}^{N, F^m, G^m, H^m}(D_{5,\epsilon})=1
\]
and hence
\[
\lim_{N\rightarrow+\infty}\widetilde{P}_{h_1^m, h_2^m, h_3^m}^{N, F^m, G^m, H^m}(D_{5,\epsilon}\cap D_{4,\epsilon}\cap\{\eta^N\in O\})=1
\]
according to Equation \eqref{equ 5.3.12}. Consequently, by Equation \eqref{equ 5.3.9} and Lemma \ref{lemma 5.1.1 moderateRatefuntionEquivalent},
\begin{align*}
P(\eta^N\in O)&\geq\exp\left(-\frac{\gamma_N^2}{N}(J_{ini}(W^m_0)+\mathcal{J}_1(W^m, F^m, G^m, H^m)+4\epsilon)\right)\\
&\text{\quad\quad\quad\quad}\times P(D_{5,\epsilon}\cap D_{4,\epsilon}\cap\{\eta^N\in O\})\\
&=\exp\left(-\frac{\gamma_N^2}{N}(J_{ini}(W^m_0)+J_{dyn}(W^m)+4\epsilon)\right)(1+o(1)).
\end{align*}
As a result,
\begin{align*}
\liminf_{N\rightarrow+\infty}\frac{N}{\gamma_N^2}\log P\left(\eta^N\in O\right)&\geq -(J_{ini}(W^m_0)+J_{dyn}(W^m)+4\epsilon)\\
&\geq -\inf_{W\in O}\left(J_{ini}(W_0)+J_{dyn}(W)\right)-6\epsilon.
\end{align*}
Since $\epsilon$ is arbitrary, let $\epsilon\rightarrow 0$ and then the proof is complete.

\qed

\subsection{Proof of Theorem \ref{Theorem 2.6 MDPofHittingTimes}}\label{subsection 5.4}
In this subsection we give the proof of Theorem \ref{Theorem 2.6 MDPofHittingTimes}. For later use, we first introduce some notations and definitions. For any linear operator $\mathcal{P}$ from $C^\infty(\mathbb{T})$ to $C^\infty(\mathbb{T})$, we define $\mathcal{P}^*$ as the conjugate of $\mathcal{P}$ from $(C^\infty(\mathbb{T}))^\prime$ to $(C^\infty(\mathbb{T}))^\prime$. That is to say, for any $f\in C^\infty(\mathbb{T})$ and $\nu\in (C^\infty(\mathbb{T}))^\prime$,
\[
\mathcal{P}^*\nu(f)=\nu(\mathcal{P}f).
\]
For any $t\geq 0$, we use $\Xi_t$ to denote the $(L(C^\infty(\mathbb{T})))^{3\times 3}$-valued element
\[
\begin{pmatrix}
-\mathcal{P}_{1,t} & 0 & -\mathcal{P}_{2,t}\\
\mathcal{P}_{1,t} & -\mathcal{P}_3 & \mathcal{P}_{2,s}\\
0 & \mathcal{P}_3 & -\mathcal{P}_4
\end{pmatrix}
\]
and use $\Xi_t^*$ to denote the $(L((C^\infty(\mathbb{T}))^\prime))^{3\times 3}$-valued element
\[
\begin{pmatrix}
-\mathcal{P}^*_{1,t} & 0 & -\mathcal{P}^*_{2,t}\\
\mathcal{P}^*_{1,t} & -\mathcal{P}_3 & \mathcal{P}^*_{2,s}\\
0 & \mathcal{P}^*_3 & -\mathcal{P}^*_4
\end{pmatrix}.
\]
We denote by $\{\Phi_t\}_{t\geq 0}$ the solution to the $(L(C^\infty(\mathbb{T})))^{3\times 3}$-valued linear ordinary differential equation
\[
\begin{cases}
&\frac{d}{dt}\Phi_t=-\Xi_t\Phi_t,\\
& \Phi_0=
\begin{pmatrix}
i_d & 0 & 0\\
0 & i_d & 0\\
0 & 0 & i_d
\end{pmatrix},
\end{cases}
\]
where $i_d$ is the identity mapping on $C^\infty(\mathbb{T})$. Then we denote by $\{\Phi_t^*\}_{t\geq 0}$ the solution to the $(L((C^\infty(\mathbb{T}))^\prime))^{3\times 3}$-valued linear ordinary differential equation
\[
\begin{cases}
&\frac{d}{dt}\Phi^*_t=-\Phi_t^*\Xi_t^*,\\
& \Phi_0^*=
\begin{pmatrix}
i_d^* & 0 & 0\\
0 & i_d^* & 0\\
0 & 0 & i_d^*
\end{pmatrix},
\end{cases}
\]
where $i_d^*$ is the identity mapping on $(C^\infty(\mathbb{T}))^\prime$. Note that $\Phi_t^*$ is not derived by taking conjugate of each element of $\Phi_t$. For any $g_1, g_2, g_3\in C^\infty(\mathbb{T})$ and $s\geq 0$, we denote by $R_{s,1}^{g_1, g_2}$, $R_{s,2}^{g_1, g_2, g_3}$ and $R_{s,3}^{g_2, g_3}$ elements in $(C^\infty(\mathbb{T}))^\prime$ such that
\[
R_{s,1}^{g_1, g_2}(f)=-\int_{\mathbb{T}^2}\lambda(u,v)\theta_s^S(u)\theta^I_s(v)f(u)(g_2(u)-g_1(u))dudv,
\]
\begin{align*}
R_{s,2}^{g_1, g_2, g_3}(f)=&\int_{\mathbb{T}^2}\lambda(u,v)\theta_s^S(u)\theta_s^I(v)f(u)(g_2(u)-g_1(u))dudv\\
&-\int_{\mathbb{T}}\psi(u)\theta_s^E(u)f(u)(g_3(u)-g_2(u))du
\end{align*}
and
\begin{align*}
R_{s,3}^{g_2, g_3}(f)=\int_{\mathbb{T}}\psi(u)\theta_s^E(u)f(u)(g_3(u)-g_2(u))du+\int_\mathbb{T}\phi(u)\theta_s^I(u)f(u)g_3(u)
\end{align*}
for any $f\in C^\infty(\mathbb{T})$. For any $\nu_1, \nu_2, \nu_3\in (C^\infty(\mathbb{T}))^\prime$ and $g_1, g_2, g_3\in C^\infty(\mathbb{T})$, we use
\[
\begin{pmatrix}
\nu_1\\
\nu_2\\
\nu_3
\end{pmatrix}
\cdot
\begin{pmatrix}
g_1\\
g_2\\
g_3
\end{pmatrix}
\]
to denote $\sum_{k=1}^3\nu_k(g_k)$.

According to above definitions, for any  $\nu_1, \nu_2, \nu_3\in (C^\infty(\mathbb{T}))^\prime$ and $g_1, g_2, g_3\in C^\infty(\mathbb{T})$,
\[
\Phi_t^*\begin{pmatrix}
\nu_1\\
\nu_2\\
\nu_3
\end{pmatrix}
\cdot
\begin{pmatrix}
g_1\\
g_2\\
g_3
\end{pmatrix}
=\begin{pmatrix}
\nu_1\\
\nu_2\\
\nu_3
\end{pmatrix}
\cdot
\Phi_t\begin{pmatrix}
g_1\\
g_2\\
g_3
\end{pmatrix}.
\]
Furthermore, Equation \eqref{equ 5.1.3} can be rewritten as
\[
\frac{d}{dt}\begin{pmatrix}
W_{t,1}\\
W_{t,2}\\
W_{t,3}
\end{pmatrix}
=\Xi_t^*\begin{pmatrix}
W_{t,1}\\
W_{t,2}\\
W_{t,3}
\end{pmatrix}+\begin{pmatrix}
R^{\tilde{F}_t, \tilde{G}_t}_{t,1}\\
R^{\tilde{F}_t, \tilde{G}_t, \tilde{F}_t}_{t,2}\\
R^{\tilde{G}_t, \tilde{F}_t}_{t,3}
\end{pmatrix}
\]
and hence the solution to Equation \eqref{equ 5.1.3} is
\begin{equation}\label{equ 5.4.1 solution to (5.1.3)}
\begin{pmatrix}
W_{t,1}\\
W_{t,2}\\
W_{t,3}
\end{pmatrix}
=(\Phi_t^*)^{-1}\begin{pmatrix}
W_{0,1}\\
W_{0,2}\\
W_{0,3}
\end{pmatrix}+\int_0^t(\Phi_t^*)^{-1}\Phi_s^*\begin{pmatrix}
R^{\tilde{F}_s, \tilde{G}_s}_{s,1}\\
R^{\tilde{F}_s, \tilde{G}_s, \tilde{F}_s}_{s,2}\\
R^{\tilde{G}_s, \tilde{F}_s}_{s,3}
\end{pmatrix}ds.
\end{equation}
According to the definition of $\mathcal{B}_{10}$, for any $\tilde{F}, \tilde{G}, \tilde{H}, \hat{F}, \hat{G}, \hat{H}\in C^{1, \infty}([0, T]\times\mathbb{T})$,
\begin{equation}\label{equ 5.4.2 equivB10}
\mathcal{B}_{10}((\tilde{F}, \tilde{G}, \tilde{H}), (\hat{F}, \hat{G}, \hat{H}))
=\int_0^T\begin{pmatrix}
R^{\tilde{F}_s, \tilde{G}_s}_{s,1}\\
R^{\tilde{F}_s, \tilde{G}_s, \tilde{F}_s}_{s,2}\\
R^{\tilde{G}_s, \tilde{F}_s}_{s,3}
\end{pmatrix}
\cdot
\begin{pmatrix}
\hat{F}_s\\
\hat{G}_s\\
\hat{H}_s
\end{pmatrix}ds.
\end{equation}
For any $g_1, g_2, g_3, h_1, h_2, h_3\in L^2(\mathbb{T})$, we define
\[
\mathcal{B}_{12}((g_1, g_2, g_3), (h_1, h_2, h_3))=\int_{\mathbb{T}}\sum_{k=1}^3\frac{g_k(u)h_k(u)}{\rho_{k-1}(u)}+\frac{(\sum_{k=1}^3h_k(u))(\sum_{k=1}^3g_k(u))}{1-\sum_{k=1}^3\rho_{k-1}(u)}du.
\]
For any $\vec{g}=(g_1, g_2, g_3)\in (C^\infty(\mathbb{T}))^3$, we define $R_4\vec{g}$ as the element in $(C^\infty(\mathbb{T}))^3$ such that
\[
(R_4\vec{g})_k(u)=\rho_{k-1}(u)(g_k-\sum_{l=1}^3g_l\rho_{l-1})(u)
\]
for $u\in \mathbb{T}$ and $k=1,2,3$. Furthermore, we define $\vec{g}_{\Phi, T}$ as the element in $(C^{1,\infty}([0, T]\times \mathbb{T}))^3$ such that
\[
(\vec{g}_{\Phi, T})_s=\Phi_s(\Phi_T)^{-1}\vec{g}
\]
for any $0\leq s\leq T$. If $W$ is the solution to Equation \eqref{lemma 5.1.3 approximately independent} with initial condition where $W_{0,k}(du)=h_k(u)du$ for $k=1,2,3$, then by Lemma \ref{lemma 5.1.1 moderateRatefuntionEquivalent} and Equations \eqref{equ 5.4.1 solution to (5.1.3)}, \eqref{equ 5.4.2 equivB10}, we have
\begin{equation}\label{equ 5.4.3 equivWT}
\begin{pmatrix}
W_{T,1}\\
W_{T,2}\\
W_{T,3}
\end{pmatrix}
\cdot
\begin{pmatrix}
f_1\\
f_2\\
f_3
\end{pmatrix}=\mathcal{B}_{12}((h_1, h_2, h_3), R_4\Phi_T^{-1}\vec{f})+\mathcal{B}_{10}((\tilde{F}, \tilde{G}, \tilde{H}), \vec{f}_{\Phi, T}).
\end{equation}

Now we can give our lemma about the precise expression of $J_{contra, T}$.

\begin{lemma}\label{lemma 5.4.1 equivJcontra}
For any $x\in \mathbb{R}$ and $T>0$,
\[
J_{contra, T}(x)=\frac{x^2}{2}\frac{1}{\mathcal{B}_{12}(R_4\Phi_T^{-1}\vec{f}, R_4\Phi_T^{-1}\vec{f})+\mathcal{B}_{10}(\vec{f}_{\Phi, T}, \vec{f}_{\Phi, T})}.
\]
Furthermore,
\[
J_{contra, T}(x)=J{ini}(W^{T,x}_0)+J_{dyn}(W^{T,x}),
\]
where $W^{T,x}$ is the solution to Equation \eqref{equ 5.1.3} in time interval $[0, T]$ with
\[
(\tilde{F}, \tilde{G}, \tilde{H})=\frac{x\vec{f}_{\Phi, T}}{\mathcal{B}_{12}(R_4\Phi_T^{-1}\vec{f}, R_4\Phi_T^{-1}\vec{f})+\mathcal{B}_{10}(\vec{f}_{\Phi, T}, \vec{f}_{\Phi, T})}
\]
and initial condition where
\[
W^{T,x}_{0,k}(du)=\frac{x(R_4\Phi_T^{-1}\vec{f})_k(u)}{\mathcal{B}_{12}(R_4\Phi_T^{-1}\vec{f}, R_4\Phi_T^{-1}\vec{f})+\mathcal{B}_{10}(\vec{f}_{\Phi, T}, \vec{f}_{\Phi, T})}du
\]
for $k=1,2,3$.

\end{lemma}

\proof[Proof of Lemma \ref{lemma 5.4.1 equivJcontra}]

If $W$ is the solution to Equation \eqref{equ 5.1.3} with initial density $\vec{h}=(h_1, h_2, h_3)$, then
\[
J_{ini}(W_0)+J_{dyn}(W)=\frac{1}{2}\left(\mathcal{B}_{12}(\vec{h}, \vec{h})+\mathcal{B}_{10}((\tilde{F}, \tilde{G}, \tilde{H}), (\tilde{F}, \tilde{G}, \tilde{H}))\right)
\]
according to Lemmas \ref{lemma 5.1.1 moderateRatefuntionEquivalent} and \ref{lemma 5.1.2}. Then, by Equation \eqref{equ 5.4.3 equivWT},
\begin{align}\label{equ 5.4.4 equivJcontra}
J_{contra, T}(x)
&=\inf\Bigg\{\frac{1}{2}\left(\mathcal{B}_{12}(\vec{h}, \vec{h})+\mathcal{B}_{10}((\tilde{F}, \tilde{G}, \tilde{H}), (\tilde{F}, \tilde{G}, \tilde{H}))\right): \vec{h}\in (C^\infty(\mathbb{T}))^3, \notag\\
&~(\tilde{F}, \tilde{G}, \tilde{H})\in \mathbb{H}, \mathcal{B}_{12}(\vec{h}, R_4\Phi_T^{-1}\vec{f})+\mathcal{B}_{10}((\tilde{F}, \tilde{G}, \tilde{H}), \vec{f}_{\Phi, T})=x\Bigg\}.
\end{align}
According to calculus of variation, under the constraint condition
\[
\mathcal{B}_{12}(\vec{h}, R_4\Phi_T^{-1}\vec{f})+\mathcal{B}_{10}((\tilde{F}, \tilde{G}, \tilde{H}), \vec{f}_{\Phi, T})=x,
\]
the objective function $\frac{1}{2}\left(\mathcal{B}_{12}(\vec{h}, \vec{h})+\mathcal{B}_{10}((\tilde{F}, \tilde{G}, \tilde{H}), (\tilde{F}, \tilde{G}, \tilde{H}))\right)$ gets minimum at some $(\vec{h}, (\tilde{F}, \tilde{G}, \tilde{H}))$ such that
\[
(\vec{h}, (\tilde{F}, \tilde{G}, \tilde{H}))=r(R_4\Phi_T^{-1}\vec{f}, \vec{f}_{\Phi, T})
\]
for some $r\in \mathbb{R}$. According to the constraint condition,
\[
r=\frac{x}{\mathcal{B}_{12}(R_4\Phi_T^{-1}\vec{f}, R_4\Phi_T^{-1}\vec{f})+\mathcal{B}_{10}(\vec{f}_{\Phi, T}, \vec{f}_{\Phi, T})}
\]
and the proof is complete.

\qed

Moderate deviations of hitting times of a family of density-dependent Markov chains are given in \cite{He2023}. Lemma \ref{lemma 5.4.1 equivJcontra} is an analogue of Proposition 2.2 of \cite{He2023}. With Lemma \ref{lemma 5.4.1 equivJcontra}, the proof of Theorem \ref{Theorem 2.6 MDPofHittingTimes} follows from an analysis similar with that in the proof of the main theorem of \cite{He2023}. So here we only give an outline of the proof of Theorem \ref{Theorem 2.6 MDPofHittingTimes}.

\proof[Proof of Theorem \ref{Theorem 2.6 MDPofHittingTimes}]

We write $\Omega_3$ as $\Omega_{3,T}$ when we need to distinguish different $T$. For sufficiently small $\delta>0$, it is easy to check that there exist closed set $C_1\subseteq \Omega_{3, \tau_c+\delta}$ and closed set $C_2\subseteq \Omega_{3,\tau_c-\delta}$ such that $\{|\tau_c^N-\tau_c|\geq \delta\}\subseteq C_1\cup C_2$ and $C_1\not\ni\{\mu_t\}_{0\leq t\leq \tau_c+\delta}, C_2\not\ni\{\mu_t\}_{0\leq t\leq \tau_c-\delta}$. Therefore, by Lemmas \ref{lemma 2.1 good LDP rate} and \ref{lemma 4.3.1},
\begin{equation*}
\limsup_{N\rightarrow+\infty}\frac{1}{N}\log P(|\tau_c^N-\tau_c|\geq \delta)<0
\end{equation*}
and hence
\begin{equation}\label{equ 5.4.5 LDPupperboundofHittingTimes}
\limsup_{N\rightarrow+\infty}\frac{N}{\gamma_N^2}\log P(|\tau_c^N-\tau_c|\geq \delta)=-\infty.
\end{equation}
According to Lagrange's mean value theorem and the fact that
\begin{align}\label{equ 5.4.6 approximateequal}
-\frac{N}{\gamma_N}&\left(\sum_{k=1}^3\mu^N_{\tau_c^N,k}(f_k)-\sum_{k=1}^3\mu_{\tau_c^N,k}(f_k)\right) \\
&=\frac{N}{\gamma_N}\left(\sum_{k=1}^3\mu_{\tau_c^N, k}(f_k)-\sum_{k=1}^3\mu_{\tau_c,k}(f_k)\right)+O(\gamma_N^{-1})=\sum_{k=1}^3\eta^N_{\tau_c^N,k}(f_k)+o(1), \notag
\end{align}
for any $x>0$ and sufficiently small $\epsilon>0$, there exists $\delta=\delta(\epsilon)>0$ such that
\begin{align*}
&\{|\tau_c^N-\tau_c|\leq \delta\}\cap \{\frac{N}{\gamma_N}(\tau_c^N-\tau_c)>x\}\\
&\subseteq \left\{\sum_{k=1}^3\eta^N_{\tau_c^N,k}(f_k)<-\zeta_1(\epsilon, x)\right\}\cap\{|\tau_c^N-\tau_c|\leq \delta\}
\end{align*}
for sufficiently large $N$, where
\[
\zeta_1(\epsilon, x)=\left(\frac{d}{dt}\sum_{k=1}^3\mu_{t,k}(f_k)\Bigg|_{t=\tau_c}-\epsilon\right)x-\epsilon.
\]
Then, by Equations \eqref{equ 5.4.5 LDPupperboundofHittingTimes} and \eqref{equ 2.10 MDP upperbound},
\begin{align}\label{equ 5.4.7}
\limsup_{N\rightarrow+\infty}\frac{N}{\gamma_N^2}\log P(\frac{N}{\gamma_N}(\tau_c^N-\tau_c)>x)
&\leq \limsup_{N\rightarrow+\infty}\frac{N}{\gamma_N^2}\log P(\eta^N\in D_6(x,\delta,\epsilon))\\
&\leq -\inf_{W\in\overline{D}_6(x,\delta, \epsilon)}(J_{ini}(W_0)+J_{dyn, \tau_c+\delta}(W)), \notag
\end{align}
where
\[
D_6(x,\delta,\epsilon)=\{W:~\inf_{\tau_c-\delta\leq s\leq \tau_c+\delta}\sum_{k=1}^3W_{s,k}(f_k)<-\zeta_1(\epsilon, x)\}.
\]
For any $W\in \overline{D}_6(x,\delta, \epsilon)$, we have $\inf_{\tau_c-2\delta\leq s\leq \tau_c+\delta}\sum_{k=1}^3W_{s,k}(f_k)\leq -\zeta_1(\epsilon, x)$. Hence there exists $s_0\in [\tau_c-2\delta\leq s\leq \tau_c+\delta]$ such that
\[
\sum_{k=1}^3W_{s_0,k}(f_k)=a\leq -\zeta_1(\epsilon, x).
\]
Then, by Lemma \ref{lemma 5.4.1 equivJcontra},
\begin{align*}
&J_{ini}(W_0)+J_{dyn, \tau_c+\delta}(W)\geq J_{ini}(W_0)+J_{dyn, s_0}(W)\\
&\geq J_{contra, s_0}(a)\\
&=\frac{a^2}{2}\frac{1}{\mathcal{B}_{12}(R_4\Phi_{s_0}^{-1}\vec{f}, R_4\Phi_{s_0}^{-1}\vec{f})+\mathcal{B}_{10}(\vec{f}_{\Phi, s_0}, \vec{f}_{\Phi, s_0})}\\
&\geq \frac{\zeta_1(\epsilon, x)^2}{2}\frac{1}{\sup_{\tau_c-2\delta\leq s\leq \tau_c+\delta}\mathcal{B}_{12}(R_4\Phi_{s}^{-1}\vec{f}, R_4\Phi_{s}^{-1}\vec{f})+\mathcal{B}_{10}(\vec{f}_{\Phi, s}, \vec{f}_{\Phi, s})}.
\end{align*}
Therefore,
\begin{align*}
&\limsup_{N\rightarrow+\infty}\frac{N}{\gamma_N^2}\log P(\frac{N}{\gamma_N}(\tau_c^N-\tau_c)>x)\\
&\leq -\frac{\zeta_1(\epsilon, x)^2}{2}\frac{1}{\sup_{\tau_c-2\delta\leq s\leq \tau_c+\delta}\mathcal{B}_{12}(R_4\Phi_{s}^{-1}\vec{f}, R_4\Phi_{s}^{-1}\vec{f})+\mathcal{B}_{10}(\vec{f}_{\Phi, s}, \vec{f}_{\Phi, s})}.
\end{align*}
Since $\epsilon$ is arbitrary, let $\epsilon\rightarrow 0$ and then
\begin{align}\label{equ 5.4.8 hittingMDPgeqUpperbound}
\limsup_{N\rightarrow+\infty}\frac{N}{\gamma_N^2}\log P(\frac{N}{\gamma_N}(\tau_c^N-\tau_c)>x)
&\leq -x^2J_{contra, \tau_c}(1)\left(\frac{d}{dt}\sum_{k=1}^3\mu_{t,k}(f_k)\Bigg|_{t=\tau_c}\right)^2 \notag\\
&=-J_{hit}(x).
\end{align}

According to an analysis similar with that leading to Equation \eqref{equ 5.4.7},
\begin{align}\label{equ 5.4.9}
\liminf_{N\rightarrow+\infty}\frac{N}{\gamma_N^2}\log P(\frac{N}{\gamma_N}(\tau_c^N-\tau_c)>x)
&\geq \liminf_{N\rightarrow+\infty}\frac{N}{\gamma_N^2}\log P(\eta^N\in D_7(x,\delta,\epsilon))\\
&\geq -\inf_{W\in D^o_7(x,\delta, \epsilon)}(J_{ini}(W_0)+J_{dyn, \tau_c+\delta}(W)), \notag
\end{align}
where
\[
D_7(x,\delta,\epsilon)=\{W:~\sup_{\tau_c-\delta\leq s\leq \tau_c+\delta}\sum_{k=1}^3W_{s,k}(f_k)<-\zeta_2(\epsilon, x)\}
\]
with $\zeta_2(\epsilon, x)=\left(\frac{d}{dt}\sum_{k=1}^3\mu_{t,k}(f_k)\Bigg|_{t=\tau_c}+\epsilon\right)x+\epsilon$.
According to the definition of $W^{T,x}$, $W^{T,x}=xW^{T,1}$.
Let
\[
b=-\frac{\zeta_2(\epsilon, x)+\epsilon}{\inf_{\tau-\delta\leq s\leq \tau_c+\delta}\sum_{k=1}^3W_{s,k}^{\tau_c+\delta,1}(f_k)},
\]
then $W^{\tau_c+\delta, b}\in D^o_7(x,\delta, \epsilon)$ and hence
\begin{align*}
\liminf_{N\rightarrow+\infty}\frac{N}{\gamma_N^2}\log P(\frac{N}{\gamma_N}(\tau_c^N-\tau_c)>x)&\geq -(J_{ini}(W^{\tau_c+\delta, b}_0)+J_{dyn, \tau_c+\delta}(W^{\tau_c+\delta, b}))\\
&=-b^2J_{contra, \tau_c+\delta}(1)
\end{align*}
by Equation \eqref{equ 5.4.9}. Since $\epsilon$ is arbitrary, let $\epsilon\rightarrow 0$ and then
\begin{align}\label{equ 5.4.10 hittingMDPgeqLowerbound}
\liminf_{N\rightarrow+\infty}\frac{N}{\gamma_N^2}\log P(\frac{N}{\gamma_N}(\tau_c^N-\tau_c)>x)
& \geq -\left(\frac{d}{dt}\sum_{k=1}^3\mu_{t,k}(f_k)\Bigg|_{t=\tau_c}\right)^2x^2J_{contra, \tau_c}(1)\notag\\
&=-J_{hit}(x).
\end{align}
By Equations \eqref{equ 5.4.8 hittingMDPgeqUpperbound} and \eqref{equ 5.4.10 hittingMDPgeqLowerbound},
\begin{equation}\label{equ 5.4.11 hittingMDPgeq}
\lim_{N\rightarrow+\infty}\frac{N}{\gamma_N^2}\log P(\frac{N}{\gamma_N}(\tau_c^N-\tau_c)>x)=-J_{hit}(x)
\end{equation}
for any $x>0$. According to an analysis similar with that leading to Equation \eqref{equ 5.4.11 hittingMDPgeq},
\begin{equation}\label{equ 5.4.12 hittingMDPleq}
\lim_{N\rightarrow+\infty}\frac{N}{\gamma_N^2}\log P(\frac{N}{\gamma_N}(\tau_c^N-\tau_c)<-x)=-J_{hit}(x)
\end{equation}
for any $x>0$. Theorem \ref{Theorem 2.6 MDPofHittingTimes} follows from Equations \eqref{equ 5.4.11 hittingMDPgeq} and \eqref{equ 5.4.12 hittingMDPleq}.

\qed

\quad

\textbf{Declaration of competing interest.}
The authors declare that they have no known competing financial interests or personal
relationships that could have appeared to influence the work reported in this paper.

\quad

\textbf{Data Availability.} Data sharing not applicable to this article as no datasets were generated or analysed during the current study.

\quad

\textbf{Acknowledgments.} The authors are grateful to the financial
support from the Fundamental Research Funds for the Central Universities with grant number 2022JBMC039.

{}
\end{document}